\documentclass[12pt]{amsart}

\usepackage[margin=1in]{geometry}

\usepackage{amsmath}
\usepackage{comment}
\usepackage{listings}
\usepackage{graphicx}
\usepackage{bbm}
\usepackage[export]{adjustbox}
\usepackage{xy}
\usepackage{tikz}
\usetikzlibrary{shapes}
\usetikzlibrary{backgrounds}
\usetikzlibrary{trees}
\usetikzlibrary{intersections}
\usetikzlibrary{knots}
\usepackage{tikz-qtree}
\usepackage{youngtab}

\newcommand{\GPA}{\text{GPA}}
\newcommand{\Rep}{\text{Rep}}

\newcommand{\overcrossing}{\begin{tikzpicture}[baseline]
\begin{knot}[clip width = 4]
	\strand (-.5,.5) -- (.5,-.5);
	\strand (-.5,-.5) -- (.5,.5);
\end{knot}
\end{tikzpicture}
}

\newcommand{\undercrossing}{\begin{tikzpicture}[baseline]
\begin{knot}[clip width = 4]
	\strand (-.5,-.5) -- (.5,.5);
	\strand (-.5,.5) -- (.5,-.5);
\end{knot}
\end{tikzpicture}}
\newcommand{\identity}{\begin{tikzpicture}[baseline, rounded corners = 5mm]
    \draw (.5,.5) -- (.1,0) -- (.5,-.5);
    \draw (-.5,.5) -- (-.1,0) -- (-.5,-.5);
\end{tikzpicture}}
\newcommand{\cupcap}{\begin{tikzpicture}[baseline, rounded corners = 5mm]
    \draw (-.5,-.5) -- (0,.1) -- (.5,-.5);
    \draw (-.5,.5) -- (0,-.1) -- (.5,.5);
\end{tikzpicture}}
\newcommand{\virtual}{\begin{tikzpicture}[baseline]
    \draw (-.5,-.5)-- (.5,.5);
    \draw (.5,-.5) -- (-.5,.5);
    \draw[fill=white] (0,0) circle (.1cm);
\end{tikzpicture}
}

\newcommand{\skeincircle}{\begin{tikzpicture}[baseline]
\draw circle (.5cm);
\end{tikzpicture}}

\newcommand{\stI}{ \begin{tikzpicture}[baseline=0cm]
 	\draw (0,.2) .. controls +(30:.3cm) .. (45:.8cm);
 	\draw (0,.2) .. controls +(150:.3cm) .. (135:.8cm);
	\draw (0,.2) -- (0,-.2);
 	\draw (0,-.2) .. controls +(-30:.3cm) .. (-45:.8cm);
 	\draw (0,-.2) .. controls +(-150:.3cm) .. (-135:.8cm);
\end{tikzpicture}
}

\newcommand{\stH}{ \begin{tikzpicture}[baseline=0cm,rotate=90]
 	\draw (0,.2) .. controls +(30:.3cm) .. (45:.8cm);
 	\draw (0,.2) .. controls +(150:.3cm) .. (135:.8cm);
	\draw (0,.2) -- (0,-.2);
 	\draw (0,-.2) .. controls +(-30:.3cm) .. (-45:.8cm);
 	\draw (0,-.2) .. controls +(-150:.3cm) .. (-135:.8cm);
\end{tikzpicture}}
\newcommand{\overcrossinga}{\begin{tikzpicture}[baseline]
\begin{knot}[clip width = 4]
	\strand[<-] (-.5,.5) -- (.5,-.5);
	\strand[->] (-.5,-.5) -- (.5,.5);
\end{knot}
\end{tikzpicture}}
\newcommand{\overcrossingb}{\begin{tikzpicture}[baseline]
\begin{knot}[clip width = 4]
	\strand[->] (-.5,.5) -- (.5,-.5);
	\strand[->] (-.5,-.5) -- (.5,.5);
\end{knot}
\end{tikzpicture}}
\newcommand{\overcrossingc}{\begin{tikzpicture}[baseline]
\begin{knot}[clip width = 4]
	\strand[->] (-.5,.5) -- (.5,-.5);
	\strand[<-] (-.5,-.5) -- (.5,.5);
\end{knot}
\end{tikzpicture}}
\newcommand{\overcrossingd}{\begin{tikzpicture}[baseline]
\begin{knot}[clip width = 4]
	\strand[<-] (-.5,.5) -- (.5,-.5);
	\strand[<-] (-.5,-.5) -- (.5,.5);
\end{knot}
\end{tikzpicture}}

\newcommand{\undercrossinga}{\begin{tikzpicture}[baseline]
\begin{knot}[clip width = 4]
	\strand[->] (-.5,-.5) -- (.5,.5);
	\strand[<-] (-.5,.5) -- (.5,-.5);
\end{knot}
\end{tikzpicture}}

\newcommand{\identitya}{\begin{tikzpicture}[baseline, rounded corners = 5mm]
\begin{knot}
    \strand[<-] (.5,.5) -- (.1,0) -- (.5,-.5);
    \strand[<-] (-.5,.5) -- (-.1,0) -- (-.5,-.5);
\end{knot}
\end{tikzpicture}}
\newcommand{\identityb}{\begin{tikzpicture}[baseline, rounded corners = 5mm]
\begin{knot}
    \strand[->] (.5,.5) -- (.1,0) -- (.5,-.5);
    \strand[<-] (-.5,.5) -- (-.1,0) -- (-.5,-.5);
\end{knot}
\end{tikzpicture}}
\newcommand{\identityc}{\begin{tikzpicture}[baseline, rounded corners = 5mm]
\begin{knot}
    \strand[<-] (.5,.5) -- (.1,0) -- (.5,-.5);
    \strand[->] (-.5,.5) -- (-.1,0) -- (-.5,-.5);
\end{knot}
\end{tikzpicture}}
\newcommand{\identityd}{\begin{tikzpicture}[baseline, rounded corners = 5mm]
\begin{knot}
    \strand[->] (.5,.5) -- (.1,0) -- (.5,-.5);
    \strand[->] (-.5,.5) -- (-.1,0) -- (-.5,-.5);
\end{knot}
\end{tikzpicture}}

\newcommand{\cupcapa}{\begin{tikzpicture}[baseline,rounded corners = 5mm]
\begin{knot}
	\strand[->] (-.5,-.5) -- (0,.1) -- (.5,-.5);
	\strand[->] (-.5,.5) -- (0,-.1) -- (.5,.5);
\end{knot}
\end{tikzpicture}}

\newcommand{\cupcapd}{\begin{tikzpicture}[baseline,rounded corners = 5mm]
\begin{knot}
	\strand[<-] (-.5,-.5) -- (0,.1) -- (.5,-.5);
	\strand[<-] (-.5,.5) -- (0,-.1) -- (.5,.5);
\end{knot}
\end{tikzpicture}}

\newcommand{\virtuala}{\begin{tikzpicture}[baseline]
    \draw[->] (-.5,-.5)-- (.5,.5);
    \draw[->] (.5,-.5) -- (-.5,.5);
    \draw[fill=white] (0,0) circle (.1cm);
\end{tikzpicture}}

\newcommand{\xdiag[1]}{\begin{tikzpicture}[baseline = -.05cm]
\draw[fill = black!40!, draw = none] (-.35,-.5) rectangle (.35,.5);
\draw (-.35,-.5) -- (-.35,.5);
\draw (.35,-.5) -- (.35,.5);
\draw[fill = white] (-.5,-.3) rectangle (.5,.3);
\node at (0,0) {\large{$#1$}};
\node at (-.7,0) {$a$};
\node at (.65,0) {$b$};
\end{tikzpicture}}

\usepackage{amssymb}
\usepackage{epstopdf}
\usepackage{amsthm}
\usepackage{enumitem}
\newtheorem{theorem}{Theorem}[subsection]
\newtheorem{customthm}{Theorem}
\newtheorem{corollary}[theorem]{Corollary}
\newtheorem{lemma}[theorem]{Lemma}
\newtheorem{proposition}[theorem]{Proposition}

\theoremstyle{definition}
\newtheorem{definition}[theorem]{Definition}
\newtheorem{example}[theorem]{Example}

\theoremstyle{remark}
\newtheorem*{remark}{Remark}

\lstdefinestyle{mystyle}{
basicstyle=\times}
\title{Skein Theories for Virtual Tangles}
\author{Joshua R. Edge}
\address{Department of Mathematics and Computer Science, Denison University, 100 W. College St., Granville, OH 43023}
\email{edgej@denison.edu}
\date{August, 10, 2020}

\begin{document}

\begin{abstract}
In this paper, we use skein-theoretic techniques to classify all virtual knot polynomials and trivalent graph invariants with certain smallness conditions. The first half of the paper classifies all virtual knot polynomials giving non-trivial invariants strictly smaller than the one given by the Higman-Sims spin model. In particular, we exhibit a family of skein theories coming from $\Rep(O(2))$ with an interesting braiding. In addition, all skein theories of oriented virtual tangles with some smallness conditions are classified. In the second half of the paper, we classify all non-trivial invariants of (perhaps non-planar) trivalent graphs coming from symmetric trivalent categories. For each of these categories, we also classify when the sub-category generated by only the trivalent vertex is braided. An interesting example of this arise from the tensor category Deligne's $S_t$. 
\end{abstract}

\maketitle

\section{Introduction}
Skein relations are tools that are used to study knots. In the late 1960s, Conway first demonstrated how the Alexander polynomial could be computed using skein relations. The Alexander polynomial was the only known knot polynomial for almost 60 years, until Jones exhibited a new knot polynomial using skein relations \cite{Jon85}. This led to the discovery of other knot polynomials using skein relations including the HOMFLY polynomial \cite{FYH85} and the Kauffman polynomial \cite{Kau90}. In \cite{Kau87}, Kauffman described how the Jones polynomial could be computed from spin models. In \cite{Jae95}, Jaeger gave a whole class of spin models for the Kauffman polynomial, including the Higman-Sims spin model.

Later in the 1990s, Kauffman introduced the idea of virtual knot theory \cite{Kau99}, which allows strands in a knot or link to cross in one of three ways:
\begin{center}
    \begin{tabular}{ccccccc}
         \begin{tikzpicture}[baseline]
\begin{knot}[clip width = 4]
	\strand (-.5,.5) -- (.5,-.5);
	\strand (-.5,-.5) -- (.5,.5);
\end{knot}
\end{tikzpicture}
 &&& \begin{tikzpicture}[baseline]
\begin{knot}[clip width = 4]
	\strand (-.5,-.5) -- (.5,.5);
	\strand (-.5,.5) -- (.5,-.5);
\end{knot}
\end{tikzpicture} &&& \begin{tikzpicture}[baseline]
    \draw (-.5,-.5)-- (.5,.5);
    \draw (.5,-.5) -- (-.5,.5);
    \draw[fill=white] (0,0) circle (.1cm);
\end{tikzpicture}
 \\
         overcrossing &&& undercrossing &&& virtual crossing
    \end{tabular}
\end{center}
The first two crossings are the standard crossings from knot theory while the latter crossing is a symmetric crossing where the strands are passing ``through" each other. In \cite{Kau99}, Kauffman gave a virtual version of the Jones polynomial using skein relations. Rather mysteriously, although the skein theory of the Jones polynomial has a skein-theoretic virtual version most other knot polynomials do not in general. The question of when knot polynomials have skein-theoretic virtual counterparts is one motivation of this paper.

When thought about in the correct way, \cite{Jae95} gives skein theoretic invariants of virtual knots and links, although some of these invariants---including the one from the Higman-Sims spin model---are complicated and their skein theory is difficult to write down. In this paper, we classify all virtual skein theories \emph{strictly smaller} than the one given by the Higman-Sims spin model. Skein theoretic classifications of this type can first be seen in \cite{BJ00}, \cite{BJ03} and also in \cite{Kup94} among others. The result of this classification can be summarized as follows:

\begin{customthm}
A non-trivial virtual knot polynomial where the dimension of the 4-box space is less than 4 is equivalent to one of the following: 

\begin{enumerate}[label = \roman{enumi}.]
\item The virtual Jones polynomial
\item The Kauffman polynomial at $d=-2$ and $a=\pm 1$ equipped with a symmetric crossing
\begin{equation*}
\begin{tikzpicture}[baseline]
    \draw (-.5,-.5)-- (.5,.5);
    \draw (.5,-.5) -- (-.5,.5);
    \draw[fill=white] (0,0) circle (.1cm);
\end{tikzpicture}
 ~=~ -\left( ~ \begin{tikzpicture}[baseline, rounded corners = 5mm]
    \draw (.5,.5) -- (.1,0) -- (.5,-.5);
    \draw (-.5,.5) -- (-.1,0) -- (-.5,-.5);
\end{tikzpicture} ~+~ \begin{tikzpicture}[baseline, rounded corners = 5mm]
    \draw (-.5,-.5) -- (0,.1) -- (.5,-.5);
    \draw (-.5,.5) -- (0,-.1) -- (.5,.5);
\end{tikzpicture} ~\right)
\end{equation*}
\item A one-parameter family of knot polynomials where the crossing has formula
\begin{equation*}
\begin{tikzpicture}[baseline]
\begin{knot}[clip width = 4]
	\strand (-.5,.5) -- (.5,-.5);
	\strand (-.5,-.5) -- (.5,.5);
\end{knot}
\end{tikzpicture}
 ~=~ \dfrac{a^{-1} - a}{2} ~ \begin{tikzpicture}[baseline, rounded corners = 5mm]
    \draw (.5,.5) -- (.1,0) -- (.5,-.5);
    \draw (-.5,.5) -- (-.1,0) -- (-.5,-.5);
\end{tikzpicture} ~-~ \dfrac{a^{-1}-a}{2} ~ \begin{tikzpicture}[baseline, rounded corners = 5mm]
    \draw (-.5,-.5) -- (0,.1) -- (.5,-.5);
    \draw (-.5,.5) -- (0,-.1) -- (.5,.5);
\end{tikzpicture} ~+~ \dfrac{a^{-1}+a}{2} ~ \begin{tikzpicture}[baseline]
    \draw (-.5,-.5)-- (.5,.5);
    \draw (.5,-.5) -- (-.5,.5);
    \draw[fill=white] (0,0) circle (.1cm);
\end{tikzpicture}

\end{equation*}
\end{enumerate}
\end{customthm}

If we now assume that all of the strands are oriented, we obtain the following classification:

\begin{customthm}
If $\mathcal{V}$ is a spherical quotient of the planar algebra of oriented virtual tangles such that the vector space of diagrams with endpoints oriented $(+,+,-,-)$ is two-dimensional, then $\mathcal{V}$ is isomorphic to either Deligne's $GL_t$ or $\Rep(GL(1))$. In both cases the planar algebra is equipped with the following braiding:
\begin{equation*}
    \overcrossinga ~=~ a\cdot \identitya
\end{equation*}
\end{customthm}

In addition to the contributions skein-theoretic calculations have in knot theory, they also provide invariants in graph theory. The planar algebras found in \cite{MPS15}, called trivalent categories, give important invariants of planar trivalent graphs. One can think of the symmetric version of these categories as simply removing the planarity condition on the graphs. The second part of this paper classifies all skein-theoretic invariants of trivalent graphs with some smallness conditions. These graph-theoretic invariants are summarized in Theorem \ref{thm:class}. In the last portion of the paper, we prove the following result:
\begin{customthm}
The non-degenerate quotient of Deligne's $S_t$ has a non-symmetric braiding on the subcategory generated by the trivalent vertex when $t\in\mathbb{C}-\{0,1,2,3\}$ and is given by the formula
\begin{equation*}
\begin{tikzpicture}[baseline]
\begin{knot}[clip width = 4]
	\strand (-.5,.5) -- (.5,-.5);
	\strand (-.5,-.5) -- (.5,.5);
\end{knot}
\end{tikzpicture}
 ~=~ (q^2-1)\cdot\begin{tikzpicture}[baseline, rounded corners = 5mm]
    \draw (.5,.5) -- (.1,0) -- (.5,-.5);
    \draw (-.5,.5) -- (-.1,0) -- (-.5,-.5);
\end{tikzpicture} ~+~ q^{-2}\begin{tikzpicture}[baseline, rounded corners = 5mm]
    \draw (-.5,-.5) -- (0,.1) -- (.5,-.5);
    \draw (-.5,.5) -- (0,-.1) -- (.5,.5);
\end{tikzpicture} ~-~(q^2+q^{-2})\stI
\end{equation*}
at $t = q^2 + 2 +q^{-2}.$ 
\end{customthm}

Those familiar with planar algebras will note that this is the formula for the braiding of $(SO(3))_q.$ This connection is interesting and will be discussed in greater detail in this paper.

\subsection{Source code}
When Mathematica is used to aid in computation, the source code for that calculation is included in the ArXiv source code. The diagrams used in this paper are included under the ``diagrams" subdirectory. 

\subsection{Acknowledgments}
The author would like to thank his advisor, Noah Snyder, for his guidance through this project. He would also like to thank Pavel Etingof for explaining the origin of the subbraiding on $\Rep(O(2))$ and Dylan Thurston for his suggestion of the $S_6$ action. In addition, he would like to acknowledge the authors of \cite{MPS15}--Scott Morrison, Emily Peters, and Noah Snyder--for including the source code for the diagrams in their papers. Lastly, he would like to thank Maxime Scott for his helpful suggestions. The author was also supported through NSF grant DMS-1454767.

\section{Planar algebra definitions}\label{background}
In this section, we record the formal definitions of a planar algebra similar to \cite{Jon99} and \cite{Pet09}. In addition, we also establish some conventions that we will use in subsequent sections. 

\subsection{Definition of a planar algebra}

\begin{definition}
A planar diagram is described as follows. The diagram is the unit disc $\mathcal{D}$ with a finite number of clockwise-numbered boundary points (possibly zero). In addition, there are a finite number of internal discs each with their own finite number of clockwise-numbered boundary points (possibly zero). The boundary points of all diagrams are paired by smooth, non-crossing curves called the strings of the diagram. In addition, there may also be a finite number of closed strings not connecting to any discs. The strings lie between the internal discs and $\mathcal{D}$, and the connected components of the interior minus the strings are called the regions of the planar diagram. To indicate the first point on the boundary of the external or an internal disc, the region immediately preceding it is starred near the relevant boundary component. 
\end{definition}
An example of a planar diagram is
\begin{equation}\label{eq:paex}
\includegraphics[valign = c, scale = .7]{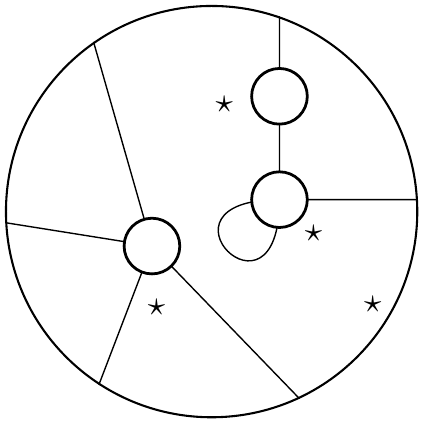}
\end{equation}
Because of the appearance of these diagrams, they are often called ``spaghetti and meatball" diagrams. We can compose diagrams $A$ and $B$ if the number of strands connecting to the outer disc of $B$ is equal to the number of strands connecting to one of the inner discs of $A$ or vice versa. An example of this composition is below:
\begin{equation}\label{comp}
\includegraphics[valign = c, scale = .7]{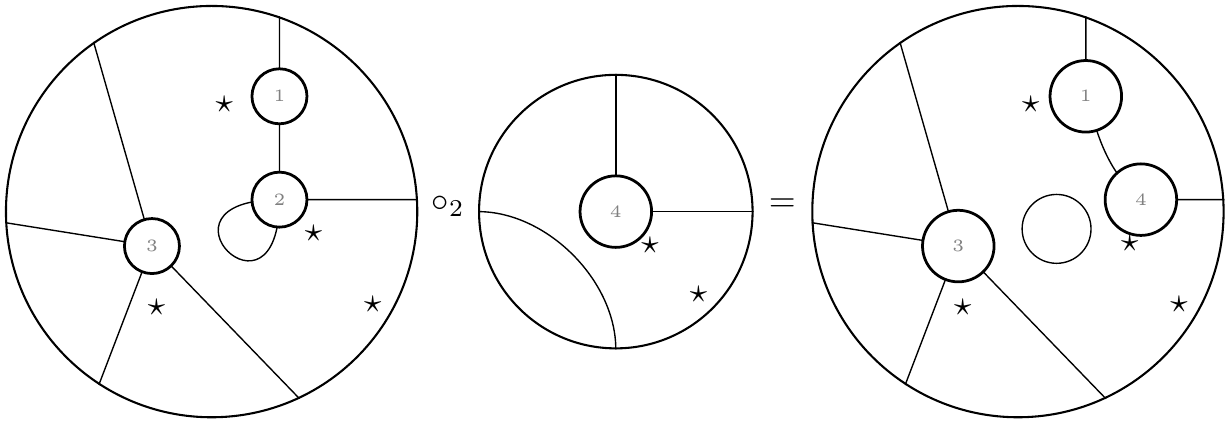}
\end{equation}
In this example, the planar diagram with inner disc 4 is inserted into inner disc 2 of the other planar diagram, giving the new planar diagram on the right.  
\begin{definition}
The planar operad, $\mathcal{P},$ is the set of isotopy classes of planar diagrams under composition as described above.  
\end{definition}
\begin{definition}
The standard involution of a planar diagram is performed in the following steps. Take any diameter of the outer disc through the starred region. Reflect the diagram over that diameter. Finally, relabel all inner and outer boundary points so that the original starred regions of each disc are preserved.
\end{definition}
As an example, consider 
\begin{equation*}
\includegraphics[valign = c, scale = .7]{diagrams/pdf/Tangle.pdf} ~~ \overset{*}{\longrightarrow} \includegraphics[valign = c, scale = .7]{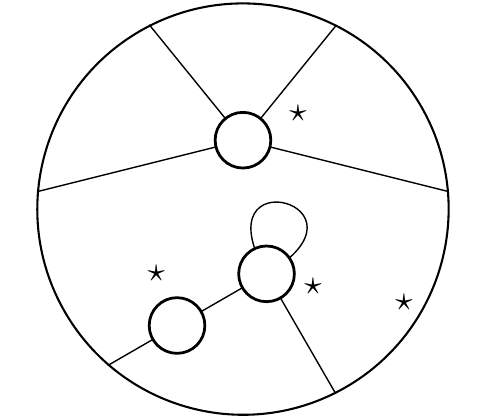}
\end{equation*}
Taking an involution can be quite cumbersome for large planar diagrams. As a solution, we draw the diagrams in what is often called ``standard'' form.
\begin{definition}\label{def:std}
A planar diagram $T$ is drawn in standard form if the input and output circles are drawn as rectangles, with strings attached only to the top or bottom and the starred region represented by the region bordering the left-hand side of the rectangle.
\end{definition}
By drawing planar diagrams in this way, an involution simply rotates the planar diagram 180 degrees but takes the complex conjugate of any coefficients. In what follows, all planar diagrams are assumed to be in standard form.
\begin{definition}
A planar algebra is a family of $\mathbbm{k}$-vector spaces $\mathcal{V}=\{V_0,V_1,V_2,\dots\},$ together with an action of the planar operad. This action assigns a multilinear map to every $T\in\mathcal{P}$ that respects the composition of planar diagrams. Moreover, if $T$ is a planar diagram with $k$ inner discs, the corresponding multilinear map is
\begin{equation*}
T:V_{i_1}\otimes V_{i_2}\otimes\dots\otimes V_{i_k} \to V_o
\end{equation*}
so that $i_m$ is the number of strands connected to the boundary of the $m$-th inner boundary disc, and $o$ is the number of strands connected to the outer boundary disc. $\mathcal{V}$ is a planar $\ast$-algebra if in addition the action of the planar operad also respects the involution on $\mathcal{P}.$
\end{definition}
Put another way, a planar algebra is a graded vector space that respects the action of the planar operad. In this paper, we will set $\mathbbm{k} = \mathbb{C}$. The dimension of a planar algebra is usually given by a sequence $a_0,a_1,a_2,\dots$ where $a_i = \dim(V_i)$ and $V_i$ is often called the $i$-box space. 

The planar operad can be thought of as the set of operations on a planar algebra. In (\ref{eq:paex}), for example, a planar algebra would determine a multilinear map
\begin{equation*}
V_2 \otimes V_4 \otimes V_4 \to V_6
\end{equation*}
Some of these operations appear often and are shown below for two elements of a planar algebra $X$ and $Y$: 
\begin{center}
\begin{tabular}{ll}
\emph{Multiplication of two diagrams:}
 &
$\begin{tikzpicture}[scale=.3, baseline = .2cm]

\draw (.5,-1) -- (.5,3);
\draw (.9,-1) -- (.9,3);
\node at (1.5,2.5) {\tiny ...};
\node at (1.5,-.5) {\tiny ...};
\draw (2.1,-1) -- (2.1,3);
\draw (2.5,-1) -- (2.5,3);

\draw[fill = white] (0,0) rectangle (3,2);
\node at (1.5,1) {\Large{{X}}};
\end{tikzpicture}
~\cdot~
\begin{tikzpicture}[scale=.3, baseline = .2cm]
\draw (.5,-1) -- (.5,3);
\draw (.9,-1) -- (.9,3);
\node at (1.5,2.5) {\tiny ...};
\node at (1.5,-.5) {\tiny ...};
\draw (2.1,-1) -- (2.1,3);
\draw (2.5,-1) -- (2.5,3);

\draw[fill = white] (0,0) rectangle (3,2);
\node at (1.5,1) {\Large{{Y}}};
\end{tikzpicture}
~ := ~
\begin{tikzpicture}[scale=.5, baseline = .4cm]
\draw (.25,-.35) -- (.25,2.7);
\draw (.45,-.35) -- (.45,2.7);
\node at (.75,1.2) {\tiny ...};
\node at (.75,-.25) {\tiny ...};
\node at (.75,2.5) {\tiny ...};
\draw (1.05,-.35) -- (1.05,2.7);
\draw (1.25,-.35) -- (1.25,2.7);

\draw[fill = white] (0,0) rectangle (1.5,1);
\draw[fill = white] (0,1.4) rectangle (1.5,2.4);
\node at (.75,.5) {X};
\node at (.75,1.9) {Y};
\end{tikzpicture}$
 \\
 \\
\emph{Capping off a diagram on the top:} 
&
$\begin{tikzpicture}[scale=.3, baseline = .2cm]

\draw (.5,-1) -- (.5,3);
\draw (.9,-1) -- (.9,3);
\node at (1.5,2.5) {\tiny ...};
\node at (1.5,-.5) {\tiny ...};
\draw (2.1,-1) -- (2.1,3);
\draw (2.5,-1) -- (2.5,3);

\draw[fill = white] (0,0) rectangle (3,2);
\node at (1.5,1) {\Large{{X}}};
\end{tikzpicture}
~\cdot~
\begin{tikzpicture}[scale = .4, baseline= -.2cm]
    \draw (.7,-1) to[out= 90, in = 180] (1.5,-.25) 
    			  to[out= 0, in = 90]  (2.3,-1);
\end{tikzpicture}
~ := ~

\begin{tikzpicture}[scale=.3, baseline = .2cm]

\draw (.5,-1) -- (.5,2.5);
\draw (.9,-1) -- (.9,2.5);
\draw (1.2,-1) -- (1.2,3);
\node at (1.8,2.5) {\tiny ...};
\node at (1.8,-.5) {\tiny ...};
\draw (2.5,-1) -- (2.5,3);
\draw (.5,2.5) to[out = 90, in = 180] (.7,2.7)
               to[out = 0, in = 90] (.9,2.5);
\draw[fill = white] (0,0) rectangle (3,2);
\node at (1.5,1) {\Large{{X}}};
\end{tikzpicture}$
\\
\\
\emph{Rotating a diagram by ``1 click:''}
&
$\rho\left(~\begin{tikzpicture}[scale=.3, baseline = .2cm]

\draw (.5,-1) -- (.5,3);
\draw (.9,-1) -- (.9,3);
\node at (1.5,2.5) {\tiny ...};
\node at (1.5,-.5) {\tiny ...};
\draw (2.1,-1) -- (2.1,3);
\draw (2.5,-1) -- (2.5,3);

\draw[fill = white] (0,0) rectangle (3,2);
\node at (1.5,1) {\Large{{X}}};
\end{tikzpicture}~\right)
~:=~
\begin{tikzpicture}[scale=.3, baseline = .2cm]

\draw (.5,-.5) -- (.5,3);
\draw (.9,-1) -- (.9,3);
\node at (1.5,2.5) {\tiny ...};
\node at (1.5,-.5) {\tiny ...};
\draw (2.1,-1) -- (2.1,3);
\draw (2.5,-1) -- (2.5,2.5);

\draw (2.5,2.5) to[out = 90, in = 180] (2.9,2.9)
				to[out = 0, in = 90]   (3.3,2.5);
\draw (3.3,2.5) -- (3.3,-1);
\draw (.5,-.5) to[out = -90, in = 0]   (.1,-.9)
			   to[out = 180, in = -90] (-.3,-.5);
\draw (-.3,-.5) -- (-.3,3);

\draw[fill = white] (0,0) rectangle (3,2);
\node at (1.5,1) {\Large{{X}}};
\end{tikzpicture}$
\end{tabular}
\end{center}
Of course, we might also want to cap a diagram off in some other way or rotate by more than 1 click, both of which are defined similarly. Note that multiplication is only defined if the number of bottom endpoints for $Y$ is the same as the number of top endpoints of $X$.
\begin{definition}\label{def:modulus}
Let $\mathcal{V}$ be a planar algebra and $f$ be the multilinear map associated with some planar diagram $T$. Then $\mathcal{V}$ is said to have modulus $d$ if for all $T \in\mathcal{P}$ inserting a closed string into a region of $T$ is equivalent to the map $d\cdot f,$ for some $d \in\mathbbm{k}^\times.$    
\end{definition}
Alternatively, one can think of a planar algebra having modulus $d$ if there exists a relation of the form:
\begin{equation*}
    \skeincircle ~=~ d
\end{equation*}
where the right-hand side of the equation is $d$ times the diagram with no strings, called the empty diagram. For this reason, we often call $d$ the circle parameter. 

\begin{definition}
A planar algebra $\mathcal{V}$ is called evaluable if $\dim V_0 = 1$, or, equivalently, any element of $V_0$ is equivalent to the empty diagram up to a scalar. 
\end{definition}

Thus, for an evaluable planar algebra we can define a map from $V_0\to\mathbbm{k}$ by sending the empty diagram to 1. This map can be used to define an inner product on each $V_n.$

\begin{definition}
Let $\mathcal{V}$ be an evaluable planar algebra. We define an inner product on $V_n$ as
\begin{center}
$\left<~
\begin{tikzpicture}[scale=.3, baseline = .2cm]

\draw (.5,-1) -- (.5,3);
\draw (.9,-1) -- (.9,3);
\node at (1.5,2.5) {\tiny ...};
\node at (1.5,-.5) {\tiny ...};
\draw (2.1,-1) -- (2.1,3);
\draw (2.5,-1) -- (2.5,3);

\draw[fill = white] (0,0) rectangle (3,2);
\node at (1.5,1) {\Large{{X}}};
\end{tikzpicture}
~,~
\begin{tikzpicture}[scale=.3, baseline = .2cm]
\draw (.5,-1) -- (.5,3);
\draw (.9,-1) -- (.9,3);
\node at (1.5,2.5) {\tiny ...};
\node at (1.5,-.5) {\tiny ...};
\draw (2.1,-1) -- (2.1,3);
\draw (2.5,-1) -- (2.5,3);

\draw[fill = white] (0,0) rectangle (3,2);
\node at (1.5,1) {\Large{{Y}}};
\end{tikzpicture}~
\right> := ~
\begin{tikzpicture}[scale=.5, baseline = .4cm]

\draw (.25,-.35) -- (.25,2.7);
\draw (.45,-.35) -- (.45,2.7);
\node at (.75,1.2) {\tiny ...};
\node at (.75,-.25) {\tiny ...};
\node at (.75,2.5) {\tiny ...};
\draw (1.05,-.35) -- (1.05,2.7);
\draw (1.25,-.35) -- (1.25,2.7);

\draw[fill = white] (0,0) rectangle (1.5,1);
\draw[fill = white] (0,1.4) rectangle (1.5,2.4);
\node at (.75,.5) {X};
\node at (.75,1.9) {Y};

\draw (1.25,2.7) to[out=90, in=180] (1.5,3)
				 to[out=0, in=90] (2,1.2)
				 to[out=-90, in=0] (1.5,-.65)
				 to[out=180, in=-90] (1.25,-.35);
\draw (1.05,2.7) to[out=90, in=180] (1.4,3.2)
				 to[out=0, in=90] (2.2,1.2)
				 to[out=-90, in=0] (1.4,-.85)
				 to[out=180, in=-90] (1.05,-.35);
\draw (.45,2.7) to[out=90, in=180] (1.2,3.6)
				 to[out=0, in=90] (2.8,1.2)
				 to[out=-90, in=0] (1.2,-1.05)
				 to[out=180, in=-90] (.45,-.35);
\draw (.25,2.7) to[out=90, in=180] (1.0,3.8)
				 to[out=0, in=90] (3.0,1.2)
				 to[out=-90, in=0] (1.0,-1.45)
				 to[out=180, in=-90] (.25,-.35);
\node at (2.5,1.2){\tiny ...};
\end{tikzpicture}$
\end{center}
\end{definition}

Of course, we could have also connected the strands on the left. If these quantities are the same, the planar algebra is said to be \emph{spherical}.

\begin{definition}
An non-zero element $v_n\in V_n$ such that $\left< v_n,w_n \right> = 0$ for all $w_n\in V_n$ is called a negligible element. 
\end{definition}

\begin{definition}\label{degenerate}
A planar algebra is \emph{non-degenerate} if it has no negligible elements.
\end{definition}

How does one find all negligible elements of $V_n?$ Given a basis for $V_n$, $\{v_1,v_2,\dots,v_k\}$, we can take the matrix of inner products, $M$, where  the $(i,j)$-th entry of $M$ is the inner product of $v_i$ and $v_j.$ The null space of this matrix would then give us the linear subspace of $V_n$ spanned by the negligible elements.

\begin{definition}\label{def:eigen}
A rotational eigenvector is a diagram $S$ such that if $\rho$ is the planar diagram that rotates the diagram one click then $\rho(S)=\zeta\cdot S,$ for some $\zeta\in\mathbb{C}.$
\end{definition}
\noindent If $\rho$ has $n$ output strands, we note that this immediately implies that $\zeta$ is an $n$-th root of unity, as $\zeta^n S =\rho^n(S) = S$. Because this space of relations is a $\mathbb{Z}/n\mathbb{Z}$ representation, the action is diagonalizable. Thus, given that $\rho^n(S) = S$ we know that $\rho(S) = \zeta\cdot S.$ 

\begin{definition}
Let $\mathcal{V}$ and $\mathcal{W}$ be planar algebras. A planar algebra map, $\Phi:\mathcal{V} \to \mathcal{W},$ is a collection of linear maps $\phi_{\pm i}: V_{\pm i} \to W_{\pm i}$ intertwining the respective actions of the planar operad.
\end{definition}
\noindent By ``intertwining" of the actions of the planar operad, we mean that if $T(v_1,\dots,v_k) = v_o$ then $T(\Phi(v_1),\dots,\Phi(v_k)) = \Phi(v_o)$ for all $T\in\mathcal{P}$ and $v_i\in\mathcal{V}$.
\begin{definition}
$\mathcal{W}$ is a sub-planar algebra of $\mathcal{V}$ if $W_i \subseteq V_i$ for all $i$ and if the inclusion maps $\{i_k:W_k \to V_k\}$ intertwine the actions of the planar operad.
\end{definition}
\begin{definition}
$\mathcal{Q}$ is a quotient of a planar algebra $\mathcal{V}$ if $Q_i$ is a quotient of $V_i$ and if the quotient maps $\{\pi_k:V_k \to Q_k\}$ intertwine the actions of the planar operad.
\end{definition}
\begin{definition}
Two planar algebras $\mathcal{V}$ and $\mathcal{W}$ are isomorphic if there is a collection of isomorphisms $\{\phi_k:V_k \to W_k\}$ that intertwine the actions of the planar operad. When the planar algebras are $\ast$-planar algebras, we also require that the isomorphism be a $\ast$-isomorphism. 
\end{definition}

\subsubsection{Skein theories}

When one defines a ring or algebra, we often describe it using generators and relations. Similarly, we can define when a set of diagrams generates a planar algebra.
\begin{definition}
A planar algebra, $\mathcal{V}$, is generated by a set of elements, $\mathcal{R},$ if for every element $v \in \mathcal{V},$ there exists a planar diagram, $T$, such that $T(r_1,\dots,r_k) = v$ for some $r_i \in \mathcal{R}.$ 
\end{definition}
When one studies all quotients of a ring or algebra, one might instead classify all two-sided ideals. Similarly, when studying the quotients of a planar algebra, we instead classify all its planar ideals.

\begin{definition}
Let $\mathcal{V}$ be a planar algebra. A planar ideal, $\mathcal{I}$, is a subset of $\mathcal{V}$ such that $\mathcal{I}$ is the kernel of some quotient map of $\mathcal{V}.$ 
\end{definition}

Equivalently, $\mathcal{I}$ is a planar ideal if it is closed under arbitrary composition with any element of the planar algebra. That is, if we let $T$ be any planar tangle, then $T(v_1, \dots, v_k)\in\mathcal{I}$ if at least one $v_i$ is an element of $\mathcal{I}$ for all possible choices of the other elements. Just like a two-sided ideal in a ring, a planar ideal can also be written using relations, which we call skein relations.
\begin{definition}
Let $\mathcal{V}$ be a planar algebra. A \emph{skein relation} is an element of the kernel of a quotient map from $\mathcal{V}$.  A \emph{skein theory} is a set of skein relations that generates a planar ideal.
\end{definition}
Thus, giving a skein theory would be an equivalent way of defining a quotient of a planar algebra. If we wanted to classify all quotients of a planar algebra, then, we could instead classify all skein theories of that planar algebra. In Section \ref{examples}, we will examine a number of well-known examples of skein theories.

\subsection{Oriented planar algebras}

Our definition of a planar algebra assumes that the strings of planar diagrams are unoriented. One could imagine defining a planar algebra, though, by giving an action on planar diagrams with oriented strings, the collection of which is appropriately called the oriented planar operad. In this case, the orientation of the strings gives orientations to the endpoints, either positive or negative. If we imagine every endpoint of a diagram on the bottom, then an endpoint will be positively oriented if the strand connecting to it is pointing away from the endpoint and negatively oriented if it is pointing toward the endpoint.
\begin{definition}
An oriented planar algebra $\mathcal{V}$ is a collection of vector spaces $V_{(z_1,\dots,z_i)}$ where $z_j\in\{+,-\}$ together with an action of the oriented planar operad.
\end{definition}

Elements of an oriented planar algebra are also drawn diagrammatically with oriented strands. The $z_i$ tell us the orientation of each endpoint. For example, elements of $V_{(+,+,-,-)}$ include any diagram that can be drawn in the following form:
\begin{center}
\includegraphics[valign = c, scale = 1.5]{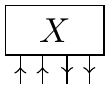}
\end{center}
Other concepts we have defined for unoriented planar algebras have oriented analogs, but those definitions are omitted here.

\subsection{Shaded planar algebras}\label{shadedplanaralgebra}

Suppose that we impose the additional requirement that diagrams in the planar operad  be checkerboard shaded. As an example, consider 
\begin{equation}\label{shadedtangle}
    T~:=~\includegraphics[valign = c]{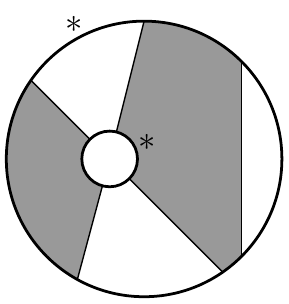}
\end{equation}
There are two main differences to note between shaded planar diagrams and their unshaded counterparts. First, in order to achieve a checkerboard shading it is clear that every disc must have an even number of boundary points. Second, the starred region can now be shaded or unshaded. Both of these differences are reflected in the definition of a shaded planar algebra.    

\begin{definition}
A shaded planar algebra is a collection of vector spaces $\mathcal{V} = \{V_{\pm i}\},$ $i \in \mathbb{Z}_{\geq 0}$ together with an action of the planar operad. 
\end{definition}

Here, $\mathcal{V}_{\pm i}$ represents the vector spaces of diagrams with $2i$ endpoints, which are also traditionally checkerboard-shaded. Thus, a shaded diagram with four endpoints is a 2-box, while an unshaded diagram with four endpoints is a 4-box. If the starred region is unshaded with $2i$ endpoints, then it is an element of $V_{+i}$. As an example, (\ref{shadedtangle}) defines the map $T:V_{-2} \to V_{+3}.$ Examples of shaded planar algebras appear in Section \ref{examples}.
\begin{definition}
Let $\mathcal{V} = \{V_{\pm i}\}_{i = 0, 1, \dots}$ be a shaded planar algebra. $\mathcal{V}^\text{even}$ is the collection of vector spaces $\{V_{+0},V_{+2}, V_{+4}, \dots\}$ (i.e. the $2i$-box spaces with unshaded starred region). 
\end{definition}
\noindent The even portion of any shaded planar algebra is necessarily unshaded by assigning
\begin{center}
    \includegraphics[valign = c]{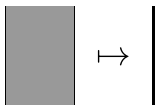}
\end{center}
Thus, when we restrict to the even portion of maps of shaded planar algebras, we can often find interesting maps of unshaded planar algebras. In Section  \ref{connection}, we will see an interesting example of how one can obtain virtual knot and link invariants from a certain map of shaded planar algebras. In that example, though, it turns out that the planar algebra in questions will be too large for us to handle, so we will use the following method to make it more manageable:

\begin{definition}\label{cutdown}
Let $\mathcal{V}$ be a shaded planar algebra and $P$ be some projection in in $\mathcal{V}_{\pm i}.$  A cut-down of $\mathcal{V}$ is the planar algebra given by the map
\begin{center}
    \includegraphics[valign = c]{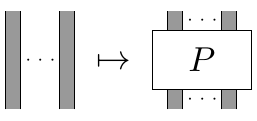}
\end{center}
when $P$ is a $+i$-box and similarly when it is a $-i$-box.
\end{definition}

When $P\in V_{\pm 2i}$ then $\mathcal{V}$ is an unshaded planar algebra, while it is still shaded if $P\in V_{\pm (2i + 1)}$. In this paper, we will only utilize projections in the 2-box space, and so all the resulting planar algebras will be unshaded.

\subsubsection{Planar algebra of a bipartite graph}\label{starn}

Another class of shaded planar algebras is the planar algebra of a bipartite graph (originally defined in \cite{Jon00}). Although we define this planar algebra in full generality here, we will only make use of the planar algebra of a special bipartite graph, $\ast_n,$ which is also called the spin planar algebra in \cite{Jon99}. 

Let $\mathcal{G}$ be a locally-finite, connected bipartite graph (not necessarily simply-laced) with edge set $E$, vertex set $U = U^+\cup U^-,$  with $|U^+| = n_1$, $|U^-| = n_2$, and $|U| = n_1 + n_2 = n$. Further, there are no edges that connect $U^+$ or $U^-$ to itself. Because $\mathcal{G}$ is a bipartite graph, its adjacency matrix is of the form 
\begin{equation*}
\left(
\begin{array}{cc}
0 & \Lambda \\
\Lambda^T & 0
\end{array}
\right)
\end{equation*}
where $\Lambda$ is a $n_1 \times n_2$ matrix and the $(i,j)$-th entry is the number of edges connect $v^+_i$ to $v^-_j.$ In this project, we care about a particular class of graphs, which we call $\ast_n$:

\begin{definition}
$\ast_n$ is the complete bipartite graph with $U^- = \{\star\}$ and $U^+ = S$, a set of $n$ even vertices. 
\end{definition}

For each $k \in \mathbb{Z}_{\geq 0}$ we define the basis of $L_{\pm k}$ to be loops in $\mathcal{G}$ of length $k,$ with one vertex of the loop declared the ``base." Loops in $L_{+k}$ are based in $U^+$ and loops in $L_{- k}$ are based in $U^-$. We will denote such a loop with a pair of functions $(\pi, \epsilon):\{0,1,\dots,k-1\} \to U\cup E,$ where the $\pi(i)$ is the $(i+1)$-th vertex of the loop and $\epsilon(i)$ is the $(i+1)$-th edge. The astute reader might note that given an edge, one can determine the corresponding vertices and recreate the loop simply from $\epsilon.$ In many of the examples that follow, however, the graphs are simply-laced and $\epsilon$ will be suppressed in favor of $\pi.$

The following definition was given in \cite{Jon00}: 
\begin{definition}
The planar algebra of a bipartite graph $\mathcal{G}$, $\GPA(\mathcal{G}),$ is a collection of vector spaces $\{GPA(\mathcal{G})_{\pm k}\}$ where $\GPA(\mathcal{G})_{\pm k} = \{f:L_{\pm k} \to \mathbb{C} \}$, the set of functionals of based loops in $\mathcal{U}^\pm$ of length $k$, together with an action of the planar operad.
\end{definition}
For any bipartite graph, $\mathcal{G},$ it is true that $L_{\pm 0} = U^\pm$. Thus, a basis for $\GPA(\mathcal{G})_{\pm0} = \{\delta_v\mid v\in U^\pm\},$ where $\delta_v$ is the Kronecker delta
\begin{equation}
\delta_v(w) :=
\left\{
\begin{array}{cl}
1 & w = v \\
0 & \text{else} 
\end{array}
\right.
\end{equation}   
We can generalize this to the following proposition:

\begin{proposition}\label{prop:basis}
For any $k$, a basis for $\GPA(\mathcal{G})_{\pm k}$ is given by $\{\delta_l\mid l \in L_{\pm k}\}.$ 
\end{proposition}

\begin{proof}
A functional $f:X\to\mathbb{C}$ is a linear mapping of every element of $X$ to a complex number. In our case, $X=L_k$, so a functional is an assignment of every loop in $L_k$ to a complex number. Then we can rewrite each functional as 
\begin{equation*}
f = \sum_{l \in L_k} \delta_l\cdot f(l)
\end{equation*}
Thus, every such functional is a linear combination of the $\delta_l.$ Additionally, since each of the $\delta_l$ takes a value of $1$ on a unique value, the set is clearly linearly independent. Thus, they must form a basis, as desired.
\end{proof}

To form this set of vector space into a planar algebra, though, we must also define an action of the planar operad. This is accomplished using ``states":
\begin{definition}
A state $\sigma$ on a planar diagram $T$ is an assignment of vertices to regions and edges to strings which is consistent with the graph structure of $\mathcal{G}$. That is, given an assignment of vertices $v$ and $w$ to neighboring regions, the edge $e$ assigned to the string separating those regions should have endpoints $v$ and $w$.
\end{definition}
We note here that for a simply-laced graph, this labeling of strings can be omitted if we instead require that adjacent regions must be labeled by adjacent vertices. If $\sigma$ is a state, then we define $\sigma|_i$ to be the based loop obtained by reading the states assigned to the regions of the $i$-th input disc with the starred region being the base, and $\sigma|_{\text{boundary circle}} = \sigma|_o$ is defined similarly.

\begin{definition}[Action 1] Let $\text{GPA}(\mathcal{G})$ be as described above, and let $T$ be a planar tangle. $T$ acts on functionals $f_1, \dots, f_k$ by sending $(f_1, \dots, f_k)$ to the functional described below: 
\begin{equation}
T(f_1,f_2,\dots,f_k)(l) = 
\sum_{\substack{\text{states on } T \\
 \text{ such that }\\
   \sigma|_o = l}}
~ \prod_{i=1}^k f_i(\sigma|_i)
\end{equation}
\end{definition}
\noindent %
As an example, consider the graph $\ast_5$
\begin{center}
\includegraphics[valign = c]{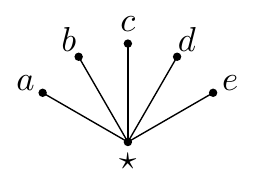}
\end{center}
and the planar tangle
\begin{center}
T ~:=~ \includegraphics[valign = c]{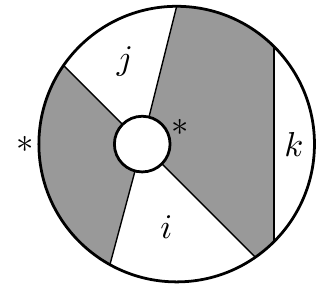}
\end{center}
where $i,j,k$ are elements of $U^+$. The shaded regions are labelled with $\star$ but are omitted to avoid confusion with the starred regions. When our graph is $\ast_n$, we will write loops as a concatenation of its vertices omitting $\star$. Thus, the loop $(\star, a, \star, b)$ would be written as $ab$. By definition of Action 1 of the planar operad, when $\delta_{ab}$ is inputted we obtain the functional
\begin{equation*}
    T(\delta_{ab})(ijk) = \left\{
    \begin{array}{cc}
        1 & i=a \text{ and } j=b  \\
        0 & \text{else} 
    \end{array}\right.
    ~ = ~\sum_{k\in S} \delta_{abk}
\end{equation*}
One can extend this result linearly to get a functional for any choice of input. By considering the action of  $\mathcal{V}$ on planar tangles with no input discs, we see that
\begin{center}
\begin{tabular}{ccccc}
     $\begin{tikzpicture}[baseline = -.1cm]
    \filldraw[fill = black!40!] (-.5,-.5) to[out = 90, in = 180] (0,-.2) to[out = 0, in = 90] (.5,-.5);
    \filldraw[fill = black!40!] (-.5,.5) to[out = -90, in = 180] (0,.2)  to[out = 0, in = -90] (.5,.5);
    \node at (-.5,0) {$a$};
    \node at (.5,0) {$b$};
\end{tikzpicture} ~=~ \left\{
    \begin{array}{cc}
        1 & a = b  \\
        0 & a \ne b 
    \end{array}
    \right.  $ && and &&
    \begin{tikzpicture}[baseline = -.1cm]
    \draw[fill = black!40!, draw = none] (-.485,-.5) rectangle (.485,.5);
    \filldraw[fill = white] (.5,.5) to[out = -135, in = 90] (.2,0) to[out = -90, in = 135]  (.5,-.5);
    \filldraw[fill = white] (-.5,.5) to[out =-45, in = 90] (-.2,0) to[out = -90, in = 45] (-.5,-.5);
    \node at (-.5,0) {$a$};
    \node at (.5,0) {$b$};
\end{tikzpicture} = 1
\end{tabular}
\end{center}
for all $a,b\in S.$ Thus, we know that 
\begin{equation}
  \begin{tikzpicture}[baseline = -.1cm]
    \filldraw[fill = black!40!] (-.5,-.5) to[out = 90, in = 180] (0,-.2) to[out = 0, in = 90] (.5,-.5);
    \filldraw[fill = black!40!] (-.5,.5) to[out = -90, in = 180] (0,.2)  to[out = 0, in = -90] (.5,.5);
\end{tikzpicture} ~=~ \sum_{i=1}^n \delta_{ii}
\end{equation}
and
\begin{equation}
  \begin{tikzpicture}[baseline = -.1cm]
    \draw[fill = black!40!, draw = none] (-.485,-.5) rectangle (.485,.5);
    \filldraw[fill = white] (.5,.5) to[out = -135, in = 90] (.2,0) to[out = -90, in = 135]  (.5,-.5);
    \filldraw[fill = white] (-.5,.5) to[out =-45, in = 90] (-.2,0) to[out = -90, in = 45] (-.5,-.5);
\end{tikzpicture} ~=~ \sum_{i,j=1}^n \delta_{ij}
\end{equation}
Let us now consider the following planar diagram: 
\begin{equation*}
\includegraphics[valign = c]{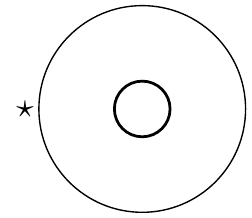}
\end{equation*} 
Note here that the center circle is not an input disc but instead a closed string.  Using the state sum formula, there are no functionals to input so for $v\in U^-$ 
\begin{equation*}
\includegraphics[valign = c]{diagrams/pdf/statesumex2.pdf} ~ (v) = \deg(v)
\end{equation*}
where $\deg(v)$ is the number of neighbors of $v$ in $\mathcal{G}$. In our example, the only $v\in U^-$ is $\star$, which has degree 5. The planar diagram with the opposite shading would give 1 for all possible choices of input in $U^+$.

The reader can verify that Action 1 respects the composition of planar tangles and, thus, gives an action of the planar operad on $\GPA(\mathcal{G}).$  We also noted, however, that the value of a closed string in a diagram was dependent upon the degree of the vertex under this action. When $\mathcal{G} = \ast_n$, however, we have that 
\begin{center}
\begin{tabular}{cccc}
$\begin{tikzpicture}[baseline]
\draw[fill = black!40!] circle (.5cm);
\end{tikzpicture} ~=~ 1$ &&& $\begin{tikzpicture}[baseline]
\draw[fill = black!40!, draw = none] (-.6,-.6) rectangle (.6,.6);
\draw[fill = white] circle (.5cm);
\end{tikzpicture} ~=~ n$
\end{tabular}
\end{center}
In order for the circle to have a fixed value, $d,$  we will redefine the action following the construction in \cite{Jon00}. Before we do this, however, we need to define the Perron-Frobenius eigenvector.
\begin{definition} The Perron-Frobenius eigenvector of a graph is a map $\omega: V(\mathcal{G}) \to \mathbb{R}$ such that 
\begin{enumerate}[label = \roman{enumi}.]
\item $\omega(v) > 0$ for all $v\in V(\mathcal{G})$
\item There is a $d\in\mathbb{R}_{\geq 0}$ such that $\displaystyle \sum_{w \text{ adj } v} \omega(w) = d\cdot \omega(v)$.
\end{enumerate}
\end{definition}
It is well-known that $\omega$ is an eigenvector of the adjacency matrix of $\mathcal{G}$ corresponding to the largest eigenvalue $d$. It is also proven that such a $d$ and function $\omega$ can be scaled to meet the above requirements for any locally finite, connected graph. For, $\ast_n$, the largest eigenvalue is $\sqrt{n}$ with corresponding eigenvector
\begin{equation*}
\omega = \left\{\sqrt{n},1,1,\dots,1 \right\}
\end{equation*} 
where $\star$ represents the first coordinate of $\omega$. We now define another action of the planar operad on $\GPA(\mathcal{G})$ using $\omega$:
\begin{definition}[Action 2] For a bipartite graph $\mathcal{G}$, a planar diagram $T$ acts on $f_1,\dots,f_k$ by
\begin{equation}
T(f_1,f_2,\dots,f_k)(l) = 
\sum_{\substack{\text{states on } T \\
 \text{ such that }\\
   \sigma|_o = l}}
c(T,\sigma) 
 \prod_{i=1}^k f_i(\sigma|_i)
\end{equation}
where  $c(T,\sigma)$ is some constant depending on $\sigma$, $T$, and $\omega,$ the Perron-Frobenius eigenvector of $\mathcal{G}$.
\end{definition}
In \cite{Jon00}, the following theorem is proven:
\begin{theorem}[\cite{Jon00}] The definition of Action 2 makes $\GPA(\mathcal{G})$ into a planar algebra. This planar algebra has modulus $d,$ the largest eigenvalue of the adjacency matrix of $\mathcal{G}$. Further, if $k=\mathbb{C}$ then this action defines a planar $C^*$-algebra.
\end{theorem}
Thus, any planar algebra defined using Action 1 can be re-normalized by using Action 2.

\subsubsection{Spin Models}\label{SpinModel}

The following definition is based on \cite{Jon99}. More details about spin models can also be found in \cite{Edg18}:
\begin{definition}
Let $\mathcal{V}$ be a shaded planar algebra. A spin model for $\mathcal{V}$ is a map of planar algebras $\Phi:\mathcal{V}\to\text{GPA}(\ast_n)$ for some $n$.
\end{definition}
Suppose $\mathcal{V}$ is generated by elements of $V_{+2}$. Recall that a basis for $\text{GPA}(\ast_n)_{2k}$ is the collection of Kronecker deltas of loops,  $\delta_{l}$, where $l$ is any $k$-tuple of elements of $S$. If $X$ is an element of $V_{+2}$ then, in general, $\Phi(X) =\displaystyle \sum_{a,b\in S} c_X(a,b)\cdot\delta_{ab},$ where $c_X(a,b) \in\mathbb{C}.$ Let $C_X$ be the $n\times n$ matrix where the $(a,b)$ entry of $C_X$ is $c_X(a,b).$ By the action of the planar operad on $\text{GPA}(\ast_n),$ we know that the matrix corresponding to \begin{tikzpicture}[baseline = -.1cm, scale = .6]
    \filldraw[fill = black!40!] (-.5,-.5) to[out = 90, in = 180] (0,-.2) to[out = 0, in = 90] (.5,-.5);
    \filldraw[fill = black!40!] (-.5,.5) to[out = -90, in = 180] (0,.2)  to[out = 0, in = -90] (.5,.5);
\end{tikzpicture} is the identity matrix; similarly, the matrix corresponding to \begin{tikzpicture}[baseline = -.1cm, scale = .6]
    \draw[fill = black!40!, draw = none] (-.475,-.5) rectangle (.475,.5);
    \filldraw[fill = white] (.5,.5) to[out = -135, in = 90] (.2,0) to[out = -90, in = 135]  (.5,-.5);
    \filldraw[fill = white] (-.5,.5) to[out =-45, in = 90] (-.2,0) to[out = -90, in = 45] (-.5,-.5);
\end{tikzpicture} is the matrix of all ones.  

To capture the fact that $\Phi$ is a map of planar algebras, we will write 
\begin{equation*}
c_X(a,b) ~:=~ \xdiag[X]  
\end{equation*}
for all $a,b\in S$ and generators $X$. More generally we can give the following definition:

\begin{definition}
Suppose $\Phi$ gives a spin model for some shaded planar algebra, $\mathcal{V}$, that is generated by 2-boxes and let $X\in V_{\pm k}$. Then a state sum with respect to $(a_1,\dots,a_k)$ is the sum over all possible labelings of the unshaded regions of $X$ such that the $i$th exterior region of $X$ (starting with the starred region and going clockwise) is $a_i.$  
\end{definition}
As an example, the state sum of 
\begin{center}
    \includegraphics[valign = c]{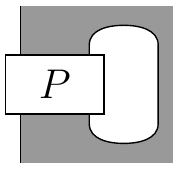}
\end{center}
with respect to $a$ is
\begin{equation*}
\includegraphics[valign = c]{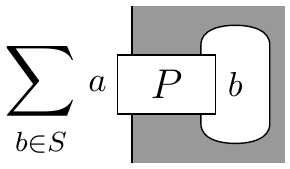}  
\end{equation*}
Just like for $c_X(a,b)$ for some generator $X$, one should think of the state sum of any diagram $Y$ with respect to $(a_1,\dots,a_i)$ as the coefficient of $\delta_{a_1\cdots a_i}$ of $\Phi(Y).$ 

\section{Examples}\label{examples}
We now discuss in detail a number of well-known planar algebras that will arise in our classifications.

\subsection{Planar algebra of unoriented tangles}
\begin{definition}\label{ex:vtangle}
The planar algebra of unoriented tangles is the planar algebra generated by $\begin{tikzpicture}[baseline = -.15cm, scale = .65]
\begin{knot}[clip width = 4]
	\strand (-.5,.5) -- (.5,-.5);
	\strand (-.5,-.5) -- (.5,.5);
\end{knot}
\end{tikzpicture}
\!$. Diagrams are considered up to isotopy. In addition, we will impose the following relations:
\begin{center}
\begin{tabular}{ccc}
Reidemeister II (R2): && $\begin{tikzpicture}[baseline=(current bounding box.center), rounded corners = 3mm]
    \draw (0,0) -- (.75,.5) -- (0,1);
    \draw (.75,0) -- (.45,.2);
    \draw (.25,.3) -- (0,.5) -- (.3,.7);
    \draw (.45,.8)--(.75,1);
\end{tikzpicture}
 = \begin{tikzpicture}[baseline, rounded corners = 5mm]
    \draw (.5,.5) -- (.1,0) -- (.5,-.5);
    \draw (-.5,.5) -- (-.1,0) -- (-.5,-.5);
\end{tikzpicture}$
\\
\\
Reidemeister III (R3): && $\begin{tikzpicture}[scale = .5, baseline]
\draw (1,1) -- (-1,-1);
\draw (1,-1) -- (.2,-.2);
\draw (-.2,.2) -- (-1,1);
\draw (0,1) to[out = -45, in = 135] (.4,.65);
\draw (.6,.45) to[out = -45, in = 90] (.8,0)
               to[out = -90, in = 45] (.6,-.45);
\draw (.4,-.65) to[out = 45, in = -135] (0,-1);       
\end{tikzpicture}
=
\begin{tikzpicture}[scale = .5, baseline]
\draw (1,1) -- (-1,-1);
\draw (1,-1) -- (.2,-.2);
\draw (-.2,.2) -- (-1,1);
\draw (0,1) to[out = 45, in = -135] (-.4,.65);
\draw (-.6,.45) to[out = -135, in = 90] (-.8,0)
                to[out = -90, in = 135] (-.6,-.45);
\draw (-.4,-.65) to[out = -45, in = 135] (0,-1);   
\end{tikzpicture}$
\end{tabular}
\end{center}
\end{definition}
We will frequently refer to these relations by R2 and R3 for simplicity. 
Reidemeister I, which is given by
\begin{equation*}
\begin{tikzpicture}[baseline=.5cm, rounded corners = 5mm]
\draw (0,1)--(.4,.6);
\draw (.6,.4)--(1,0)--(0,0)--(1,1);
\end{tikzpicture} ~=~ \begin{tikzpicture}[baseline = .5 cm, rounded corners = 6mm]
    \draw (0,1)--(.5,-.1) -- (1,1);
\end{tikzpicture} 
\end{equation*}
is \emph{not} required to be satisfied. In all of our following examples, however, we will require that R1 be satisfied up to a scalar. A relation of this form is often called a twist relation:
\begin{equation}\label{eq:R1}
\begin{tikzpicture}[baseline=.5cm, rounded corners = 5mm]
\draw (0,1)--(.4,.6);
\draw (.6,.4)--(1,0)--(0,0)--(1,1);
\end{tikzpicture} ~=~ a\cdot\begin{tikzpicture}[baseline = .5 cm, rounded corners = 6mm]
    \draw (0,1)--(.5,-.1) -- (1,1);
\end{tikzpicture} 
\end{equation}
Note that (\ref{eq:R1}) and R2 imply that 
\begin{equation*}
\begin{tikzpicture}[baseline=.5cm, rounded corners = 5mm]
\draw (0,1)--(1,0) -- (0,0)--(.4,.4);
\draw (.6,.6)--(1,1);
\end{tikzpicture} ~=~ a^{-1}\cdot\begin{tikzpicture}[baseline = .5 cm, rounded corners = 6mm]
    \draw (0,1)--(.5,-.1) -- (1,1);
\end{tikzpicture} 
\end{equation*}

Clearly, this planar algebra has trivial odd-box spaces. The even-box spaces, on the other hand, are infinite dimensional! This is evident even in the 0-box space, as every non-equivalent knot and link are linearly independent. Rather than studying the planar algebra of unoriented tangles itself, we will instead look at quotients of this planar algebra by studying their skein theories. A few well-known skein theories of the planar algebra of unoriented tangles are the Jones and Kauffman polynomials, originally discussed in \cite{Jon85} and \cite{Kau90}.

\subsubsection{Braided planar algebras}
Before we discuss these quotients, we give a general definition of a braided planar algebra using the relations above.
\begin{definition}\label{def:naturality}
A planar algebra, $\mathcal{V},$ is said to be \emph{braided} if there exists an element of $V_4$ that satisfies a twist relation, R2, and R3. In addition, the braiding must satisfy the following naturality condition:
\begin{equation*}
\includegraphics[valign = c,scale = 1.25]{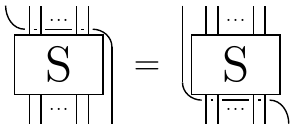}
\end{equation*}
for any element $S\in \mathcal{V}.$
\end{definition}
We will classify all such braidings of symmetric trivalent planar algebras in Section \ref{tri}. While some planar algebras have elements that satisfy the Reidemeister relations, they do not always act naturally with every element in the planar algebra. In these instances, although the entire planar algebra is not braided, it is true that the planar algebra has a subplanar algebra that is braided.
\begin{definition}\label{subbraid}
A planar algebra, $\mathcal{V}$ is \emph{sub-braided} if it contains diagrams \begin{tikzpicture}[baseline = -.15cm, scale = .65]
\begin{knot}[clip width = 4]
	\strand (-.5,.5) -- (.5,-.5);
	\strand (-.5,-.5) -- (.5,.5);
\end{knot}
\end{tikzpicture}
 and \begin{tikzpicture}[baseline = -.15cm, scale = .65]
\begin{knot}[clip width = 4]
	\strand (-.5,-.5) -- (.5,.5);
	\strand (-.5,.5) -- (.5,-.5);
\end{knot}
\end{tikzpicture} satisfying a twist relation, R2, and R3 and the naturality condition on some sub-planar algebra $\mathcal{W}\subseteq\mathcal{V}$. 
\end{definition}

\begin{example}[Temperley-Lieb Jones polynomial/Kauffman bracket planar algebra]\label{Jones}

The Jones polynomial planar algebra---also called the Temperley-Lieb-Jones (TLJ) planar algebra---is a quotient of the planar algebra of tangles with the following relations:
\begin{center}
\begin{tabular}{cccc}
$\skeincircle ~=~ -(A^2+A^{-2})$ &&& $\begin{tikzpicture}[baseline=.5cm, rounded corners = 5mm]
\draw (0,1)--(.4,.6);
\draw (.6,.4)--(1,0)--(0,0)--(1,1);
\end{tikzpicture} = -A^{-3}\begin{tikzpicture}[baseline = .5 cm, rounded corners = 6mm]
    \draw (0,1)--(.5,-.1) -- (1,1);
\end{tikzpicture}$
\end{tabular}
\end{center}
\begin{equation*}
\begin{tikzpicture}[baseline]
\begin{knot}[clip width = 4]
	\strand (-.5,.5) -- (.5,-.5);
	\strand (-.5,-.5) -- (.5,.5);
\end{knot}
\end{tikzpicture}
 = A\begin{tikzpicture}[baseline, rounded corners = 5mm]
    \draw (.5,.5) -- (.1,0) -- (.5,-.5);
    \draw (-.5,.5) -- (-.1,0) -- (-.5,-.5);
\end{tikzpicture} + A^{-1}\begin{tikzpicture}[baseline, rounded corners = 5mm]
    \draw (-.5,-.5) -- (0,.1) -- (.5,-.5);
    \draw (-.5,.5) -- (0,-.1) -- (.5,.5);
\end{tikzpicture}
\end{equation*}
In the last relation, we will often call ~\begin{tikzpicture}[baseline={([yshift=-1ex]current bounding box.center)}, rounded corners = 5mm, scale = .4]
    \draw (.5,.5) -- (.1,0) -- (.5,-.5);
    \draw (-.5,.5) -- (-.1,0) -- (-.5,-.5);
\end{tikzpicture}~ the identity while we denote ~\begin{tikzpicture}[baseline={([yshift=-1ex]current bounding box.center)}, rounded corners = 5mm, scale = .4]
    \draw (-.5,-.5) -- (0,.1) -- (.5,-.5);
    \draw (-.5,.5) -- (0,-.1) -- (.5,.5);
\end{tikzpicture}~ the cupcap. It can easily be verified that this planar algebra is braided. It can also be shown that if $\mathcal{V}$ is a quotient of the planar algebra of tangles with initial dimensions $1,0,1,0,2$ then $\mathcal{V}$ is a specialization of the TLJ planar algebra. 
\end{example}

\begin{example}[The Kauffman/Dubrovnik polynomial planar algebra]\label{dubkau}

The Kauffman/Dubrovnik polynomial planar algebra is a generalization of the TLJ planar algebra. The skein relations are 
\begin{center}
$\begin{tikzpicture}[baseline]
\begin{knot}[clip width = 4]
	\strand (-.5,.5) -- (.5,-.5);
	\strand (-.5,-.5) -- (.5,.5);
\end{knot}
\end{tikzpicture}
 ~\pm~ \begin{tikzpicture}[baseline]
\begin{knot}[clip width = 4]
	\strand (-.5,-.5) -- (.5,.5);
	\strand (-.5,.5) -- (.5,-.5);
\end{knot}
\end{tikzpicture} = z \left( ~ \begin{tikzpicture}[baseline, rounded corners = 5mm]
    \draw (.5,.5) -- (.1,0) -- (.5,-.5);
    \draw (-.5,.5) -- (-.1,0) -- (-.5,-.5);
\end{tikzpicture} \pm \begin{tikzpicture}[baseline, rounded corners = 5mm]
    \draw (-.5,-.5) -- (0,.1) -- (.5,-.5);
    \draw (-.5,.5) -- (0,-.1) -- (.5,.5);
\end{tikzpicture} ~ \right)$
\end{center}
\begin{center}
\begin{tabular}{ccc}
$\skeincircle = d$	&& $\begin{tikzpicture}[baseline=.5cm, rounded corners = 5mm]
\draw (0,1)--(.4,.6);
\draw (.6,.4)--(1,0)--(0,0)--(1,1);
\end{tikzpicture} = a\cdot\begin{tikzpicture}[baseline = .5 cm, rounded corners = 6mm]
    \draw (0,1)--(.5,-.1) -- (1,1);
\end{tikzpicture}$
\end{tabular}
\end{center}
The $+$ case is called the Kauffman polynomial planar algebra and the $-$ case is called the Dubrovnik polynomial planar algebra. By capping off the first equation on the bottom, we get
\begin{equation*}
\begin{tikzpicture}[baseline=.5cm, rounded corners = 5mm]
\draw (0,1)--(1,0) -- (0,0)--(.4,.4);
\draw (.6,.6)--(1,1);
\end{tikzpicture} ~\pm~ \begin{tikzpicture}[baseline=.5cm, rounded corners = 5mm]
\draw (0,1)--(.4,.6);
\draw (.6,.4)--(1,0)--(0,0)--(1,1);
\end{tikzpicture} = z \left( ~ \begin{tikzpicture}[baseline = .5 cm, rounded corners = 6mm]
    \draw (0,1)--(.5,-.1) -- (1,1);
\end{tikzpicture} \pm d\cdot\begin{tikzpicture}[baseline = .5 cm, rounded corners = 6mm]
    \draw (0,1)--(.5,-.1) -- (1,1);
\end{tikzpicture} ~\right)
\end{equation*}
\noindent Thus, $z(1\pm d) = a^{-1} \pm a,$ which gives a two-parameter family of relations in terms of $a$ and $d$. If one sets $a = -A^{-3}$ and $d = -(A^2 - A^{-2}),$ note that one does not obtain the TLJ planar algebra. At these values of the parameters, the Kauffman polynomial planar algebra is still defined but is degenerate. See \cite{BW90} and \cite{MT90} for more information on this. The quotient of that planar algebra by the negligible elements of this planar algebra, however, is isomorphic to the TLJ planar algebra.
\end{example}

A fact that is proven in \cite{MPS11} and that we will employ in Section \ref{virtualskein} is recorded here:

\begin{theorem}\label{thm:kauffman}
Let $\mathcal{V}$ be a quotient of the planar algebra of tangles such that $\dim V_4 = 2$ or $3$. Then there exists a relation of the form
\begin{equation*}
\begin{tikzpicture}[baseline]
\begin{knot}[clip width = 4]
	\strand (-.5,.5) -- (.5,-.5);
	\strand (-.5,-.5) -- (.5,.5);
\end{knot}
\end{tikzpicture}
 ~\pm~ \begin{tikzpicture}[baseline]
\begin{knot}[clip width = 4]
	\strand (-.5,-.5) -- (.5,.5);
	\strand (-.5,.5) -- (.5,-.5);
\end{knot}
\end{tikzpicture} = z \left( ~ \begin{tikzpicture}[baseline, rounded corners = 5mm]
    \draw (.5,.5) -- (.1,0) -- (.5,-.5);
    \draw (-.5,.5) -- (-.1,0) -- (-.5,-.5);
\end{tikzpicture} \pm \begin{tikzpicture}[baseline, rounded corners = 5mm]
    \draw (-.5,-.5) -- (0,.1) -- (.5,-.5);
    \draw (-.5,.5) -- (0,-.1) -- (.5,.5);
\end{tikzpicture} ~ \right)
\end{equation*}
\end{theorem}

\begin{proof}
Since $\dim V_4 = 2$ or $3$, we know there must be a relation amongst the two crossings and the Temperley-Lieb diagrams which we can generically write as
\begin{equation*}
a_1 \cdot \begin{tikzpicture}[baseline]
\begin{knot}[clip width = 4]
	\strand (-.5,.5) -- (.5,-.5);
	\strand (-.5,-.5) -- (.5,.5);
\end{knot}
\end{tikzpicture}
 ~+~ a_2 \cdot \begin{tikzpicture}[baseline]
\begin{knot}[clip width = 4]
	\strand (-.5,-.5) -- (.5,.5);
	\strand (-.5,.5) -- (.5,-.5);
\end{knot}
\end{tikzpicture} ~=~ a_3 \cdot \begin{tikzpicture}[baseline, rounded corners = 5mm]
    \draw (.5,.5) -- (.1,0) -- (.5,-.5);
    \draw (-.5,.5) -- (-.1,0) -- (-.5,-.5);
\end{tikzpicture} ~+~ a_4 \cdot \begin{tikzpicture}[baseline, rounded corners = 5mm]
    \draw (-.5,-.5) -- (0,.1) -- (.5,-.5);
    \draw (-.5,.5) -- (0,-.1) -- (.5,.5);
\end{tikzpicture}
\end{equation*}
By rotating the equation however, we see that $a_1 = \pm a_2$ and $a_3 = \pm a_4$. If $a_1 = a_2 = 0$, then the cupcap is a multiple of the identity, which implies that $\dim V_4 = 1.$ Thus, $a_1 \ne 0$, and so we can rewrite the equation as 
\begin{equation*}
\begin{tikzpicture}[baseline]
\begin{knot}[clip width = 4]
	\strand (-.5,.5) -- (.5,-.5);
	\strand (-.5,-.5) -- (.5,.5);
\end{knot}
\end{tikzpicture}
 ~\pm~ \begin{tikzpicture}[baseline]
\begin{knot}[clip width = 4]
	\strand (-.5,-.5) -- (.5,.5);
	\strand (-.5,.5) -- (.5,-.5);
\end{knot}
\end{tikzpicture} ~=~  \dfrac{a_3}{a_1}\cdot\left(~ \begin{tikzpicture}[baseline, rounded corners = 5mm]
    \draw (.5,.5) -- (.1,0) -- (.5,-.5);
    \draw (-.5,.5) -- (-.1,0) -- (-.5,-.5);
\end{tikzpicture} ~\pm~ \begin{tikzpicture}[baseline, rounded corners = 5mm]
    \draw (-.5,-.5) -- (0,.1) -- (.5,-.5);
    \draw (-.5,.5) -- (0,-.1) -- (.5,.5);
\end{tikzpicture} ~ \right)
\end{equation*}
Setting $z = a_3/a_1$ gives our desired formula.
\end{proof}

\subsection{Planar algebra of unoriented virtual tangles} \label{virtualtangles}
Another planar algebra that we will study is the planar algebra of unoriented virtual tangles. Virtual knot theory was first discussed in \cite{Kau99} while virtual tangles can be found in \cite{SW99}. A virtual tangle is a tangle that includes the following diagram as a possible crossing:
\begin{center}
\begin{tikzpicture}[baseline]
    \draw (-.5,-.5)-- (.5,.5);
    \draw (.5,-.5) -- (-.5,.5);
    \draw[fill=white] (0,0) circle (.1cm);
\end{tikzpicture}

\end{center}
in addition to the normal crossings. One can think of the planar algebra of tangles as a collection of planar tetravalent graphs, with crossings and boundary points serving as vertices. The inclusion of the virtual crossing merely drops the condition of planarity on these graphs. The virtual crossing is neither an overcrossing nor an undercrossing, so one could think of the strands as literally passing through each other. Alternatively, \cite{Bro16} explains that one could imagine this non-planar graph as being embedded on some punctured surface in which the given diagram could be drawn without the need for self-intersection.

\begin{definition}\label{def:symmetric}
A planar algebra, $\mathcal{V},$ is said to be \emph{symmetric} if there exists an element of $V_4$, \begin{tikzpicture}[baseline = -.15cm, scale = .65]
    \draw (-.5,-.5)-- (.5,.5);
    \draw (.5,-.5) -- (-.5,.5);
    \draw[fill=white] (0,0) circle (.1cm);
\end{tikzpicture}
, that satisfies the following relations:
\begin{center}
\begin{tabular}{ccccccccc}
$\begin{tikzpicture}[baseline=.5cm, rounded corners = 5mm]
\draw (0,1)--(1,0) -- (0,0)--(1,1);
\draw[fill = white] (.5,.5) circle (.1cm);
\end{tikzpicture} = \begin{tikzpicture}[baseline = .5 cm, rounded corners = 6mm]
    \draw (0,1)--(.5,-.1) -- (1,1);
\end{tikzpicture}$
&&
$\begin{tikzpicture}[baseline={([yshift=-.5ex]current bounding box.center)}, rounded corners = 5mm, scale = .7]
\draw[name path=line 1] (-.5,-.5) -- (.5,.5);
\draw[name path=line 2] (.5,-.5) -- (-.5,.5);
\draw[name path=line 3] (-.5,.5) -- (.5,1.5);
\draw[name path=line 4] (.5,.5) -- (-.5,1.5);
	\draw[fill=white,name intersections={of=line 1 and line 2}]
    (intersection-1) circle (.095cm);
	\draw[fill=white,name intersections={of=line 3 and line 4}]
    (intersection-1) circle (.095cm);
\end{tikzpicture} = \begin{tikzpicture}[baseline, rounded corners = 5mm]
    \draw (.5,.5) -- (.1,0) -- (.5,-.5);
    \draw (-.5,.5) -- (-.1,0) -- (-.5,-.5);
\end{tikzpicture}$
&&
$\begin{tikzpicture}[baseline={([yshift=-.5ex]current bounding box.center)}, rounded corners = 5mm]
    \draw[name path=line 1] (-.45,-.5) -- (0,.3) -- (.45,.5);
    \draw[name path=line 2] (0,.5) -- (.3,0) -- (0,-.5);
    \draw[name path=line 3] (.45,-.5) -- (-.6,0) -- (-.45,.5);
    \draw[fill=white,name intersections={of=line 1 and line 2}]
    (intersection-1) circle (.095cm);
    \draw[fill=white,name intersections={of=line 1 and line 3}]
    (intersection-1) circle (.095cm);
    \draw[fill=white,name intersections={of=line 2 and line 3}]
    (intersection-1) circle (.095cm);
\end{tikzpicture}
=
\begin{tikzpicture}[baseline={([yshift=-.5ex]current bounding box.center)}, rounded corners = 5mm, xscale = -1]
    \draw[name path=line 1] (-.45,-.5) -- (0,.3) -- (.45,.5);
    \draw[name path=line 2] (0,.5) -- (.3,0) -- (0,-.5);
    \draw[name path=line 3] (.45,-.5) -- (-.6,0) -- (-.45,.5);
    \draw[fill=white,name intersections={of=line 1 and line 2}]
    (intersection-1) circle (.095cm);
    \draw[fill=white,name intersections={of=line 1 and line 3}]
    (intersection-1) circle (.095cm);
    \draw[fill=white,name intersections={of=line 2 and line 3}]
    (intersection-1) circle (.095cm);
\end{tikzpicture} $ \\
Virtual Reidemeister I && Virtual Reidemeister II && Virtual Reidemeister III
\end{tabular}
\end{center}
In addition, the \begin{tikzpicture}[baseline = -.15cm, scale = .65]
    \draw (-.5,-.5)-- (.5,.5);
    \draw (.5,-.5) -- (-.5,.5);
    \draw[fill=white] (0,0) circle (.1cm);
\end{tikzpicture}
 must satisfy the following naturality condition:
\begin{equation*}
\includegraphics[valign = c,scale = 1.25]{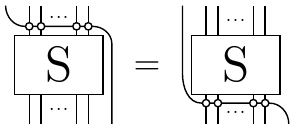}
\end{equation*}
for any element $S\in \mathcal{V}.$
\end{definition}

We will abbrieviate the virtual Reidemeister relations as vR1, vR2, and vR3, respectively. A classic example of a symmetric planar algebra generated by $\begin{tikzpicture}[baseline = -.15cm, scale = .65]
    \draw (-.5,-.5)-- (.5,.5);
    \draw (.5,-.5) -- (-.5,.5);
    \draw[fill=white] (0,0) circle (.1cm);
\end{tikzpicture}
$ is Deligne's $O_t$, which is also known as the Brauer category. One can find more information about the Brauer category in \cite{FY89} and \cite{RT90} (Brauer originally defined the algebra in \cite{Bra37}.).

\begin{example}[Deligne's $O_t$]\label{ex:ot}
Deligne's $O_t$ (found in \cite{EGNO15}(Section 12.9.3), \cite{Del90}, and \cite{Del02}) is a symmetric planar algebra generated by \begin{tikzpicture}[baseline = -.15cm, scale = .65]
    \draw (-.5,-.5)-- (.5,.5);
    \draw (.5,-.5) -- (-.5,.5);
    \draw[fill=white] (0,0) circle (.1cm);
\end{tikzpicture}
 with modulus $t\in\mathbb{C}^\times$.

A priori, it is not clear that the above relations are enough to reduce every closed diagram to a multiple of the empty diagram. This fact is proven in \cite{RT90}, however. Because $O_t$ is generated by just a symmetric crossing, we also get the following universal property:

\begin{proposition}\label{prop:universal}
Let $\mathcal{V}$ be a symmetric planar algebra with modulus $t$. Then there exists a map
\begin{equation*}
    O_t \hookrightarrow \mathcal{V}
\end{equation*}
\end{proposition}

\begin{proof}
Let $\mathcal{V}$ be as described above, and $X$ be the element of $V_4$ satisfying the virtual Reidemeister relations and naturality conditions. Then the map $\begin{tikzpicture}[baseline = -.15cm, scale = .65]
    \draw (-.5,-.5)-- (.5,.5);
    \draw (.5,-.5) -- (-.5,.5);
    \draw[fill=white] (0,0) circle (.1cm);
\end{tikzpicture}
\mapsto X$ is a map of planar algebras because $X$ satisfies the virtual Reidemeister relations and $\mathcal{V}$ has modulus $t.$ Since this map is clearly injective, there must exist a map $O_t \hookrightarrow \mathcal{V},$ as desired.
\end{proof}
\end{example}

\begin{definition}
The planar algebra of unoriented virtual tangles is the planar algebra generated by 
\[ \begin{tikzpicture}[baseline]
\begin{knot}[clip width = 4]
	\strand (-.5,-.5) -- (.5,.5);
	\strand (-.5,.5) -- (.5,-.5);
\end{knot}
\end{tikzpicture} \hspace{.5cm} \text{and} \hspace{.5cm} \begin{tikzpicture}[baseline]
    \draw (-.5,-.5)-- (.5,.5);
    \draw (.5,-.5) -- (-.5,.5);
    \draw[fill=white] (0,0) circle (.1cm);
\end{tikzpicture}
 \]
Diagrams are equivalent up to isotopy, a finite number of Reidemeister and virtual Reidemeister II and III moves, inaddition to the following relations:
\begin{center}
\begin{tabular}{ccc}
$\begin{tikzpicture}[baseline=.5cm, rounded corners = 5mm]
\draw (0,1)--(1,0) -- (0,0)--(1,1);
\draw[fill = white] (.5,.5) circle (.1cm);
\end{tikzpicture} = \pm\begin{tikzpicture}[baseline = .5 cm, rounded corners = 6mm]
    \draw (0,1)--(.5,-.1) -- (1,1);
\end{tikzpicture}$
&&
$\begin{tikzpicture}[scale = .5, baseline,rounded corners = 5mm]
	\draw[white, name path=line 1] (-1,-1)-- (1,1);
    \draw[white, name path=line 2] (1,-1) -- (-1,1);
    \draw[white, name path=line 3] (0,-1) -- (1,0) -- (0,1);
    \draw (-1,-1)-- (-.1,-.1);
    \draw (.1,.1) -- (1,1);
    \draw (1,-1) -- (-1,1);
    \draw (0,-1) -- (1,0) -- (0,1);
	\draw[fill=white,name intersections={of=line 1 and line 3}]
    (intersection-1) circle (.15cm);
    \draw[fill=white,name intersections={of=line 2 and line 3}]
    (intersection-1) circle (.15cm);
\end{tikzpicture}
=
\begin{tikzpicture}[scale = .5, baseline,rounded corners = 5mm, xscale = -1]
	\draw[white, name path=line 1] (-1,-1)-- (1,1);
    \draw[white, name path=line 2] (1,-1) -- (-1,1);
    \draw[white, name path=line 3] (0,-1) -- (1,0) -- (0,1);
    \draw (1,-1)-- (.1,-.1);
    \draw (-.1,.1) -- (-1,1);
    \draw (-1,-1) -- (1,1);
    \draw (0,-1) -- (1,0) -- (0,1);
	\draw[fill=white,name intersections={of=line 1 and line 3}]
    (intersection-1) circle (.15cm);
    \draw[fill=white,name intersections={of=line 2 and line 3}]
    (intersection-1) circle (.15cm);
\end{tikzpicture}$\\
Virtual Twist Relation && Mixed Reidemeister III
\end{tabular}
\end{center}
\end{definition}

Note that a priori we do not require vR1 to be satisfied on the nose. By multiplying the virtual crossing by itself and then capping, however, we see that vR1 must be satisfied up to a sign. Moreover, we have the following fact:

\begin{lemma}\label{lemma:vtwist}
The planar algebra of unoriented virtual tangles has an automorphism generated by
\begin{equation*}
\begin{tikzpicture}[baseline]
    \draw (-.5,-.5)-- (.5,.5);
    \draw (.5,-.5) -- (-.5,.5);
    \draw[fill=white] (0,0) circle (.1cm);
\end{tikzpicture}
 ~\mapsto~ -\begin{tikzpicture}[baseline]
    \draw (-.5,-.5)-- (.5,.5);
    \draw (.5,-.5) -- (-.5,.5);
    \draw[fill=white] (0,0) circle (.1cm);
\end{tikzpicture}

\end{equation*}
\end{lemma}

\begin{proof}
Let $\Phi$ be the map of planar algebras generated by sending the virtual crossing to its negative. It is clear that this map is a map of planar algebras and that it is injective. Considering the map is also an idempotent, we see that the inverse map is also injective and thus the map is a planar algebra isomorphism.
\end{proof}

Thus, any skein theory in which the virtual twist parameter is $+1$ can be transformed into a skein theory in which it is $-1$ by simply replacing the virtual crossing with its negative in every relation. As such, we will choose for the parameter to be $+1$ in what follows. Thus, the planar algebra of unoriented virtual tangles is symmetric.
\begin{definition}\label{trivial}
A quotient of the planar algebra of unoriented virtual tangles is said to give a \emph{fully-flat} invariant if the following relation holds:
\begin{equation*}
\begin{tikzpicture}[baseline]
\begin{knot}[clip width = 4]
	\strand (-.5,.5) -- (.5,-.5);
	\strand (-.5,-.5) -- (.5,.5);
\end{knot}
\end{tikzpicture}
 ~=~ \pm\begin{tikzpicture}[baseline]
    \draw (-.5,-.5)-- (.5,.5);
    \draw (.5,-.5) -- (-.5,.5);
    \draw[fill=white] (0,0) circle (.1cm);
\end{tikzpicture}

\end{equation*}

\end{definition}
\begin{remark}
A planar algebra is \emph{flat} if the overcrossing and undercrossing are equal to each other. A planar algebra is fully-flat if the two crossings are equal \emph{and} they are both equal to the virtual crossing up to a sign. Thus, fully-flatness is a stronger condition than flatness. We will see examples of planar algebras that are flat yet not fully-flat (The virtual TLJ planar algebra in Example \ref{ex:vjones} at $A = \pm 1$ for instance.). Quotients of the planar algebra of unoriented virtual tangles that are fully-flat are already well-understood and are omitted from this classification. 
\end{remark}

In addition to the above relations, there is also an important non-relation of note which appears later in this paper that is aptly named the ``forbidden move:"
\begin{equation}\label{eq:forbidden}
\begin{tikzpicture}[scale = .5, baseline,rounded corners = 5mm]
\begin{knot}[clip width = 3]
    \strand (-1,-1)-- (1,1);
    \strand (1,-1) -- (-1,1);
    \strand (0,-1) -- (1,0) -- (0,1);
\end{knot}
	\draw[fill = white]  (0,0) circle [radius = .15cm];
\end{tikzpicture}
~\ne~
\begin{tikzpicture}[scale = .5, baseline,rounded corners = 5mm, xscale = -1]
\begin{knot}[clip width = 3]
    \strand (-1,-1)-- (1,1);
    \strand (1,-1) -- (-1,1);
    \strand (0,-1) -- (1,0) -- (0,1);
\end{knot}
	\draw[fill = white]  (0,0) circle [radius = .15cm];
\end{tikzpicture}
\end{equation}

The forbidden move implies that the virtual crossing is not natural with the actual crossing in the planar algebra of unoriented virtual tangles. In tensor category terms, \cite{Bro16} explains that although there is a unique symmetric monoidal functor from the category of tangles to the category of symmetric tangles, that functor is not braided. We explain this phenomenon in planar algebra terms in the following theorem:
\begin{theorem}[\cite{Bro16}]\label{thm:symmetric}
Let $\mathcal{T}_{virt}$ be the planar algebra of unoriented virtual tangles, $\mathcal{P}$ be a symmetric planar algebra, and $\mathcal{A}$ be a braided sub-planar algebra of $\mathcal{P}$. Then there exists a canonical map from
\[ \mathcal{T}_{virt} \to \mathcal{P} \]
\end{theorem}
\subsubsection{$S_{2n}$ actions on elements generated by a virtual crossing}

Suppose that $\mathcal{V}$ is a quotient of the planar algebra of virtual tangles that is generated by the virtual crossing. Then $S_{2n}$ acts on $V_{2n}$ by permutation of the endpoints. Two strings crossing is represented by $\begin{tikzpicture}[baseline = -.15cm, scale = .65]
    \draw (-.5,-.5)-- (.5,.5);
    \draw (.5,-.5) -- (-.5,.5);
    \draw[fill=white] (0,0) circle (.1cm);
\end{tikzpicture}
\!$. Moreover, this action forms an $S_{2n}$ $\mathbb{C}$-representation. An important example that will come up in Section \ref{virtualskein} is the following:

\begin{proposition}\label{prop:s2n} Let $\mathcal{V}$ be a quotient of the planar algebra of virtual tangles generated by $\begin{tikzpicture}[baseline = -.15cm, scale = .65]
    \draw (-.5,-.5)-- (.5,.5);
    \draw (.5,-.5) -- (-.5,.5);
    \draw[fill=white] (0,0) circle (.1cm);
\end{tikzpicture}
\!$. Then the permutation action of $S_6$ on the diagrams of $V_6$ forms a 15-dimensional $\mathbb{C}$-representation of $S_6$. Moreover, this representation is a direct sum of the trivial representation 
a 5-dimensional representation, and a 9-dimensional representation.
\end{proposition}

\begin{proof}
Any quotient of the planar algebra of virtual tangles generated by the virtual crossing has that $\dim V_6 \leq 15$. Thus, we obtain the following spanning set:
\begin{center}
\includegraphics[valign = c]{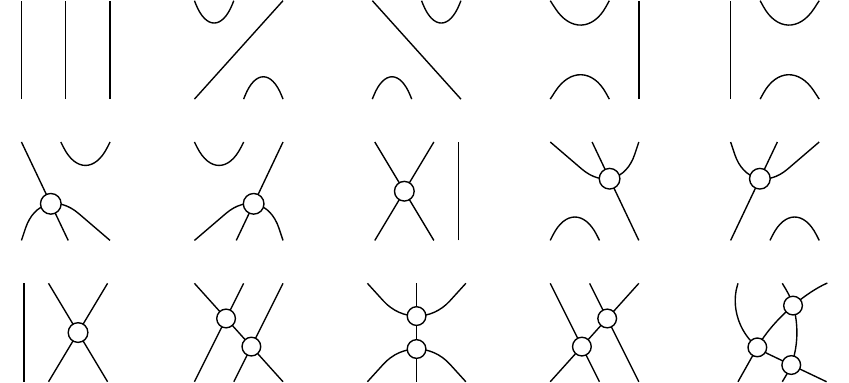}
\end{center}
It is clear the permuting the vertices of these diagrams gives an action of $S_6$ on $V_6$ and that the resulting representation of $S_6$ is 15-dimensional. By examining the characters of this representation, we take the inner product of this representation with each of the irreducible representations of $S_6$. From the first orthogonality theorem, we see that the 15-dimensional representation is reducible as a direct sum of three irreducible representations: the trivial representation, a 5-dimensional representation, and a 9-dimensional representation corresponding to the following Young diagrams:
\begin{center}
\begin{tabular}{ccccc}
\yng(6) && \yng(2,2,2) && \yng(4,2) \\
\end{tabular}
\end{center}
which completes our proof.
\end{proof}

Knowing how this action decomposes into irreducible representations is extremely helpful in determine the possible dimensions of a box space. Let $\mathcal{V}$ be a quotient of the planar algebra of virtual tangles that is generated by the virtual crossing and $\mathcal{W}$ be a quotient of $\mathcal{V}$. Thus, $S_{2n}$ acts on $W_{2n}$ by a corresponding quotient representation. Since $\mathcal{V}$ and $\mathcal{W}$ are semi-simple, any element appearing in the kernel of a surjective map between them must be equal to zero. In other words, the sub-space of relations that hold in $\mathcal{W}$ but not $\mathcal{V}$ is given by the kernel of a surjective map between them. Because $\mathcal{V}$ is semi-simple, $\mathcal{W}$ is also a sub-planar algebra of $\mathcal{V}$.  Moreover, the $S_{2n}$ action on $W$ must be a sub-representation of the $S_{2n}$ action on $\mathcal{V}$. Thus, if we can determine the invariant subspace of an irreducible component of the $S_{2n}$ action on $\mathcal{V}$, then we have found a space of possible relations for some quotient $\mathcal{W}$. If we know how the latter representation decomposes as irreducibles, then we know what the possible representations for the action on $W_{2n}$ are and thus what the possible box-space dimensions are. 

\begin{corollary}\label{cor:boxdim}
Let $\mathcal{V}$ be a quotient of the planar algebra of virtual tangles generated by the virtual crossing and suppose that $\dim V_4 = 3$. Then $\dim V_6 \in \{10, 15\}$.
\end{corollary}

\begin{proof}
By Proposition \ref{prop:s2n}, we know that $S_6$ acts on $V_6$ and forms a 15-dimensional representation. Thus, $S_6$ must act on any quotient of $V_6$. Moreover, the corresponding representation of $S_6$ must be a quotient of that 15-dimensional representation. Since we have required that $\dim V_4 = 3$, this implies that the cupcap, identity, and virtual crossing form a basis for $V_4$.

By inspection, the trivial representation corresponds to the subspace of $V_6$ generated by the sum of the 15 elements of $V_6$ shown above. Thus, if the corresponding quotient representation of some sub-planar algebra did not include the trivial representation, this element would be equal to zero. Capping the resulting equation on the left, though, yields 
\begin{equation}\label{eq:d14}
    4 \cdot \begin{tikzpicture}[baseline, rounded corners = 5mm]
    \draw (-.5,-.5) -- (0,.1) -- (.5,-.5);
    \draw (-.5,.5) -- (0,-.1) -- (.5,.5);
\end{tikzpicture} + (d + 3) \begin{tikzpicture}[baseline, rounded corners = 5mm]
    \draw (.5,.5) -- (.1,0) -- (.5,-.5);
    \draw (-.5,.5) -- (-.1,0) -- (-.5,-.5);
\end{tikzpicture} + (d + 2) \begin{tikzpicture}[baseline]
    \draw (-.5,-.5)-- (.5,.5);
    \draw (.5,-.5) -- (-.5,.5);
    \draw[fill=white] (0,0) circle (.1cm);
\end{tikzpicture}
 = 0
\end{equation}
Thus, since for all $d \in \mathbb{C}$ the cupcap, identity, and virtual crossing cannot simultaneously vanish, they must be linearly dependent, a contradiction. Thus, any quotient of $\mathcal{V}$ with $\dim V_4 = 3$ must have $\dim V_6 \in \{1,6, 10,15\}$. When $\dim V_6 = 1$ this implies that every diagram is a multiple of the identity, which implies that $\dim V_4 = 1$, a contradiction. 

Suppose that $\dim V_6 = 6$. Since the cupcap and identity are linearly independent, this implies a relation of the following form:
\begin{equation*}
    \includegraphics{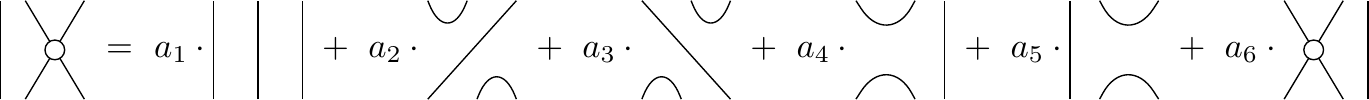}
\end{equation*}
By capping this off on the left, we obtain the equation:
\begin{equation*}
d \cdot \begin{tikzpicture}[baseline]
    \draw (-.5,-.5)-- (.5,.5);
    \draw (.5,-.5) -- (-.5,.5);
    \draw[fill=white] (0,0) circle (.1cm);
\end{tikzpicture}
 = (da_1 + a_3 + a_6) \cdot \begin{tikzpicture}[baseline, rounded corners = 5mm]
    \draw (.5,.5) -- (.1,0) -- (.5,-.5);
    \draw (-.5,.5) -- (-.1,0) -- (-.5,-.5);
\end{tikzpicture} + (da_2+a_4 + a_5) \cdot \begin{tikzpicture}[baseline, rounded corners = 5mm]
    \draw (-.5,-.5) -- (0,.1) -- (.5,-.5);
    \draw (-.5,.5) -- (0,-.1) -- (.5,.5);
\end{tikzpicture}     
\end{equation*}
Since $d \ne 0$, this implies $\dim V_4 < 3$, a contradiction. Thus, $\dim V_6 \in \{10,15\}$, as desired.
\end{proof}

We can also imagine elements of $S_{2n}$ as diagrams in this planar algebra:
\begin{definition}\label{def:sigman}
Let $\sigma$ be an element of $S_n.$ Then define 
\begin{equation*}
     \large
\begin{tikzpicture}[scale=.3, baseline = .2cm]

\draw (.5,-1) -- (.5,3);
\draw (.9,-1) -- (.9,3);
\node at (1.5,2.5) {\tiny ...};
\node at (1.5,-.5) {\tiny ...};
\draw (2.1,-1) -- (2.1,3);
\draw (2.5,-1) -- (2.5,3);

\draw[fill = white] (0,0) rectangle (3,2);
\node at (1.5,1) {$\sigma$};

\end{tikzpicture}     
\end{equation*}
to be the diagram (in standard form) obtained by connecting the $i$-th point on the bottom of the diagram to the $\sigma(i)$-th point on the top of the diagram. Any point at which strands must cross each other is represented by $\begin{tikzpicture}[baseline = -.15cm, scale = .65]
    \draw (-.5,-.5)-- (.5,.5);
    \draw (.5,-.5) -- (-.5,.5);
    \draw[fill=white] (0,0) circle (.1cm);
\end{tikzpicture}
\!$.  
\end{definition}

As an example of this when $\sigma = (123)$ for $\sigma\in S_3$ we have the diagram  
\begin{center}
    \begin{tikzpicture}[baseline={([yshift=-.5ex]current bounding box.center)}, rounded corners = 5mm, yscale = -1]
	\draw[name path=line 1] (.45,-.5) -- (-.05,.5);
    \draw[name path=line 2] (-.45,.5) -- (.05,-.5);
    \draw[name path=line 3] (-.45,-.5) -- (.45,.5);
    \draw[fill=white,name intersections={of=line 1 and line 3}]
    (intersection-1) circle (.095cm);
    \draw[fill=white,name intersections={of=line 2 and line 3}]
    (intersection-1) circle (.095cm);
\end{tikzpicture}\!.   
\end{center}
\begin{definition}\label{def:wedge}
Let $\mathcal{V}$ be the planar algebra of unoriented virtual tangles. We define the element $\wedge^{m+1}$ to be
\begin{equation*}
\includegraphics[valign = c]{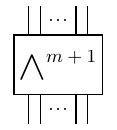} ~:=~ \sum_{\sigma\in S_{m+1}} \dfrac{\text{sgn}(\sigma)}{(m+1)!}\cdot\large
\begin{tikzpicture}[scale=.3, baseline = .2cm]

\draw (.5,-1) -- (.5,3);
\draw (.9,-1) -- (.9,3);
\node at (1.5,2.5) {\tiny ...};
\node at (1.5,-.5) {\tiny ...};
\draw (2.1,-1) -- (2.1,3);
\draw (2.5,-1) -- (2.5,3);

\draw[fill = white] (0,0) rectangle (3,2);
\node at (1.5,1) {$\sigma$};

\end{tikzpicture}
\end{equation*}
\end{definition}

As the name suggests, this diagram is related to the exterior product, as arbitrary composition of $\wedge^{m+1}$ with the virtual crossing is the negative of $\wedge^{m+1}$. If we consider these elements to be maps of representations (see Section \ref{repg} for more information on this), this implies that $\wedge^{m+1}(x_1\otimes x_2\otimes\dots\otimes x_{m + 1}) = \text{sgn}(\sigma)\cdot \wedge^{m+1}(x_{\sigma(1)}\otimes x_{\sigma(2)}\otimes\dots\otimes x_{\sigma(m + 1)}),$ for any permutation of the indices, $\sigma.$  

With this definition in hand, we are now ready to state the following result, which is a consequence of the first fundamental theorem for orthogonal groups due to \cite{Wey39}. (It is also a consequence of the second fundamental theorem due to \cite{LeZ11}).

\begin{theorem}\label{thm:wedge}
Let $\mathcal{V}$ be a quotient of the planar algebra of unoriented virtual tangles. Then $\mathcal{V} \cong \Rep(O(m))$ for $m,n\in\mathbb{Z}$ if
\begin{enumerate}[label = \roman{enumi}.]

\item $\skeincircle ~=~ m$ and

\item $\includegraphics[valign = c]{diagrams/pdf/wedgediagram.pdf} ~=~ 0$,

\end{enumerate}
where $\wedge^{m + 1}$ is as defined in Definition \ref{def:wedge}. Furthermore, these relations generate all other relations. 
\end{theorem}

Thus, in order to identify a quotient of the planar algebra of unoriented virtual tangles as $\Rep(O(t))$ at $t\in\mathbb{Z}_{> 0}$ we need only check that the planar algebra has $d = t$ and that $\wedge^{t+1} = 0.$ The last statement in this theorem is also quite powerful, as any relation between any diagrams in any box space can be obtained from the above two relations.

One well-known quotient of the planar algebra of unoriented virtual tangles is the virtual TLJ planar algebra (originally in \cite{Kau99}).

\begin{example}[The virtual TLJ planar algebra planar algebra]\label{ex:vjones}
The virtual TLJ planar algebra is a quotient of the planar algebra of unoriented virtual tangles subject to the virtual Reidemeister moves and the following relations:
\begin{center}
\begin{tabular}{ccc}
$\skeincircle ~=~ -(A^2+A^{-2})$ && $\begin{tikzpicture}[baseline=.5cm, rounded corners = 5mm]
\draw (0,1)--(.4,.6);
\draw (.6,.4)--(1,0)--(0,0)--(1,1);
\end{tikzpicture} = -A^{-3}\begin{tikzpicture}[baseline = .5 cm, rounded corners = 6mm]
    \draw (0,1)--(.5,-.1) -- (1,1);
\end{tikzpicture}$
\end{tabular}
\end{center}
\begin{equation*}
\begin{tikzpicture}[baseline]
\begin{knot}[clip width = 4]
	\strand (-.5,.5) -- (.5,-.5);
	\strand (-.5,-.5) -- (.5,.5);
\end{knot}
\end{tikzpicture}
 = A\begin{tikzpicture}[baseline, rounded corners = 5mm]
    \draw (.5,.5) -- (.1,0) -- (.5,-.5);
    \draw (-.5,.5) -- (-.1,0) -- (-.5,-.5);
\end{tikzpicture} + A^{-1}\begin{tikzpicture}[baseline, rounded corners = 5mm]
    \draw (-.5,-.5) -- (0,.1) -- (.5,-.5);
    \draw (-.5,.5) -- (0,-.1) -- (.5,.5);
\end{tikzpicture}
\end{equation*}
\end{example}

Note that although this skein theory seems identical to the actual TLJ planar algebra, we now have the virtual crossing as a generator. To see that the planar algebra is evaluable (and thus gives a virtual link invariant), we use the following evaluation algorithm: For any diagram of $V_0$, we start by using the crossing relation to remove any crossings from the diagram. This process leaves the diagram as a sum of diagrams with only virtual crossings. From the note in Example \ref{ex:ot}, we see it is then possible to reduce any such diagram to a multiple of the empty diagram.

The similarities in the evaluation algorithms of $O_t$ and the virtual TLJ planar algebra are not coincidental. It is well known that the non-virtual TLJ planar algebra is closely related to the quantum group $SU(2)_q$. Moreover, by Theorem \ref{thm:symmetric} we have an inclusion map 
\[ 
SU(2)_q \hookrightarrow O_t
\]
when $t = q^2 + 1 +q^{-2}.$ This relationship is intimately related to the underlying groups themselves.

\subsection{Quotients of the planar algebra of oriented virtual tangles}

\begin{definition}
The planar algebra of oriented virtual tangles is the planar algebra generated by 
\begin{center}
   \begin{tabular}{ccccc}
 \overcrossinga && \undercrossinga && \virtuala
\end{tabular} 
\end{center}
In addition, two diagrams will be equivalent if one can be obtained from the other via the oriented virtual Reidemeister moves. The action of the oriented planar operad is the usual insertion action.
\end{definition}

Here the oriented virtual Reidemeister relations are the relations obtained from the usual virtual Reidemeister relations with all possible orientations of the strands. Note that the generators of the planar algebra of oriented virtual tangles have two positively oriented endpoints and two negatively oreinted endpoints. Thus, the number of positively and negatively oriented endpoints must always be the same. Two examples of quotients of the planar algebra of oriented virtual tangles are $\Rep(S^1,a)$ and Deligne's $GL_t$.

\begin{example}\label{def:s1}
Consider $\Rep(S^1,a)$, where $S^1$ is the circle group and $a\in\mathbb{C}^\times$ represents the formula for the braiding below. More information about this planar algebra can be found in \cite{DGNO10}. It is given by the following relations:
\begin{center}
\includegraphics[valign=c]{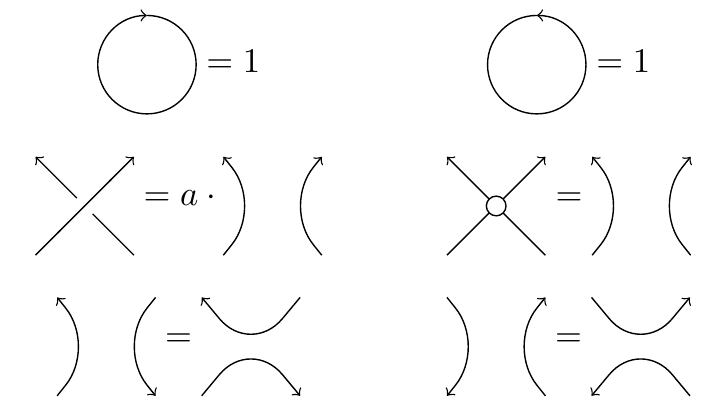}
\end{center}
The reason this planar algebra includes $\Rep(S^1)$ is detailed in Section \ref{repg}.
\end{example}

We leave it to the reader to check that the crossing satisfies the Reidemeister relations. Another example of a quotient of the planar algebra of oriented virtual tangles is Deligne's $GL_t,$ found in \cite{Del02}:

\begin{example}[Deligne's $GL_t$]\label{ex:glt}
This planar algebra is generated by the oriented symmetric crossing and has modulus $t$. Because it is generated by the oriented symmetric crossing, every diagram must have the same number of positive and negative endpoints. The planar algebra $GL_t^a$ represents the planar algebra $GL_t$ equipped with the braiding
\begin{equation*}
    \overcrossinga ~=~ a\cdot\virtuala
\end{equation*}
for some $a \in\mathbb{C}^\times$. Unlike in the unoriented case, one does not obtain one crossing from the other by rotating, so the planar algebra is not necessarily fully-flat. Note that at $t = \pm 1$, $GL_t^a$ is degenerate for all $a.$ At $t = 1$, its non-degenerate quotient is $\Rep(GL(1),a) \simeq \Rep(S^1,a).$ 
\end{example}

\subsubsection{Planar algebras of the category $\Rep(G)$}\label{repg}

One class of planar algebras that is oriented in general is the one arising from the category $\Rep(G).$ We have already seen such an example in Example \ref{def:s1}.

\begin{definition}
Let $G$ be a group. The category $\Rep(G)$ is a tensor category where $\otimes$ is defined as the tensor product of representations and $\mathbbm{1}$ is the trivial representation of $G$. Objects are tensor-generated by the irreducible, finite-dimensional representations of $G$ and morphisms are maps of representations.       
\end{definition}

For more information on tensor categories, refer to \cite{EGNO15}. We define the \emph{planar algebra} $\Rep(G)$ as follows:

\begin{definition}\label{def:repg}
Let $G$ be a group and $X$ be a finite dimensional representation of $G$. Define $X^+ = X$ and $X^- = X^*.$ Then the oriented planar algebra $\Rep(G,X)$ is a collection of vector spaces $\Rep(G,X)_{(z_1,\dots,z_i)} := \text{Hom}(\mathbbm{1}, X^{z_1}\otimes\dots\otimes X^{z_i}).$
\end{definition}

The reader might rightly object to our definition, as we are choosing a particular representation. If the representation is faithful and irreducible, however, every irreducible representation will occur as a \emph{summand} of some tensor power of $X$ and $X^*.$ Thus, every representation will appear as a \emph{projection} in some box space of $\mathcal{V}$, so we do not ``lose" much information by only considering maps involving the tensor powers of $X$ and $X^*$. When $G$ has a standard representation, we will often omit $X$ and write ``$\Rep(G)$." 

Since we have chosen to work over $\mathbb{C}$, Schur's lemma guarantees for any simple objects $X_i$ and $X_j$ that
\[\text{Hom}\left(X_i, X_j\right) \cong
\left\{\begin{array}{cl}
c\cdot id_{X_i} & i = j \text{ and } c\in\mathbb{C} \\
0 & \text{else}
\end{array}\right.\]

Thus, for some fixed value $c$, we indicate $c\cdot id_{X}$ by an oriented vertical strand and the $c^{-1} \cdot id_{X^\ast}$ by a vertical strand with the opposite orientation. If the representation is self-dual then we omit the orientation of the strands. Assume for brevity that $X$ is self-dual. The pairing and co-pairing of $X$ and $X^\ast$ can then be represented by a cap and cup, respectively. Thus, it is easy to check that the value of the circle, the map $d: \mathbb{C} \to X\otimes X \to \mathbb{C}$ is simply multiplication by the dimension of $X$. The tensor product of maps of representations is simply concatenation in the planar algebra. Also note that because $\Rep(G,X)$ is a \emph{symmetric} tensor category (i.e. $X\otimes Y \cong Y\otimes X$ for all $X$ and $Y$), we also have a map 
\begin{equation*}
\begin{tikzpicture}[baseline]
    \draw (-.5,-.5)-- (.5,.5);
    \draw (.5,-.5) -- (-.5,.5);
    \draw[fill=white] (0,0) circle (.1cm);
\end{tikzpicture}
: X\otimes X \to X\otimes X  
\end{equation*}
given by $x_1 \otimes x_2 \mapsto x_2 \otimes x_1.$ The reader can check that this map is a map of representations.

In order to give a presentation of $\Rep(G,X)$, we need ``enough" maps to tensor-generate every other map. That is, we need a set of maps such that any map from $\mathbbm{1} \to X^{\otimes n}$ can be obtained by arbitrary composition of those maps. For example, $\Rep(O(2))$ is generated by $\begin{tikzpicture}[baseline = -.15cm, scale = .65]
    \draw (-.5,-.5)-- (.5,.5);
    \draw (.5,-.5) -- (-.5,.5);
    \draw[fill=white] (0,0) circle (.1cm);
\end{tikzpicture}
\!$.

Now that we have defined the planar algebra $\Rep(G)$, we can prove the following corollary to Proposition \ref{prop:universal}:
\begin{corollary}\label{otmaps}
Let $G$ be a group and $X$ be a faithful, irreducible, self-dual representation with $\dim X = t$. Then $\Rep(G,X)$ is a symmetric planar algebra.
\end{corollary}
\begin{proof}
Since $\Rep(G,X)$ has a symmetric crossing that satisfies the naturality conditions of Definition \ref{def:symmetric}, this proposition is a consequence of Proposition \ref{prop:universal}.
\end{proof}

\subsection{Examples of shaded planar algebras}

\begin{example}[$\text{GPA}(\Gamma)$] 
 The planar algebra $\text{GPA}(\Gamma)$, described in Section \ref{starn} is a shaded planar algebra.
\end{example}

Recall that for $\text{GPA}(\ast_n)$,  $\text{GPA}(\ast_n)_{\pm 2i}$ is isomorphic as a vector space to linear combinations of Kronecker deltas of loops of length $i$. If we fix an indexing of the vertices of $S$, then a two-box in $\text{GPA}(\ast_n)_{+2}$ is a complex matrix where the $(i,j)$-th entry is the the coefficient of $\delta_{ij}$. Suppose $\mathcal{V}$ is a planar algebra generated by a 2-box, $P$. Recall that a spin model is a map from $\mathcal{V}$ to $\text{GPA}(\ast_n)$. Then giving a spin model for $\mathcal{V}$ is the same as specifying a complex matrix as the image of $P$. This leads to the following definition:

\begin{definition}
Let $\mathcal{V}$ be a shaded planar algebra generated by a 2-box, $P$. Then a graph $\Gamma$ gives a spin model for $\mathcal{V}$ if the image of $P$ in the spin model is the adjacency matrix for $\Gamma$. 
\end{definition}
See \cite{Edg18} for more information on what types of graphs can appear for certain planar algebras.

\subsubsection{Group/Subgroup planar algebra.}\label{connection}

An important sub-planar algebra of $\GPA(\ast_n)$ is the group/subgroup planar algebra. Before we define it we make the following definition based on the work of \cite{Gup08}:
\begin{definition}
Let $H\leq G$ be finite groups with $[G:H] = m$ and coset representatives $\{g_i H\}_{i=1\dots,m}.$ The diagonal action of $G$ on a $k$-tuple of cosets is defined by 
\begin{equation*}
 g\cdot(g_{i_1}H \cdots g_{i_k}H) =  (g g_{i_1})H \cdots (g g_{i_k})H 
\end{equation*}
\end{definition}
\begin{definition}
Let $H\leq G$ be finite groups with $[G:H] = m$. Suppose we have $\GPA(\ast_m)$ where the even vertices are labelled by coset representatives $S = \{g_i H\}_{i=1\dots,m}.$ Then the group-subgroup planar algebra, $\text{PA}(H \leq G),$ is the sub-planar algebra of $\GPA(\ast_m)$ that is invariant under the diagonal action of $G.$ That is, for any linear functional $f\in V_k$
\begin{equation*}
f(g_{i_1}H \cdots g_{i_k}H) =  f(g(g_{i_1}H \cdots g_{i_k}H)) = f((g g_{i_1})H \cdots (g g_{i_k})H )
\end{equation*}
for all $g\in G$ and $g_i H \in S$. 
\end{definition}
Because $\text{PA}(H \leq G)$ is shaded, we will refer to diagrams with $2i$ boundary points as $i$-boxes and $\text{PA}(H \leq G)^\text{even}$ as the collection of vector spaces $\{\text{PA}(H \leq G)_{+2i}\}_i = 0,1,\dots.$
\begin{theorem}[Theorem 5.12 in \cite{Gup08}]
Given a finite group $G$, a subgroup $H$ such that $[G:H] = n$, and
an outer action $\alpha$ of $G$ on the hyperfinite $II_1$-factor $R$, the planar algebra of the subgroup-subfactor $R \rtimes H \subset R \rtimes G$ is isomorphic to the G-invariant planar sub-algebra of $\GPA(\ast_n).$
\end{theorem}
\noindent Standard translation between subfactors and tensor categories gives us the following corollary:
\begin{corollary}\label{thm:trans}
$\text{PA}(H \leq G)^\text{even} \cong \Rep(G, X)$ where $X = \text{Ind}_{H}^{G}(\mathbbm{1}_H)$, the induction of the trivial representation from $H$ to $G.$  
\end{corollary}
Another way to think about $\text{Ind}_{H}^{G}(\mathbbm{1}_H)$ is that it is the permutation representation of $G$ acting on the set of left cosets of $H$ in $G$, so we will call it $\text{Perm}_H^G$ for short. Let $\mathcal{V}$ be generated by a 2-box. Suppose we have a spin model for $\Phi:\mathcal{V}\to\text{GPA}(\ast_n)$ given by a graph $\Gamma$. $\Gamma$ is \emph{transitive} if every pair of vertices is equivalent under some element of its automorphism group. Let $G = \text{Aut}(\Gamma)$ and $H$ be the subgroup of $G$ fixing a special vertex $x$. Then if $\Gamma$ is transitive the coset $g_i H$ can be thought of as the set of all automorphisms sending $x$ to some vertex $g_i(x)$ of $\Gamma$. In this way, we can identify all the vertices of $\Gamma$ to distinct cosets. Thus, we see that necessarily $[G:H] = n$ and that the action of $G$ on the left cosets of $H$ is transitive. Hence, $\text{PA}(H \leq G)$ sits inside $\GPA(\ast_n)$. With this in mind, we can prove the following theorem:

\begin{theorem}\label{thm:gup}
Let $\mathcal{V}$ be a planar algebra generated by a 2-box and $\Phi:\mathcal{V}\to\GPA(\ast_n)$ be a spin model for $\mathcal{V}$. Suppose that the matrix of weights for this map is given by the adjacency matrix of some graph $\Gamma.$ Let $G= \text{Aut}(\Gamma),$ the group of graph automorphisms of $\Gamma$, and $H$ be the subgroup of graph automorphisms fixing a chosen vertex $x$. If the graph $\Gamma$ is transitive then $\Phi(\mathcal{V}) \subseteq \text{PA}(H \leq G).$ 
\end{theorem}
\begin{proof}
Let $\mathcal{V}$ be generated by a single 2-box which we will call $X$. Suppose $\Gamma$ is transitive and gives a map $\Phi:\mathcal{V}\to\GPA(\ast_n)$. Then by the note above $\text{PA}(H \leq G) \subseteq \text{GPA}(\ast_n)$ and we will identify the vertices of $\Gamma$ with the cosets $g_i H$. To show $\Phi(\mathcal{V}) \subseteq \text{PA}(H \leq G)$ we need only prove that $\Phi(X) \in \text{PA}(H \leq G).$ Since $\displaystyle X\mapsto \sum_{i,j} \delta_{g_i H g_j H},$ the set of Kronecker deltas of edges in $\Gamma$, we see that the diagonal action of $g$ on $X$ is
\begin{equation*}
  g\sum_{i,j} \delta_{g_i H g_j H} = \sum_{i,j} g\cdot\delta_{g_i H g_j H} = \sum_{i,j} \delta_{(g g_i)H(g g_j)H} = \sum_{k,l} \delta_{g_k H g_l H} = \sum_{i,j} \delta_{g_i H g_j H}
\end{equation*}
Thus $\Phi(X)$ is invariant under the diagonal action, and so $\Phi(X)\in\text{PA}(H \leq G)$ by definition. Since $X$ generates $\mathcal{V}$, $\Phi(\mathcal{V}) \in\text{PA}(H \leq G),$ as desired. 
\end{proof}

Unfortunately, $\text{Perm}^G_H$ is never irreducible when $H \ne G$, which means that the box spaces of $\Rep(G,\text{Perm}^G_H)$ can grow very quickly. To make the planar algebra more manageable to study, we will consider a cut-down (See Definition \ref{cutdown}) of it instead.

When $P\in V_{\pm 2i}$ then $\mathcal{V}$ is an unshaded planar algebra, while it is still shaded if $P\in V_{\pm (2i + 1)}$. In this paper, we will only utilize projections in the 2-box space, and so all the resulting planar algebras will be unshaded. If $\mathcal{V}$ is an unoriented planar algebra generated by a 2-box with a spin model given by a transitive graph $\Gamma$, then we know that $\mathcal{V}$ is a sub-planar algebra of $\text{PA}(H \leq G),$ which is isomorphic to $\Rep(G,\text{Perm}^G_H)$. If we cut down $\mathcal{V}$ by an appropriate projection in the 2-box space, however, the even portion of the cut-down is then isomorphic to $\Rep(G, X),$ where $X$ is some summand of $\text{Perm}^G_H$.

Let $\Gamma$ be a transitive graph, $G = \text{Aut}(\Gamma),$ and $H$ be the automorphism group fixing a chosen point. Then in this case we can describe the decomposition of $\text{Perm}^G_H$ into irreducible representations. This exercise is equivalent to finding the projections in the 2-box space of $\text{PA}(H \leq G)$. Recall from the definition of $\text{PA}(H \leq G)$ that every element of the 2-box space is of the form
\begin{equation*}
    \sum_{g \in G } g \cdot \delta_{H g_i H}
\end{equation*}
where the $g_i$ are specified coset representatives. Thus, $\delta_{H g_i H}$ appears as a summand in every element of the 2-box space, and we can think of describing the action of $G$ on the 2-box space as sending $Hg_i H \mapsto Hg_j H.$ Further, we can think of these $H g_i H$ as representing the relationship between every point of $\Gamma$ and our chosen fixed point (i.e. adjacent, non-adjacent, or equal), by identifying $x$ with the left coset $1 \cdot H = H$. Since $G$ is the automorphism group of $\Gamma$, every element of $G$ must preserve the above relationships. Thus for all $g \in G,$ $g \cdot (HH) = HH$, and for $g_i \ne 1,$ $g \cdot (H g_i H) = H g_j H,$ where $g_j H$ is adjacent to $H$ if and only if $g_i H$ is. Thus, in this way we see that there are at most three invariant subspaces of the 2-box space of $\text{PA}(H \leq G).$ Moreover, by the transitivity of $\Gamma$ these invariant subspaces cannot contain any non-trivial, proper invariant subspaces. 

Suppose $\Gamma$ is a non-complete graph with at least one edge. Translating the above information to $Rep(G,\text{Perm}^G_H),$ using Corollary \ref{thm:trans}, we know when $\Gamma$ is transitive that
\begin{equation*}
    \text{Perm}^G_H ~=~ \mathbbm{1}_G \oplus X_k \oplus X_{n-k-1}     
\end{equation*}
with $\dim X_k = k$ and $\dim X_{n-k-1} = n-k-1$ and $k$ being the number of neighbors of our chosen vertex. The number of non-neighbors is thus $n-k-1.$ Note that transitivity implies that $\Gamma$ is regular, and hence the above decomposition is independent of the choice of vertex. In the case where $\Gamma$ is the unique graph with 0 or 1 vertex, $G = H = \{1\},$ so $\text{Perm}^G_H = \mathbbm{1}_G.$ When $\Gamma$ is a complete or empty graph of at least two vertices, we note that $\text{Perm}^G_H = \text{Perm}^{S_n}_{S_{n-1}} = \mathbbm{1}_{S_n} \oplus X_{n-1},$ where $X_{n-1}$ is the standard representation for $S_n$.

\begin{proposition}\label{prop:symmetric}
Suppose $\mathcal{V}$ is a planar algebra that has a spin model given by a transitive graph, $\Gamma$. Let $G = \text{Aut}(\Gamma)$ and $H$ be the automorphism group fixing a chosen point. Further, assume that $\mathcal{V}^\text{even}$ is braided. Then $\Rep(G,X)$ takes a map from the planar algebra of unoriented virtual tangles, where $X$ is an irreducible summand of the permutation representation of $G$ induced from $H$ appearing exactly once.
\end{proposition}

\begin{proof} 
Suppose a planar algebra, $\mathcal{V},$ has a spin model given by the adjacency matrix of some graph $\Gamma$ and suppose that $\Gamma$ is transitive. Then we know from Corollary \ref{thm:gup} that the even portion of this planar algebra is a sub-planar algebra of $\text{PA}(H \leq G),$ where $G = \text{Aut}(\Gamma)$ and $H$ is the subgroup of $G$ fixing some vertex $x.$ By Corollary \ref{thm:trans} we know that $\text{PA}(H \leq G) \simeq \Rep(G, \text{Perm}^G_H)$. Because $\Gamma$ is a transitive graph, Perm is the sum of at most three distinct representations, including $\mathbbm{1}_G.$ Let $X$ be an irreducible summand of $\text{Perm}^G_H$ appearing exactly once. Since the characters of $\text{Perm}^G_H$ are integral, $\text{Perm}^G_H$ is necessarily self-dual. Moreover, since this irreducible representation of $\text{Perm}^G_H$ appears exactly once, we know it is also self-dual.

Let $\mathcal{W}$ be the cut-down of $\mathcal{V}$ by the appropriate projection of $V_{+2}$ such that $\mathcal{W}^\text{even} \simeq \Rep(G, X)$. It is clear that $\mathcal{V}^\text{even}$ being braided implies that $\mathcal{W}^\text{even}$ is braided also. Moreover, because $X$ is self-dual, $\Rep(G,X)$ is an unoriented planar algebra, and thus a symmetric planar algebra by Corollary \ref{otmaps}. Since $\mathcal{W}^\text{even}$ is braided, we know that $\Rep(G,X)$ has a sub-braiding. Thus by Theorem \ref{thm:symmetric}, $\Rep(G,X)$ must take a map from the planar algebra of unoriented virtual tangles, as desired.
\end{proof}

Kuperberg \cite{Kup97} reformulates the classification of spin models for the Kauffman polynomial \cite{Jae95} in terms of the combinatorial $B_2$ spider. In particular, he proved that $\Gamma$ gave a spin model for the $B_2$ spider if and only if it gave a spin model for a Kauffman polynomial planar algebra. The $B_2$ spider is generated by two types of strands, which are denoted by single and doubled strands. If one takes the sub-planar algebra generated by only the doubled strands, she would obtain a Kauffman polynomial planar algebra, which is clearly braided. Moreover, it is the cut-down of the Kauffman polynomial planar algebra in \cite{Jae95} with spin model given by $\Gamma$. Since every graph that gives a spin model for the Kauffman polynomial is necessarily transitive, the even portion of this sub-planar algebra of the $B_2$ spider is a sub-planar algebra of $\Rep(\text{Aut}(\Gamma),X),$ where $X$ is now some irreducible representation given by the cut-down. As any Kauffman polynomial planar algebra is braided, $\Rep(\text{Aut}(\Gamma),X)$ must take a map from the planar algebra of unoriented virtual tangles by Proposition \ref{prop:symmetric}. Two examples of this are given below:

\begin{example}[The pentagon spin model for the Kauffman polynomial]\label{ex:pent}
The pentagon gives a spin model for the $B_2$ spider when the value of $Q$ in \cite{Kup97} is equal to $e^{i\pi/5}$. The pentagon has automorphism group $D_{10},$ the dihedral group of 10 elements, so the even portion of the sub-planar algebra of the $B^2$ spider given by the doubled strands at this value of $Q$ is a sub-planar algebra of $\Rep(D_{10},X),$ where $X$ is one of the two 2-dimensional representations of $D_{10}$. In this case, the choice of $X$ is irrelevant as they both appear as summands of $\text{Perm}^G_H$. Moreover, the dimensions of $\Rep(D_{10},X)$ and $B_{2}^\text{even}$ show that the two are actually isomorphic. By Corollary \ref{otmaps} and the fact that $\Rep(D_{10})$ is non-degenerate, we know that there exists the following map:
\begin{equation*}
    \Rep(O(2)) \hookrightarrow \Rep(D_{10},X)
\end{equation*}
We note here that if $H \hookrightarrow G$ then $\Rep(G) \hookrightarrow \Rep(H).$ Because, $B_2^\text{even}$ is braided, we know by Proposition \ref{prop:symmetric} that $\Rep(D_{10})$ is a quotient of the planar algebra of unoriented virtual tangles. Thus, we can give a skein-theoretic description of $\Rep(D_{10})$ by determining which quotient it corresponds to. As we will see in Section \ref{virtualskein}, $\Rep(D_{10})$ is isomorphic to $\Rep(O(2),e^{i\pi/5}),$ the quotient of the planar algebra of unoriented virtual tangles with $\dim V_0 = \dim V_2 = 1$ and $\dim V_4 = 3$, subject to the following relations:
\begin{equation*}
\begin{tikzpicture}[baseline]
\begin{knot}[clip width = 4]
	\strand (-.5,.5) -- (.5,-.5);
	\strand (-.5,-.5) -- (.5,.5);
\end{knot}
\end{tikzpicture}
 ~=~ \dfrac{Q^{-1} - Q}{2}\cdot\begin{tikzpicture}[baseline, rounded corners = 5mm]
    \draw (.5,.5) -- (.1,0) -- (.5,-.5);
    \draw (-.5,.5) -- (-.1,0) -- (-.5,-.5);
\end{tikzpicture} ~+~ \dfrac{Q^{-1} - Q}{2}\cdot\begin{tikzpicture}[baseline, rounded corners = 5mm]
    \draw (-.5,-.5) -- (0,.1) -- (.5,-.5);
    \draw (-.5,.5) -- (0,-.1) -- (.5,.5);
\end{tikzpicture} ~+~ \dfrac{1 + \sqrt{5}}{4}\cdot\begin{tikzpicture}[baseline]
    \draw (-.5,-.5)-- (.5,.5);
    \draw (.5,-.5) -- (-.5,.5);
    \draw[fill=white] (0,0) circle (.1cm);
\end{tikzpicture}

\end{equation*}
\begin{equation*}
\skeincircle ~=~ 2 
\end{equation*}
where $Q = e^{i\pi/5}.$
\end{example}

\begin{example}[The Higman-Sims spin model for the Kauffman polynomial]\label{ex:hs}
Now consider the specialization of the $B_2$ spider with spin model given by the Higman-Sims graph, which occurs at $Q = \tau^2,$ the square of the golden ratio. The automorphism group of this graph is called $HS.2.$ Because the Higman-Sims graph is transitive, we know that the image of the even portion of the sub-planar algebra of the $B_2$ spider generated by the doubled strand at $Q = \tau^2$ is contained in $\Rep(HS.2,X_k)$, where the summand of $\text{Perm}^G_H$ is $22$-dimensional. Like with the pentagon we obtain a map
\begin{equation*}
    \Rep(O(22)) \hookrightarrow \Rep(HS.2,X_k)
\end{equation*}
Unlike with the pentagon, however, the even portion of the sub-planar algebra of $B_2$ generated by the doubled strand is \emph{not} isomorphic to $\Rep(HS.2,X_k),$ which can again be seen by comparing the dimensions of the planar algebras. In particular, $\Rep(HS.2)$ contains an additional 5-box. As such, a skein theoretic description of $\Rep(HS.2)$ is not possible using the same methods as in the case of the pentagon.  For more spin models arising from the Kaufman polynomial planar algebra, please see \cite{Edg18}.

Since $B_2^\text{even}$ is braided, we know $\Rep(HS.2)$ takes a map from the planar algebra of unoriented virtual tangles by Proposition \ref{prop:symmetric}. By inspection, the image of this planar algebra has a 4-dimensional 4-box space. Further, no skein-theoretic description of this sub-planar algebra has yet been given. The dimension bounds in Section \ref{virtualskein} were chosen to include all quotients of the planar algebra of unoriented virtual tangles with 4-box space strictly smaller than $\Rep(HS.2)$. It would be interesting to expand the classification to $\dim V_4 = 4$, so that we can include this sub-planar algebra of $\Rep(HS.2)$. 
\end{example}

\subsection{Examples of symmetric trivalent planar algebras.}\label{trivalent}

For all previous planar algebras that have been presented, it has been true that the odd-box spaces have been trivial. Now suppose that we consider the unshaded planar algebra generated by a trivalent vertex
\begin{center}
\begin{tikzpicture}[baseline = -.05cm, rounded corners = 5mm]
    \draw (0,-.6) -- (0,0)
    	  (0,0) -- (-.5,.5)
    	  (0,0) -- (.5,.5);
\end{tikzpicture}
\end{center}
which is rotationally invariant. All non-degenerate quotients of the planar algebra of \emph{planar} trivalent graphs with initial box space dimension bounds $1,0,1,1,4,11,41$ were classified in \cite{MPS15}, which the authors call a \emph{trivalent category}, as these planar algebras arise from tensor categories in the same manner as the planar algebra $\Rep(G)$ does.
\begin{theorem}[\cite{MPS15}]\label{thm:triclass}
The following table is an exhaustive list of non-degenerate trivalent planar algebras with initial box space dimension bounds $1,0,1,1,4,11,41$:
\begin{center}
\begin{tabular}{l | l}
Dimension bounds & Name \\
\hline 
$1,0,1,1,2,\dots$ & $\text{SO}(3)_{\zeta_5}$ \\
$1,0,1,1,3,\dots$ & $\text{SO}(3)_q$ or $OSp(1|2)$ \\
$1,0,1,1,4,8,\dots$ & $ABA$ \\
$1,0,1,1,4,9,\dots$ & $(G_2)_{\zeta_{20}}$ \\
$1,0,1,1,4,10,\dots$ & $(G_2)_q$ \\
$1,0,1,1,4,11,37,\dots$ & $H3$
\end{tabular}
\end{center}
where $ABA$ is a sub-planar algebra of the free product of $TL(\sqrt{dt^{-1}})\ast TL(t)$ and $H3$ is the fusion category found in \cite{GS12}.
\end{theorem}

If we also include \begin{tikzpicture}[baseline = -.15cm, scale = .65]
    \draw (-.5,-.5)-- (.5,.5);
    \draw (.5,-.5) -- (-.5,.5);
    \draw[fill=white] (0,0) circle (.1cm);
\end{tikzpicture}
 as a generator and impose the virtual Reidemeister moves, we can really think of this planar algebra as being the set of graphs with only trivalent and univalent vertices. In this case, the virtual crossing would simply represent two edges intersecting. This leads us to our next definition:

\begin{definition}
The planar algebra of trivalent graphs is a collection of vector spaces $\mathcal{V}=\{V_0,V_1,\dots\}$ where $V_n$ consists of simply-laced graphs (i.e. no multiple edges) with $n$ numbered univalent vertices and some number of trivalent vertices. Two graphs are considered equivalent if there exists a graph isomorphism between them that respects the numbering of the univalent vertices.
\end{definition}

Note that the graph isomorphism property implies that the planar algebra is symmetric. We will study quotients of this planar algebra with dimension bounds $1,0,1,1,4,11,41.$ 

\begin{definition}\label{def:vtri}
A symmetric trivalent planar algebra is a quotient of the planar algebra of trivalent graphs with $\dim V_1 =0$ and $\dim V_0 = \dim V_2 = \dim V_3 = 1$ that satisfies the virtual Reidemeister moves and the following additional relations:
\begin{center}
\begin{tabular}{cccc}
$\skeincircle = d = t - 1$ && $\begin{tikzpicture}[baseline]
 	\draw (0,-.6) -- (0,-.2);
    \draw (0,-.2) to[out=45,in=-90] (.35,.2) to[out=90,in=0] (0,.6)
    			  to[out=180,in=90] (-.35,.2) to[out=-90,in=135] (0,-.2);
\end{tikzpicture}
 ~=~ 0$\\
\\
$\begin{tikzpicture}[baseline]
 	\draw (0,.75) -- (0,.45);
	\draw (0,-.75) -- (0,-.45);
    \draw (0,-.45) .. controls (-.25,0) .. (0,.45);
    \draw (0,-.45) .. controls (.25,0) .. (0,.45);
\end{tikzpicture}
~=~
c_1 \cdot ~
\begin{tikzpicture}[baseline]
 	\draw (0,.75) -- (0,-.75);
\end{tikzpicture}
$ && \includegraphics[valign = c,scale = 1.25]{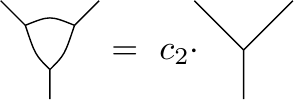} \\
\\
$\begin{tikzpicture}[scale = .7, baseline = .1cm]
    \draw (0,-.6) -- (0,0)
    	  (0,0) -- (-.4,.4)
    	  (0,0) -- (.4,.4);
    \draw (-.4,.4) to[out = 135, in = -135] (-.4, .6)
    			   to[out = 45, in = -135]  (.5,1.1);
    \draw (.4,.4)  to[out = 45, in = -45]   (.4, .6)
   				   to[out = 135, in = -45]  (-.5,1.1);
   	\draw[fill=white] (0,.825) circle (.075cm);			   
\end{tikzpicture}
=

\begin{tikzpicture}[baseline = -.1cm]
    \draw (0,-.6) -- (0,0)
    	  (0,0) -- (-.5,.5)
    	  (0,0) -- (.5,.5);
\end{tikzpicture}
$

&&

$
\begin{tikzpicture}[baseline = -.1cm]
    \draw (0,-.6) -- (0,0)
    	  (0,0) -- (-.5,.5)
    	  (0,0) -- (.5,.5);
   \draw  (0,.5) to[out=-90,in=45] (-.2,.2)
   				  to[out=225, in=90] (-.35,-.6);
   \draw[fill=white] (-.2,.2) circle (.075cm);
\end{tikzpicture}
~
=
~
\begin{tikzpicture}[baseline = -.1cm]
    \draw (0,-.6) -- (0,0)
    	  (0,0) -- (-.5,.5)
    	  (0,0) -- (.5,.5);
     \draw (0,.5)  to[out=-90,in=75] (.25,.25)
     				to[out=-90,in=45] (0,-.35)
     				to[out=225,in=90] (-.35,-.6);
     \draw[fill=white] (.25,.25) circle (.075cm);
	 \draw[fill=white] (0,-.35) circle (.075cm);			
\end{tikzpicture}
$
\end{tabular}
\end{center}
\end{definition}

\begin{remark}
In \cite{MPS15}, a trivalent planar algebra is generated by only the trivalent vertex. A symmetric trivalent planar algebra, though, is generated by \emph{both} the trivalent vertex and the symmetric crossing. Thus, while every symmetric planar algebra appearing in \cite{MPS15} is a symmetric trivalent planar algebra, not every symmetric trivalent planar algebra appears in \cite{MPS15}. 
\end{remark}

We will first classify all symmetric trivalent planar algebras and then classify all their sub-braidings. Some examples of planar algebras that do satisfy this property are given presently:

\begin{example}[$\Rep(SO(3))$]\label{ex:so3}
Here the chosen faithful, irreducible representation is the standard representation of $SO(3)$. Thus, $d=3$. It has initial box space dimensions $1,0,1,1,3$ and--in addition to the relations mentioned above--has the following relations:
\begin{center}
\begin{tabular}{c}
$\begin{tikzpicture}[baseline]
    \draw (-.5,-.5)-- (.5,.5);
    \draw (.5,-.5) -- (-.5,.5);
    \draw[fill=white] (0,0) circle (.1cm);
\end{tikzpicture}
 ~=~ \begin{tikzpicture}[baseline, rounded corners = 5mm]
    \draw (-.5,-.5) -- (0,.1) -- (.5,-.5);
    \draw (-.5,.5) -- (0,-.1) -- (.5,.5);
\end{tikzpicture} ~-~ 2\cdot \stH$ 
\\
\\
$\stI ~-~ \stH ~=~ \begin{tikzpicture}[baseline, rounded corners = 5mm]
    \draw (.5,.5) -- (.1,0) -- (.5,-.5);
    \draw (-.5,.5) -- (-.1,0) -- (-.5,-.5);
\end{tikzpicture} ~-~ \begin{tikzpicture}[baseline, rounded corners = 5mm]
    \draw (-.5,-.5) -- (0,.1) -- (.5,-.5);
    \draw (-.5,.5) -- (0,-.1) -- (.5,.5);
\end{tikzpicture}$\\
\end{tabular}
\end{center}
The latter relation is usually referred to as an ``I=H'' relation. The fact that the first relation is actually rotationally invariant is implied by the I=H relation. 
\end{example}

\begin{example}[$\Rep(OSp(1|2))$]\label{ex:osp}
In this case the irreducible representation is the standard representation of $OSp(1|2)$. Thus, the circle parameter is $1-2 = -1$ and--in addition to the relations mentioned above-has the following relations:
\begin{center}
\begin{tabular}{c}
$\begin{tikzpicture}[baseline]
    \draw (-.5,-.5)-- (.5,.5);
    \draw (.5,-.5) -- (-.5,.5);
    \draw[fill=white] (0,0) circle (.1cm);
\end{tikzpicture}
 ~=~ -2\cdot\begin{tikzpicture}[baseline, rounded corners = 5mm]
    \draw (.5,.5) -- (.1,0) -- (.5,-.5);
    \draw (-.5,.5) -- (-.1,0) -- (-.5,-.5);
\end{tikzpicture} ~-~ \begin{tikzpicture}[baseline, rounded corners = 5mm]
    \draw (-.5,-.5) -- (0,.1) -- (.5,-.5);
    \draw (-.5,.5) -- (0,-.1) -- (.5,.5);
\end{tikzpicture} ~+~ 2\cdot \stH$
\\
\\
$\stI ~-~ \stH ~=~ -\dfrac{1}{2}\left[~\begin{tikzpicture}[baseline, rounded corners = 5mm]
    \draw (.5,.5) -- (.1,0) -- (.5,-.5);
    \draw (-.5,.5) -- (-.1,0) -- (-.5,-.5);
\end{tikzpicture} ~-~ \begin{tikzpicture}[baseline, rounded corners = 5mm]
    \draw (-.5,-.5) -- (0,.1) -- (.5,-.5);
    \draw (-.5,.5) -- (0,-.1) -- (.5,.5);
\end{tikzpicture}~\right]$
\end{tabular}
\end{center}
Again, the fact that the formula for the virtual crossing is actually rotationally invariant is implied by the given I=H relation.
\end{example}

\begin{example}[$\Rep(G_2)$]
In this case the representation is the standard, 7-dimensional representation of the compact real Lie group $G_2.$ Thus, the circle parameter is $7$. It does not have an I=H relation but does have the following relation for the virtual crossing:
\begin{equation*}
\begin{tikzpicture}[baseline]
    \draw (-.5,-.5)-- (.5,.5);
    \draw (.5,-.5) -- (-.5,.5);
    \draw[fill=white] (0,0) circle (.1cm);
\end{tikzpicture}
 ~=~ \dfrac{1}{2}\left[~\begin{tikzpicture}[baseline, rounded corners = 5mm]
    \draw (.5,.5) -- (.1,0) -- (.5,-.5);
    \draw (-.5,.5) -- (-.1,0) -- (-.5,-.5);
\end{tikzpicture} ~+~ \begin{tikzpicture}[baseline, rounded corners = 5mm]
    \draw (-.5,-.5) -- (0,.1) -- (.5,-.5);
    \draw (-.5,.5) -- (0,-.1) -- (.5,.5);
\end{tikzpicture} \right] ~-~ 3\left[~ \stI + \stH ~\right]
\end{equation*}
\end{example}

\begin{example}[$\Rep(S_{3})$]\label{ex:s3}
This planar algebra is generated by the faithful two-dimensional representation of $S_3,$ the permutation group on 3 letters. Thus, it has circle parameter $2$ and the following relations:
\begin{equation*}
\begin{tikzpicture}[baseline]
    \draw (-.5,-.5)-- (.5,.5);
    \draw (.5,-.5) -- (-.5,.5);
    \draw[fill=white] (0,0) circle (.1cm);
\end{tikzpicture}
 ~=~ \stI ~+~ \stH
\end{equation*}
\begin{equation*}
\stI ~-~ \stH ~=~ \begin{tikzpicture}[baseline, rounded corners = 5mm]
    \draw (.5,.5) -- (.1,0) -- (.5,-.5);
    \draw (-.5,.5) -- (-.1,0) -- (-.5,-.5);
\end{tikzpicture} ~-~ \begin{tikzpicture}[baseline, rounded corners = 5mm]
    \draw (-.5,-.5) -- (0,.1) -- (.5,-.5);
    \draw (-.5,.5) -- (0,-.1) -- (.5,.5);
\end{tikzpicture}
\end{equation*}
Note that while the planar algebra $\Rep(SO(3))$ and $\Rep(S_3)$ have the same I=H relation, they have different formulas for the virtual crossing, and so they are non-isomorphic planar algebras.
\end{example}

The previous example gives rise to a general class of planar algebras called Deligne's $S_t$ found in \cite{EGNO15} (Section 9.12.1) and originally in \cite{Del07}. When $t$ is a non-negative integer, $S_t$ is degenerate and has $\Rep(S_t)$ as its only non-degenerate quotient, which is proven in the following theorem from \cite{CO11} and \cite{Del07}:

\begin{theorem}[\cite{CO11} and \cite{Del07}]\label{thm:dim}
At generic values of $t\in\mathbb{C}$, the planar algebra $S_t$ has initial box space dimensions $1,0,1,1,4,11,41$ and is non-degenerate. When $t$ is a non-negative integer, $S_t$ with the aforementioned dimensions exists but is a degenerate planar algebra. The only non-degenerate quotient of that planar algebra is isomorphic to the planar algebra $\Rep(S_t)$ with $n$-box space dimensions equal to $\dim\text{Hom}(\mathbbm{1},X^{\otimes n}),$ where $X$ is the standard $t-1$ dimensional representation of $S_t.$ 
\end{theorem}

We should note here that at $t=0,$ the non-degenerate quotient of $S_0$ is isomorphic to $\Rep(OSp(1|2)).$ The definition given below of Deligne's $S_t$ from \cite{EGNO15}{Section 9.12.1} using the standard representation and not the permutation representation (See \cite{DO14} for more information). This is why the circle value is $t-1$ and not $t$.

\begin{example}(Deligne's $S_t$)\label{ex:St}
 Deligne's $S_t$ is a symmetric trivalent planar algebras subject to the additional I=H relation
\begin{equation*}
\stI - \stH ~=~ \dfrac{1}{t-2}\left( ~ \begin{tikzpicture}[baseline, rounded corners = 5mm]
    \draw (.5,.5) -- (.1,0) -- (.5,-.5);
    \draw (-.5,.5) -- (-.1,0) -- (-.5,-.5);
\end{tikzpicture} - \begin{tikzpicture}[baseline, rounded corners = 5mm]
    \draw (-.5,-.5) -- (0,.1) -- (.5,-.5);
    \draw (-.5,.5) -- (0,-.1) -- (.5,.5);
\end{tikzpicture} ~ \right)
\end{equation*}
It turns out that this planar algebra has an interesting sub-braiding on it given by 
\begin{equation*}
\begin{tikzpicture}[baseline]
\begin{knot}[clip width = 4]
	\strand (-.5,.5) -- (.5,-.5);
	\strand (-.5,-.5) -- (.5,.5);
\end{knot}
\end{tikzpicture}
 ~=~ (q^2-1)\cdot\begin{tikzpicture}[baseline, rounded corners = 5mm]
    \draw (.5,.5) -- (.1,0) -- (.5,-.5);
    \draw (-.5,.5) -- (-.1,0) -- (-.5,-.5);
\end{tikzpicture} ~+~ q^{-2}\begin{tikzpicture}[baseline, rounded corners = 5mm]
    \draw (-.5,-.5) -- (0,.1) -- (.5,-.5);
    \draw (-.5,.5) -- (0,-.1) -- (.5,.5);
\end{tikzpicture} ~-~(q^2+q^{-2})\stI
\end{equation*}
with $t = q^2 + 2 + q^{-2}.$ We will prove that this is a sub-braiding and classify all such sub-braidings of symmetric trivalent planar algebras in Section \ref{tri}. Those familiar with the quantum group $SO(3)_q$--also called the second-colored TLJ planar algebra (See \cite{CDS95} and \cite{KL94})--will recognize this formula as the standard formula for the braiding on that planar algebra. The fact that it appears as a braiding for $S_t$ arises from the fact that there is exists a map
\[ 
SO(3)_q \hookrightarrow S_t 
\]
when $t = q^2 + 2 + q^{-2}.$ Given this relationship, $S_t$ is often called the virtual second-colored TLJ planar algebra.  
\end{example}

\section{Classification of simple virtual skein theories for unoriented virtual tangles} \label{virtualskein}
In this section, we would like to classify certain quotients of the planar algebra of unoriented virtual tangles, described in Section \ref{virtualtangles}. Let $\mathcal{V}:=\{V_i\}_{i=0,1\dots}$ be such a quotient. Then $\mathcal{V}$ inherits the relations from the planar algebra of unoriented virtual tangles. That is, $\mathcal{V}$ must satisfy the Reidemeister relations and their virtual counterparts. 

Moreover, we will require that $\dim V_0 = \dim V_2 = 1$ and that $\dim V_4 \leq 3$. This implies that $\mathcal{V}$ has the following relations:
\begin{center}
\begin{tabular}{ccccccc} 
$\skeincircle ~ = ~ d$ &&& $\begin{tikzpicture}[baseline=.5cm, rounded corners = 5mm]
\draw (0,1)--(.4,.6);
\draw (.6,.4)--(1,0)--(0,0)--(1,1);
\end{tikzpicture} ~ = ~ a \cdot \begin{tikzpicture}[baseline = .5 cm, rounded corners = 6mm]
    \draw (0,1)--(.5,-.1) -- (1,1);
\end{tikzpicture}$
&&& $\begin{tikzpicture}[baseline=.5cm, rounded corners = 5mm]
\draw (0,1)--(1,0) -- (0,0)--(1,1);
\draw[fill = white] (.5,.5) circle (.1cm);
\end{tikzpicture} ~=~ \begin{tikzpicture}[baseline = .5 cm, rounded corners = 6mm]
    \draw (0,1)--(.5,-.1) -- (1,1);
\end{tikzpicture}$ \\
circle parameter &&& twist relation &&& virtual Reidemeister I
\end{tabular}
\end{center}
where $d,a\in\mathbb{C}^\times$. Given that $\dim V_4 \leq 3$, then there must be at least one relation among  
\begin{center}
\begin{tabular}{cccccccccc}
\begin{tikzpicture}[baseline, rounded corners = 5mm]
    \draw (-.5,-.5) -- (0,.1) -- (.5,-.5);
    \draw (-.5,.5) -- (0,-.1) -- (.5,.5);
\end{tikzpicture} && \begin{tikzpicture}[baseline, rounded corners = 5mm]
    \draw (.5,.5) -- (.1,0) -- (.5,-.5);
    \draw (-.5,.5) -- (-.1,0) -- (-.5,-.5);
\end{tikzpicture} && \begin{tikzpicture}[baseline]
\begin{knot}[clip width = 4]
	\strand (-.5,.5) -- (.5,-.5);
	\strand (-.5,-.5) -- (.5,.5);
\end{knot}
\end{tikzpicture}
 && \begin{tikzpicture}[baseline]
\begin{knot}[clip width = 4]
	\strand (-.5,-.5) -- (.5,.5);
	\strand (-.5,.5) -- (.5,-.5);
\end{knot}
\end{tikzpicture} && \begin{tikzpicture}[baseline]
    \draw (-.5,-.5)-- (.5,.5);
    \draw (.5,-.5) -- (-.5,.5);
    \draw[fill=white] (0,0) circle (.1cm);
\end{tikzpicture}

\end{tabular}
\end{center}
Notice that Examples \ref{ex:ot} and \ref{ex:vjones} are examples of such a quotient while a notable non-example is Example \ref{ex:hs}. Our goal for this section will be to prove the following theorem:

\setcounter{subsection}{1}

\begin{theorem}\label{thm:virtual}
Let $\mathcal{V}$ be a non-fully-flat quotient of the planar algebra of unoriented virtual tangles such that $\dim V_0 = \dim V_2 = 1$ and $\dim V_4 \leq 3$. Then $\mathcal{V}$ is isomorphic to one of the following:
\begin{enumerate}[label = \roman{enumi}.]
\item The virtual TLJ planar algebra
\item The Kauffman polynomial planar algebra at $a = \pm 1$ and $d = -2$ equipped with the virtual crossing:
\begin{equation*}
 \begin{tikzpicture}[baseline]
    \draw (-.5,-.5)-- (.5,.5);
    \draw (.5,-.5) -- (-.5,.5);
    \draw[fill=white] (0,0) circle (.1cm);
\end{tikzpicture}
 ~=~ -\left(~ \begin{tikzpicture}[baseline, rounded corners = 5mm]
    \draw (.5,.5) -- (.1,0) -- (.5,-.5);
    \draw (-.5,.5) -- (-.1,0) -- (-.5,-.5);
\end{tikzpicture} ~+~ \begin{tikzpicture}[baseline, rounded corners = 5mm]
    \draw (-.5,-.5) -- (0,.1) -- (.5,-.5);
    \draw (-.5,.5) -- (0,-.1) -- (.5,.5);
\end{tikzpicture} ~\right)   
\end{equation*}
\item $\Rep(O(2),a)$, where $a\in\mathbb{C}^\times$ denotes that the planar algebra is equipped with a braiding given by the formula
\begin{equation}
\begin{tikzpicture}[baseline]
\begin{knot}[clip width = 4]
	\strand (-.5,.5) -- (.5,-.5);
	\strand (-.5,-.5) -- (.5,.5);
\end{knot}
\end{tikzpicture}
 = \dfrac{a^{-1} - a}{2} ~ \begin{tikzpicture}[baseline, rounded corners = 5mm]
    \draw (.5,.5) -- (.1,0) -- (.5,-.5);
    \draw (-.5,.5) -- (-.1,0) -- (-.5,-.5);
\end{tikzpicture} - \dfrac{a^{-1}-a}{2} ~ \begin{tikzpicture}[baseline, rounded corners = 5mm]
    \draw (-.5,-.5) -- (0,.1) -- (.5,-.5);
    \draw (-.5,.5) -- (0,-.1) -- (.5,.5);
\end{tikzpicture} + \dfrac{a^{-1}+a}{2} ~ \begin{tikzpicture}[baseline]
    \draw (-.5,-.5)-- (.5,.5);
    \draw (.5,-.5) -- (-.5,.5);
    \draw[fill=white] (0,0) circle (.1cm);
\end{tikzpicture}
 
\end{equation}
\end{enumerate}
\end{theorem}
For more on the term ``fully-flat," see Definition \ref{trivial}.

\begin{lemma}\label{casea}
Let $\mathcal{V}$ be a quotient of the planar algebra of unoriented virtual tangles. If the cupcap and identity are linearly dependent in $\mathcal{V}$, then the planar algebra is fully-flat.
\end{lemma}

\begin{proof}
Suppose that
\begin{equation}\label{eq:casea}
\begin{tikzpicture}[baseline, rounded corners = 5mm]
    \draw (.5,.5) -- (.1,0) -- (.5,-.5);
    \draw (-.5,.5) -- (-.1,0) -- (-.5,-.5);
\end{tikzpicture} ~=~ z\cdot\begin{tikzpicture}[baseline, rounded corners = 5mm]
    \draw (-.5,-.5) -- (0,.1) -- (.5,-.5);
    \draw (-.5,.5) -- (0,-.1) -- (.5,.5);
\end{tikzpicture} 
\end{equation} 
By multiplying both sides of the equation by the virtual crossing, we obtain a relation relating the virtual crossing to a multiple of the cupcap (after using vR1). Performing a similar operation with the actual crossing, we obtain a similar relation with the crossing. This implies that the actual and virtual crossings are multiples of each other, which implies that the planar algebra is fully-flat, as desired.  
\end{proof}

\begin{lemma}\label{lemma:virtual}
Let $\mathcal{V}$ be a quotient of the planar algebra of unoriented virtual tangles in which the cupcap and the identity are linearly independent. If there is a relation between the virtual crossing, cupcap, and identity, it must be of the form
\begin{equation}\label{twodim}
\begin{tikzpicture}[baseline]
    \draw (-.5,-.5)-- (.5,.5);
    \draw (.5,-.5) -- (-.5,.5);
    \draw[fill=white] (0,0) circle (.1cm);
\end{tikzpicture}
 = -\left( ~\begin{tikzpicture}[baseline, rounded corners = 5mm]
    \draw (.5,.5) -- (.1,0) -- (.5,-.5);
    \draw (-.5,.5) -- (-.1,0) -- (-.5,-.5);
\end{tikzpicture} + \begin{tikzpicture}[baseline, rounded corners = 5mm]
    \draw (-.5,-.5) -- (0,.1) -- (.5,-.5);
    \draw (-.5,.5) -- (0,-.1) -- (.5,.5);
\end{tikzpicture} ~ \right)
\end{equation}
Moreover, in this instance the circle parameter is $-2$.
\end{lemma}

\begin{proof}
Because the virtual crossing is rotationally invariant, any relation must be of the form
\begin{equation}\label{eq:caseb}
\begin{tikzpicture}[baseline]
    \draw (-.5,-.5)-- (.5,.5);
    \draw (.5,-.5) -- (-.5,.5);
    \draw[fill=white] (0,0) circle (.1cm);
\end{tikzpicture}
 = z\left( ~\begin{tikzpicture}[baseline, rounded corners = 5mm]
    \draw (.5,.5) -- (.1,0) -- (.5,-.5);
    \draw (-.5,.5) -- (-.1,0) -- (-.5,-.5);
\end{tikzpicture} + \begin{tikzpicture}[baseline, rounded corners = 5mm]
    \draw (-.5,-.5) -- (0,.1) -- (.5,-.5);
    \draw (-.5,.5) -- (0,-.1) -- (.5,.5);
\end{tikzpicture} ~ \right)
\end{equation}
where $z\in\mathbb{C}^\times$. Squaring both sides of (\ref{eq:caseb}) and applying vR2 to the left-hand side gives the equation
\begin{equation*}
\begin{tikzpicture}[baseline, rounded corners = 5mm]
    \draw (.5,.5) -- (.1,0) -- (.5,-.5);
    \draw (-.5,.5) -- (-.1,0) -- (-.5,-.5);
\end{tikzpicture} ~ = ~  z^2\cdot \begin{tikzpicture}[baseline, rounded corners = 5mm]
    \draw (.5,.5) -- (.1,0) -- (.5,-.5);
    \draw (-.5,.5) -- (-.1,0) -- (-.5,-.5);
\end{tikzpicture} ~ + ~ z^2\cdot (2+d) \cdot\begin{tikzpicture}[baseline, rounded corners = 5mm]
    \draw (-.5,-.5) -- (0,.1) -- (.5,-.5);
    \draw (-.5,.5) -- (0,-.1) -- (.5,.5);
\end{tikzpicture}
\end{equation*}
Since the identity and the cupcap were assumed to be linearly independent, this tells us that $z^2 = 1$ and that $z^2\cdot(2+d) = 0$. Since $z^2 = 1 \neq 0,$ it must be the case that $d = -2$. Using this value of $d$ and capping off (\ref{eq:caseb}) at the top, we see that $z = -1,$ as desired.
\end{proof}

The above two lemmas allow us to prove the following result:

\begin{proposition}\label{cor:caseb}
All non-fully-flat quotients of the planar algebra of unoriented virtual tangles with $\dim V_4$  less than or equal to $3$ such that there is a linear dependence among the virtual crossing, cupcap and the identity are isomorphic to either the virtual TLJ planar algebra at $A=\pm 1$ or the Kauffman polynomial planar algebra at $d = -2$ and $a = \pm 1$ equipped with a virtual crossing given by the formula
\begin{equation*}
\begin{tikzpicture}[baseline]
    \draw (-.5,-.5)-- (.5,.5);
    \draw (.5,-.5) -- (-.5,.5);
    \draw[fill=white] (0,0) circle (.1cm);
\end{tikzpicture}
 ~=~ -\left(~\begin{tikzpicture}[baseline, rounded corners = 5mm]
    \draw (.5,.5) -- (.1,0) -- (.5,-.5);
    \draw (-.5,.5) -- (-.1,0) -- (-.5,-.5);
\end{tikzpicture} ~+~ \begin{tikzpicture}[baseline, rounded corners = 5mm]
    \draw (-.5,-.5) -- (0,.1) -- (.5,-.5);
    \draw (-.5,.5) -- (0,-.1) -- (.5,.5);
\end{tikzpicture}\right)
\end{equation*}
\end{proposition}
\begin{proof}
Let $\mathcal{V}$ be such a quotient. Then Lemma \ref{casea} tells us that the cupcap and identity must be linearly independent since $\mathcal{V}$ is not fullly-flat. By Lemma \ref{lemma:virtual}, if there is a relation among the identity, cupcap, virtual crossing, it must be (\ref{twodim}) and the circle parameter must be $-2$. 

Suppose this relation holds. If $\dim V_4 = 2$ then the cupcap and identity form a basis for $V_4$. This implies a relation of the form
\begin{equation*}
\begin{tikzpicture}[baseline]
\begin{knot}[clip width = 4]
	\strand (-.5,.5) -- (.5,-.5);
	\strand (-.5,-.5) -- (.5,.5);
\end{knot}
\end{tikzpicture}
 ~=~ \pm\left(~\begin{tikzpicture}[baseline, rounded corners = 5mm]
    \draw (.5,.5) -- (.1,0) -- (.5,-.5);
    \draw (-.5,.5) -- (-.1,0) -- (-.5,-.5);
\end{tikzpicture} ~+~ \begin{tikzpicture}[baseline, rounded corners = 5mm]
    \draw (-.5,-.5) -- (0,.1) -- (.5,-.5);
    \draw (-.5,.5) -- (0,-.1) -- (.5,.5);
\end{tikzpicture}\right)
\end{equation*}
which implies that $\mathcal{V}$ is fully-flat. Thus $\dim V_4 \ne 2.$

Suppose $\dim V_4 = 3.$ Then by Theorem \ref{thm:kauffman} we have a relation of the form
\begin{equation}\label{eq:kauffman2}
\begin{tikzpicture}[baseline]
\begin{knot}[clip width = 4]
	\strand (-.5,.5) -- (.5,-.5);
	\strand (-.5,-.5) -- (.5,.5);
\end{knot}
\end{tikzpicture}
 ~\pm~ \begin{tikzpicture}[baseline]
\begin{knot}[clip width = 4]
	\strand (-.5,-.5) -- (.5,.5);
	\strand (-.5,.5) -- (.5,-.5);
\end{knot}
\end{tikzpicture} = z\cdot\left[ \begin{tikzpicture}[baseline, rounded corners = 5mm]
    \draw (.5,.5) -- (.1,0) -- (.5,-.5);
    \draw (-.5,.5) -- (-.1,0) -- (-.5,-.5);
\end{tikzpicture} ~\pm~ \begin{tikzpicture}[baseline, rounded corners = 5mm]
    \draw (-.5,-.5) -- (0,.1) -- (.5,-.5);
    \draw (-.5,.5) -- (0,-.1) -- (.5,.5);
\end{tikzpicture} \right]
\end{equation}
Capping (\ref{eq:kauffman2}) on the top tells us that $a^{-1} - a = z$ in the $-$ case and that $a^{-1} - a = -z$ in the $+$ case. Since $\mathcal{V}$ is a quotient of the planar algebra of unoriented virtual tangles, we know that it must satisfy the mixed Reidemeister III relation
\begin{equation}
\begin{tikzpicture}[scale = .5, baseline,rounded corners = 5mm]
	\draw[white, name path=line 1] (-1,-1)-- (1,1);
    \draw[white, name path=line 2] (1,-1) -- (-1,1);
    \draw[white, name path=line 3] (0,-1) -- (1,0) -- (0,1);
    \draw (-1,-1)-- (-.1,-.1);
    \draw (.1,.1) -- (1,1);
    \draw (1,-1) -- (-1,1);
    \draw (0,-1) -- (1,0) -- (0,1);
	\draw[fill=white,name intersections={of=line 1 and line 3}]
    (intersection-1) circle (.15cm);
    \draw[fill=white,name intersections={of=line 2 and line 3}]
    (intersection-1) circle (.15cm);
\end{tikzpicture}
=
\begin{tikzpicture}[scale = .5, baseline,rounded corners = 5mm, xscale = -1]
	\draw[white, name path=line 1] (-1,-1)-- (1,1);
    \draw[white, name path=line 2] (1,-1) -- (-1,1);
    \draw[white, name path=line 3] (0,-1) -- (1,0) -- (0,1);
    \draw (1,-1)-- (.1,-.1);
    \draw (-.1,.1) -- (-1,1);
    \draw (-1,-1) -- (1,1);
    \draw (0,-1) -- (1,0) -- (0,1);
	\draw[fill=white,name intersections={of=line 1 and line 3}]
    (intersection-1) circle (.15cm);
    \draw[fill=white,name intersections={of=line 2 and line 3}]
    (intersection-1) circle (.15cm);
\end{tikzpicture}
\end{equation}
By expanding the virtual crossing on both sides of this formula using (\ref{twodim}), we obtain the following relation:
\begin{equation}\label{mR31}
\includegraphics[valign = c]{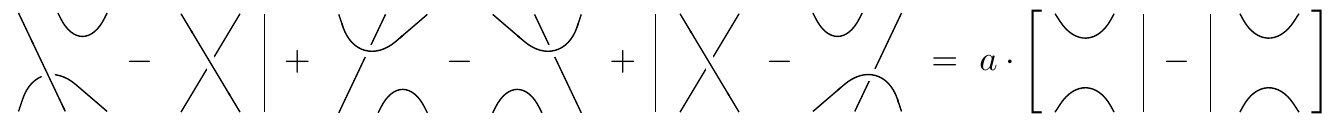}
\end{equation}
By solving the mixed Reidemeister move with the other crossing, one obtains a similar relation: 
\begin{equation}\label{mR32}
\includegraphics[valign = c]{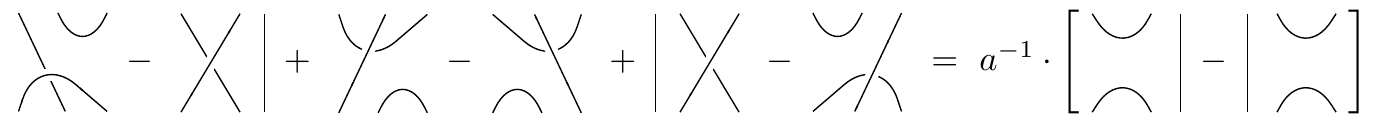}
\end{equation}
If we rotate and rearrange (\ref{mR31}), however, we obtain the negative of the left-hand side of (\ref{mR32}). This implies that $a=a^{-1}$ or that $a = \pm 1$. In the $-$ case, $z=0$ and the two crossings are equal. Thus, the actual crossing is symmetric or anti-symmetric. A quick calculation shows that in this skein theory, the forbidden move is actually satisfied. Because the forbidden move is satisfied, this planar algebra is isomorphic to the virtual TLJ planar algebra at $A = \pm 1$ with the isomorphism given by $\begin{tikzpicture}[baseline = -.15cm, scale = .65]
\begin{knot}[clip width = 4]
	\strand (-.5,.5) -- (.5,-.5);
	\strand (-.5,-.5) -- (.5,.5);
\end{knot}
\end{tikzpicture}
 \mapsto \begin{tikzpicture}[baseline = -.15cm, scale = .65]
    \draw (-.5,-.5)-- (.5,.5);
    \draw (.5,-.5) -- (-.5,.5);
    \draw[fill=white] (0,0) circle (.1cm);
\end{tikzpicture}
 $ and $\begin{tikzpicture}[baseline = -.15cm, scale = .65]
    \draw (-.5,-.5)-- (.5,.5);
    \draw (.5,-.5) -- (-.5,.5);
    \draw[fill=white] (0,0) circle (.1cm);
\end{tikzpicture}
 \mapsto \begin{tikzpicture}[baseline = -.15cm, scale = .65]
\begin{knot}[clip width = 4]
	\strand (-.5,.5) -- (.5,-.5);
	\strand (-.5,-.5) -- (.5,.5);
\end{knot}
\end{tikzpicture}
.$

In the $+$ case, we have that $d=-2.$ We obtain the relation $z=\mp 2.$ This gives the following relation:
\begin{equation*}
\begin{tikzpicture}[baseline]
\begin{knot}[clip width = 4]
	\strand (-.5,.5) -- (.5,-.5);
	\strand (-.5,-.5) -- (.5,.5);
\end{knot}
\end{tikzpicture}
 ~+~ \begin{tikzpicture}[baseline]
\begin{knot}[clip width = 4]
	\strand (-.5,-.5) -- (.5,.5);
	\strand (-.5,.5) -- (.5,-.5);
\end{knot}
\end{tikzpicture} ~=~ \pm 2 \left( ~ \begin{tikzpicture}[baseline, rounded corners = 5mm]
    \draw (.5,.5) -- (.1,0) -- (.5,-.5);
    \draw (-.5,.5) -- (-.1,0) -- (-.5,-.5);
\end{tikzpicture} ~+~ \begin{tikzpicture}[baseline, rounded corners = 5mm]
    \draw (-.5,-.5) -- (0,.1) -- (.5,-.5);
    \draw (-.5,.5) -- (0,-.1) -- (.5,.5);
\end{tikzpicture} ` \right)
\end{equation*}
This planar algebra is isomorphic to the Kauffman polynomial planar algebra at $a = \pm 1$ and $d = -2$. It is degenerate, however, as 

\begin{equation*}
\begin{tikzpicture}[baseline]
\begin{knot}[clip width = 4]
	\strand (-.5,.5) -- (.5,-.5);
	\strand (-.5,-.5) -- (.5,.5);
\end{knot}
\end{tikzpicture}
 \pm  \left( ~ \begin{tikzpicture}[baseline, rounded corners = 5mm]
    \draw (.5,.5) -- (.1,0) -- (.5,-.5);
    \draw (-.5,.5) -- (-.1,0) -- (-.5,-.5);
\end{tikzpicture} ~+~ \begin{tikzpicture}[baseline, rounded corners = 5mm]
    \draw (-.5,-.5) -- (0,.1) -- (.5,-.5);
    \draw (-.5,.5) -- (0,-.1) -- (.5,.5);
\end{tikzpicture} ` \right)
\end{equation*}
is in the null space of the inner-product matrix of $V_4$.

Thus, the only non-fully flat quotient of the planar algebra of unoriented virtual tangles with a dependence between the virtual crossing, cupcap, and the identity is the virtual TLJ planar algebra at $A = \pm 1$ or the Kauffman polynomial planar algebra with $d=-2$ and $a = \pm 1,$ as desired.
\end{proof}

If there is a relation between the virtual crossing, identity, and cupcap, then we know that the associated planar algebra is either a specialization of the virtual TLJ planar algebra at $A = \pm 1$ or the actual Kauffman polynomial planar algebra at $a=\pm 1$ and $d = -2$. As mentioned above, this planar algebra is degenerate in the sense of Definition \ref{degenerate}. If we were to quotient by the negligibles, we would obtain the virtual TLJ planar algebra at $A= \pm 1$.

\begin{lemma}\label{lemma:virtualjones}
Suppose that $\mathcal{V}$ is a quotient of the planar algebra of unoriented virtual tangles with initial box space dimensions $1,0,1,0,3$. Then if $\mathcal{V}$ has a relation of the form
\begin{equation*}
 \begin{tikzpicture}[baseline]
\begin{knot}[clip width = 4]
	\strand (-.5,.5) -- (.5,-.5);
	\strand (-.5,-.5) -- (.5,.5);
\end{knot}
\end{tikzpicture}
 ~ = ~ x\cdot \begin{tikzpicture}[baseline, rounded corners = 5mm]
    \draw (.5,.5) -- (.1,0) -- (.5,-.5);
    \draw (-.5,.5) -- (-.1,0) -- (-.5,-.5);
\end{tikzpicture} ~ + ~ y\cdot \begin{tikzpicture}[baseline, rounded corners = 5mm]
    \draw (-.5,-.5) -- (0,.1) -- (.5,-.5);
    \draw (-.5,.5) -- (0,-.1) -- (.5,.5);
\end{tikzpicture}
\end{equation*}
then $\mathcal{V}$ is isomorphic to the virtual TLJ planar algebra. 
\end{lemma}

\begin{proof}
By rotating the above equation and squaring we see that $y = x^{-1}$ and that the circle parameter $d = -(x^2+x^{-2}).$ By capping off the above equation we see that the twist parameter $a = -x^{-3}.$ Clearly, these are the exact skein relations of the virtual TLJ planar algebra described in Example \ref{ex:vjones}, and so the two planar algebras must be isomorphic. 
\end{proof}

\begin{proposition}\label{sixbox}
Let $\mathcal{V}$ be any quotient of the planar algebra of unoriented virtual tangles with initial box space dimensions $1,0,1,0,3$ with the following relation
\begin{equation}\label{eq:boxspace}
\begin{tikzpicture}[baseline]
\begin{knot}[clip width = 4]
	\strand (-.5,.5) -- (.5,-.5);
	\strand (-.5,-.5) -- (.5,.5);
\end{knot}
\end{tikzpicture}
 ~ = ~ x\cdot\begin{tikzpicture}[baseline, rounded corners = 5mm]
    \draw (.5,.5) -- (.1,0) -- (.5,-.5);
    \draw (-.5,.5) -- (-.1,0) -- (-.5,-.5);
\end{tikzpicture} ~ - ~ x\cdot\begin{tikzpicture}[baseline, rounded corners = 5mm]
    \draw (-.5,-.5) -- (0,.1) -- (.5,-.5);
    \draw (-.5,.5) -- (0,-.1) -- (.5,.5);
\end{tikzpicture} ~ + ~ z\cdot\begin{tikzpicture}[baseline]
    \draw (-.5,-.5)-- (.5,.5);
    \draw (.5,-.5) -- (-.5,.5);
    \draw[fill=white] (0,0) circle (.1cm);
\end{tikzpicture}

\end{equation}
where $x$ and $z$ are not simultaneously 0. Then $\dim V_6 \in \{10,15\}.$  Furthermore, when $\dim V_6 = 15,$ the planar algebra is fully-flat. When $\dim V_6 = 10$, we have that $\wedge^3 = 0.$
\end{proposition}

\begin{proof}
Let $\mathcal{V}$ be as described above. By multiplying (\ref{eq:boxspace}) with its 1-click rotation we obtain, we obtain the equation
\begin{equation*}
\begin{tikzpicture}[baseline, rounded corners = 5mm]
    \draw (.5,.5) -- (.1,0) -- (.5,-.5);
    \draw (-.5,.5) -- (-.1,0) -- (-.5,-.5);
\end{tikzpicture} ~ = ~ (z^2 - x^2)\cdot\begin{tikzpicture}[baseline, rounded corners = 5mm]
    \draw (.5,.5) -- (.1,0) -- (.5,-.5);
    \draw (-.5,.5) -- (-.1,0) -- (-.5,-.5);
\end{tikzpicture}
~ + ~ (2x^2 - d x^2 )\cdot\begin{tikzpicture}[baseline, rounded corners = 5mm]
    \draw (-.5,-.5) -- (0,.1) -- (.5,-.5);
    \draw (-.5,.5) -- (0,-.1) -- (.5,.5);
\end{tikzpicture}
\end{equation*}
Because the cupcap and the identity are linearly independent by assumption, we see that $z^2 - x^2 = 1$ and $x^2(2-d) = 0$. Thus, either $x = 0$ or $d = 2$. If $x = 0$, then $z = \pm 1$ which implies that $\mathcal{V}$ is fully-flat. 

Suppose now that $d = 2$ and let us consider $V_6$. Since $\dim V_4 = 3$, we know that any basis of $V_6$ can be written using only the virtual crossing, the cupcap, and the identity. In addition, vR3 tells us that both $6$-box diagrams with three virtual crossings are equal. Hence, we obtain the spanning set listed in Proposition \ref{prop:s2n}. Suppose that this spanning set were a basis. This would imply that $\dim V_6 = 15$. By expanding R3 using (\ref{eq:boxspace}) we obtain
\begin{center}
\includegraphics[valign = c]{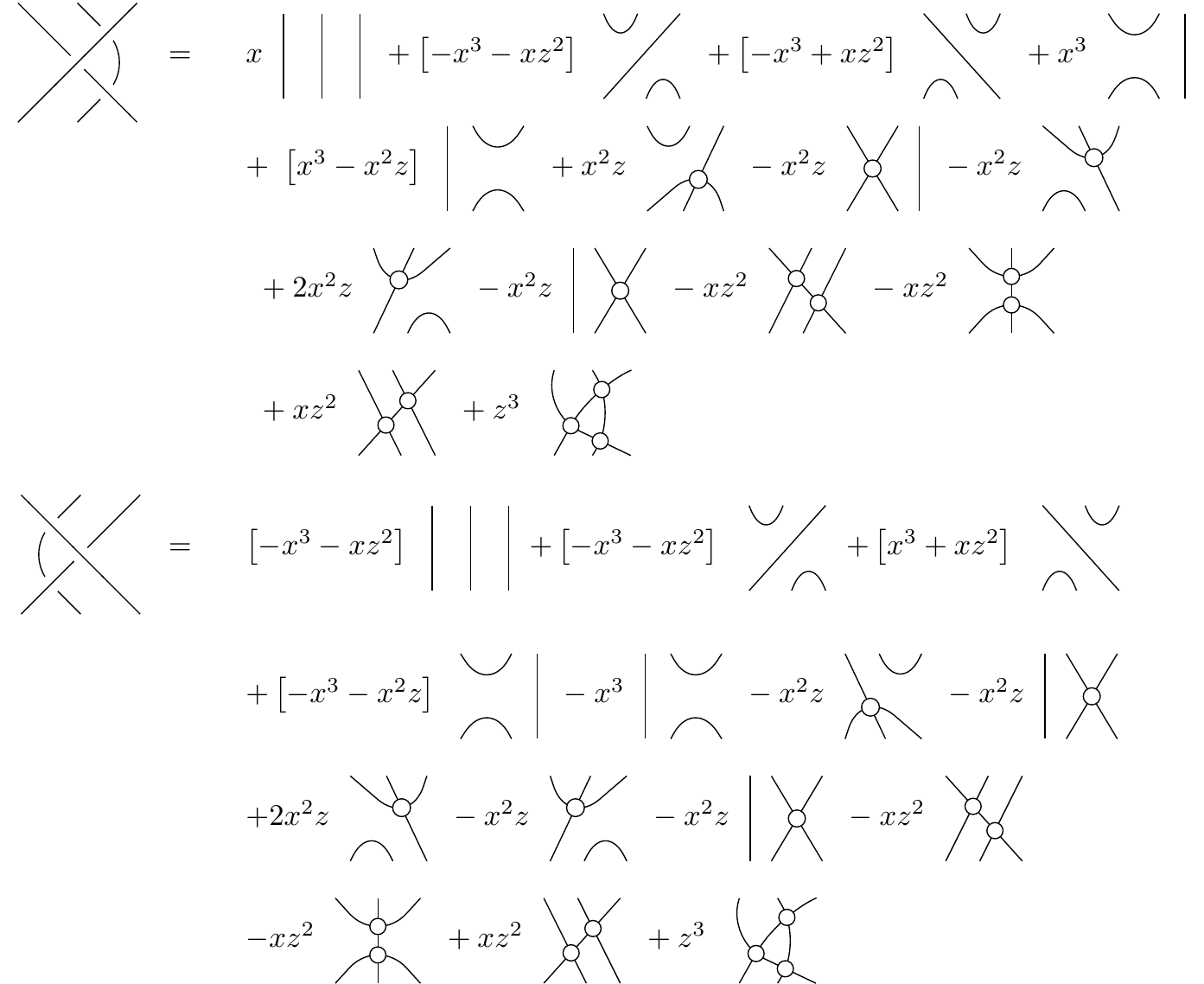}
\end{center}
Since these diagrams are all linearly independent, setting the right-hand sides of both equations above equal to each other gives a system of 15 equations. By inspection, the only solutions for $(x,z)$ are $\{(0,1), (0,-1), (0,0)\}.$ Since $x$ and $z$ are not simultaneously $0$, then $x=0$ and $z=\pm 1$ both of which imply $\mathcal{V}$ is fully-flat.

Suppose that the spanning set is not a basis. By Corollary \ref{cor:boxdim}, we know that the dimension of $V_6$ is $10$ as $\dim V_4 = 3$. In this case, the invariant subspace of the $5$-dimensional irreducible representation of $S_6$ gives us new relations in $V_6$. By inspection, that subspace gives us the following rotational eigenvector (see Definition \ref{def:eigen}):

\begin{equation}\label{eq:2rotation}
\includegraphics[valign = c]{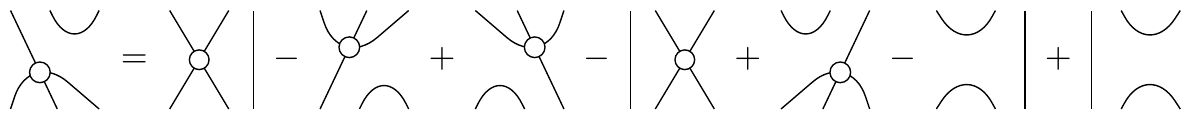}
\end{equation}
By multiplying both sides of (\ref{eq:2rotation}) by certain 6-boxes, we can obtain relations simplifying the following diagrams to terms with fewer crossings:
\begin{center}
    \includegraphics[valign = c]{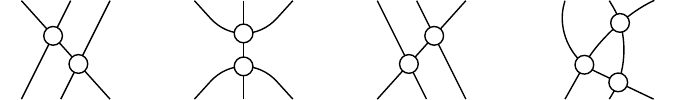}
\end{center}
In particular, we obtain the relation
\begin{equation*}
\includegraphics[valign = c]{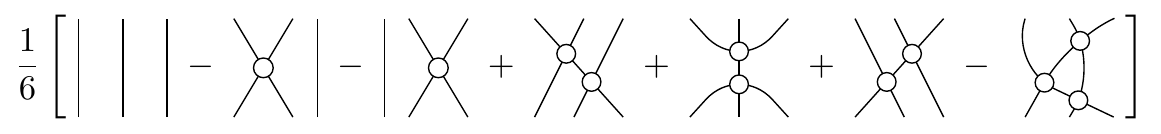} ~=~ 0
\end{equation*}
Since the left-hand side is $\wedge^3$ by Definition \ref{def:wedge}, we see that $\wedge^3 = 0$ in this case. This completes our proof. 
\end{proof}

We note here that the results of Proposition \ref{sixbox} followed from the work done in \cite{CH15}, but we give a quicker, diagrammatic proof of the result in our case. Nevertheless, this proposition tells us that if we have a planar algebra with the above crossing formula, it is either fully-flat or has the special property that $\wedge^3 = 0.$ This last fact will be key in proving the following proposition:

\begin{proposition}\label{prop:o2}
Let $\mathcal{V}$ be a quotient of the planar algebra of unoriented virtual tangles with $\dim V_4 = 3$ subject to the following relations:
\begin{equation*}
\skeincircle = 2, \hspace{2cm} \begin{tikzpicture}[baseline=.5cm, rounded corners = 5mm]
\draw (0,1)--(.4,.6);
\draw (.6,.4)--(1,0)--(0,0)--(1,1);
\end{tikzpicture} = a\begin{tikzpicture}[baseline = .5 cm, rounded corners = 6mm]
    \draw (0,1)--(.5,-.1) -- (1,1);
\end{tikzpicture}
\end{equation*}
\begin{equation}
\begin{tikzpicture}[baseline]
\begin{knot}[clip width = 4]
	\strand (-.5,.5) -- (.5,-.5);
	\strand (-.5,-.5) -- (.5,.5);
\end{knot}
\end{tikzpicture}
 = \dfrac{a^{-1} - a}{2}~\begin{tikzpicture}[baseline, rounded corners = 5mm]
    \draw (.5,.5) -- (.1,0) -- (.5,-.5);
    \draw (-.5,.5) -- (-.1,0) -- (-.5,-.5);
\end{tikzpicture} - \dfrac{a^{-1}-a}{2}~\begin{tikzpicture}[baseline, rounded corners = 5mm]
    \draw (-.5,-.5) -- (0,.1) -- (.5,-.5);
    \draw (-.5,.5) -- (0,-.1) -- (.5,.5);
\end{tikzpicture} + \dfrac{a^{-1}+a}{2}~\begin{tikzpicture}[baseline]
    \draw (-.5,-.5)-- (.5,.5);
    \draw (.5,-.5) -- (-.5,.5);
    \draw[fill=white] (0,0) circle (.1cm);
\end{tikzpicture}
 
\end{equation}
Then $\mathcal{V}$ is non-degenerate and is isomorphic to $\Rep(O(2),a)$, where $a$ denotes the above braiding formula. 
\end{proposition}

The proof of this proposition will appear in the next section. In this particular instance this family of possible skein theories has another strange property. The forbidden move (see Section \ref{virtualtangles} for the definition) is satisfied in this skein theory if and only if R2 and R3 are, which can be seen by using the so-called ``Kauffman trick:" 
\begin{equation*}
\includegraphics[valign = c]{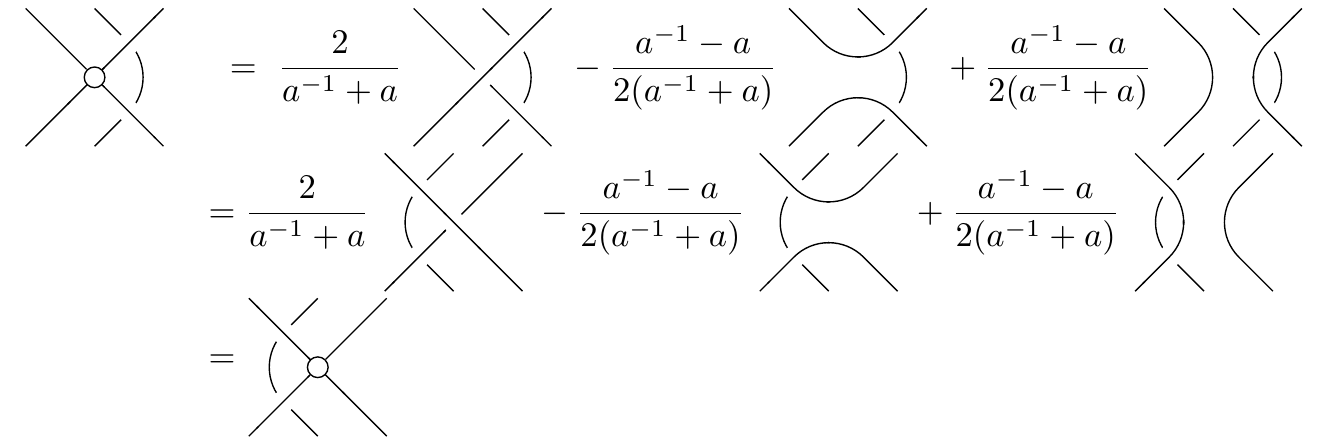}
\end{equation*}
\subsection{Proof of Proposition \ref{prop:o2} and Theorem \ref{thm:virtual}}

In order to show that this planar algebra is isomorphic to $\Rep(O(2),a)$, we will make an identification of the unoriented diagrams in this family of skein theories with oriented ones (Thanks to Pavel Etingof for his suggestion that led us to this trick.). 

\begin{proposition}\label{subalg}
$\Rep(S^1,a)$ in Example \ref{def:s1} contains $\Rep(O(2))$ as a sub-planar algebra. 
\end{proposition}

\begin{proof}
Let $\mathcal{V}$ be $\Rep(O(2))$ and $f: \mathcal{V}\to\Rep(S^1,a)$ be the map generated by
\begin{equation*}
\begin{tikzpicture}[baseline={([yshift=-.5ex]current bounding box.center)}]
\draw (0,-.5) -- (0,.5);
\end{tikzpicture}
~~~\mapsto~~~ \begin{tikzpicture}[baseline={([yshift=-.5ex]current bounding box.center)}]
\draw[->] (0,-.5) -- (0,.5);
\end{tikzpicture}
~~~+~~~
\begin{tikzpicture}[baseline]
\draw[<-] (0,-.5) -- (0,.5);
\end{tikzpicture}
\end{equation*}
Better put, $f$ is the map that sends any unoriented diagram to the sum of all possible orientations of the strands. For example, the overcrossing would map to
\begin{equation*}
\begin{tikzpicture}[baseline]
\begin{knot}[clip width = 4]
	\strand (-.5,.5) -- (.5,-.5);
	\strand (-.5,-.5) -- (.5,.5);
\end{knot}
\end{tikzpicture}
 ~\mapsto~ \overcrossinga ~+~ \overcrossingb ~+~ \overcrossingc ~+~ \overcrossingd 
\end{equation*}
We now show that this map is a map of planar algebras. In order to do this, we check that all relations in $\mathcal{V}$ are preserved under $f$. Since by Theorem \ref{thm:wedge} the circle parameter and the $\wedge^3 =0$ relation generate all the others, we need only check those relations still hold in the image of $f$. The circle relation is clearly preserved under $f$. By using a computer to aid in the simplification process (See $\mathtt{wedge3.nb}$ for the details of this computation), it is easy to see that the image of
\begin{equation}\label{eq:wedge}
\includegraphics[valign = c]{diagrams/pdf/Wedge3.pdf}
\end{equation}
is 0 under $f$. By looking at the image of the identity, cupcap, and virtual crossing, it is clear that $f$ is an injective map. Thus, $f(\mathcal{V})$ is isomorphic to $\Rep(O(2))$ by Theorem \ref{thm:wedge}, as desired.
\end{proof}

Thus, we see that $\Rep(S^1,a)$ contains $\Rep(O(2))$ as a sub-planar algebra. This relationship is rather surprising. Recall that $O(2) \cong S^1 \rtimes_{\alpha} \mathbb{Z}/2\mathbb{Z},$ where $\alpha$ is the automorphism that sends $e^{i\theta}$ to $e^{-i\theta}.$ Thus, any quotient of $S^1$ that is invariant under switching the direction of the arrows must contain $O(2)$, which explains this exotic braiding.

\begin{proof}[Proof of Proposition \ref{prop:o2}]
Let $\mathcal{V}$ be the planar algebra in Proposition \ref{prop:o2} and $a\in\mathbb{C}^\times.$ Let $f$ be as in Proposition \ref{subalg}, in which it was shown that $f$ is an injective map whose image is isomorphic to $\Rep(O(2))$. In order to show that this image is isomorphic to $\Rep(O(2),a)$, we also need to show that the formula for the crossing is preserved under $f$. Using the relations on these diagrams from the definition of $\Rep(S^1,a)$, we note that
\begin{equation}\label{eq:crossing2}
\begin{split}
\begin{tikzpicture}[baseline]
\begin{knot}[clip width = 4]
	\strand (-.5,.5) -- (.5,-.5);
	\strand (-.5,-.5) -- (.5,.5);
\end{knot}
\end{tikzpicture}
 &~\mapsto~ \overcrossinga ~+~ \overcrossingb ~+~ \overcrossingc ~+~ \overcrossingd 
\\
&~=~ a^{-1}\identitya + a\identityb + a^{-1}\identityc + a\identityd
\\
&~=~ a^{-1}\left[~ \identitya + \identityd ~\right] + a\left[~ \identityb + \identityc ~\right]
\end{split}
\end{equation}
Simplifying the right hand side gives
\begin{equation}\label{eq:oriented}
\begin{split}
\dfrac{a^{-1}-a}{2} &~ \begin{tikzpicture}[baseline, rounded corners = 5mm]
    \draw (.5,.5) -- (.1,0) -- (.5,-.5);
    \draw (-.5,.5) -- (-.1,0) -- (-.5,-.5);
\end{tikzpicture} -\dfrac{a^{-1}-a}{2}~\begin{tikzpicture}[baseline, rounded corners = 5mm]
    \draw (-.5,-.5) -- (0,.1) -- (.5,-.5);
    \draw (-.5,.5) -- (0,-.1) -- (.5,.5);
\end{tikzpicture} + \dfrac{a^{-1}+a}{2}\begin{tikzpicture}[baseline]
    \draw (-.5,-.5)-- (.5,.5);
    \draw (.5,-.5) -- (-.5,.5);
    \draw[fill=white] (0,0) circle (.1cm);
\end{tikzpicture}

\\
& \mapsto \dfrac{a^{-1}}{2}\cdot\left[~2\identitya + 2\identityd ~\right] + \dfrac{a}{2}\cdot\left[~2\cupcapa + 2\cupcapd ~\right]
\\
&
= a^{-1}\cdot\left[~\identitya + \identityd ~\right] + a\cdot\left[~\identityb + \identityc ~\right]
\end{split}
\end{equation}
Hence, (\ref{eq:crossing2}) and (\ref{eq:oriented}) are equal and so the formula for the crossing is preserved under $f$. 

Thus, since $f$ is injective $\mathcal{V}$ is isomorphic to a subplanar algebra of $\Rep(S^1,a)$. By Proposition \ref{subalg} and the preservation of the crossing formula under $f$, we know that $\mathcal{V} \cong \Rep(O(2),a),$ as desired. 
\end{proof}

Now that we have shown that $\Rep(O(2),a)$ is a sub-quotient of the planar algebra of unoriented virtual tangles, we can easily prove the following theorem: 

\begin{theorem}\label{thm:dimthree}
Let $\mathcal{V}$ with $\dim V_4 = 3$ be a non-fully flat quotient of the planar algebra of unoriented virtual tangles such that the virtual crossing, identity, and cupcap are linearly independent. Then $\mathcal{V}$ is isomorphic to one of the following:
\begin{enumerate}[label = \roman{enumi}.]
\item The virtual TLJ planar algebra
\item $\Rep(O(2),a)$.
\end{enumerate}
\end{theorem}

\begin{proof}
Assume that the virtual crossing, the cupcap, and the identity are linearly independent and $\dim V_4 = 3$. Since these diagrams now form a basis, we know that there must be a relation of the form   
\begin{equation}\label{eq:3dskein}
 \begin{tikzpicture}[baseline]
\begin{knot}[clip width = 4]
	\strand (-.5,.5) -- (.5,-.5);
	\strand (-.5,-.5) -- (.5,.5);
\end{knot}
\end{tikzpicture}
 ~ = ~ x\begin{tikzpicture}[baseline, rounded corners = 5mm]
    \draw (.5,.5) -- (.1,0) -- (.5,-.5);
    \draw (-.5,.5) -- (-.1,0) -- (-.5,-.5);
\end{tikzpicture} ~ + ~ y\begin{tikzpicture}[baseline, rounded corners = 5mm]
    \draw (-.5,-.5) -- (0,.1) -- (.5,-.5);
    \draw (-.5,.5) -- (0,-.1) -- (.5,.5);
\end{tikzpicture} ~ + ~ z\begin{tikzpicture}[baseline]
    \draw (-.5,-.5)-- (.5,.5);
    \draw (.5,-.5) -- (-.5,.5);
    \draw[fill=white] (0,0) circle (.1cm);
\end{tikzpicture}

\end{equation}
If $x = y = 0$ then $\mathcal{V}$ is fully-flat. If $z=0$ then by Lemma \ref{lemma:virtualjones}, we know that our quotient is isomorphic to the virtual TLJ planar algebra. Thus, let us assume that $z \ne 0$ and that one of $x$ or $y$ is non-zero. By rotating (\ref{eq:3dskein}) and multiplying the two equations together we obtain the following equation:
\begin{equation*}\label{eq:squaring3d}
\begin{tikzpicture}[baseline, rounded corners = 5mm]
    \draw (.5,.5) -- (.1,0) -- (.5,-.5);
    \draw (-.5,.5) -- (-.1,0) -- (-.5,-.5);
\end{tikzpicture} ~ = ~ (xy+z^2)\cdot\begin{tikzpicture}[baseline, rounded corners = 5mm]
    \draw (.5,.5) -- (.1,0) -- (.5,-.5);
    \draw (-.5,.5) -- (-.1,0) -- (-.5,-.5);
\end{tikzpicture}
~ + ~ (x^2 + y^2 + dxy + yz + xz)\cdot\begin{tikzpicture}[baseline, rounded corners = 5mm]
    \draw (-.5,-.5) -- (0,.1) -- (.5,-.5);
    \draw (-.5,.5) -- (0,-.1) -- (.5,.5);
\end{tikzpicture} ~ + ~(xz+yz)\cdot\begin{tikzpicture}[baseline]
    \draw (-.5,-.5)-- (.5,.5);
    \draw (.5,-.5) -- (-.5,.5);
    \draw[fill=white] (0,0) circle (.1cm);
\end{tikzpicture}

\end{equation*}
Since we assumed that these were linearly independent, we have the following equations:
\begin{align*}
  xy + z^2 &=~ 1 \\
  xz + yz &=~ 0 \\
  x^2 + y^2 + dxy + yz + xz &=~ 0
\end{align*}
Since $z \ne 0$, we know that $x = -y$ from the second equation. From the first equation, we know that $z^2 - x^2 = 1$. Additionally, the last equation tells us that $x^2 + x^2 - dx^2 =
x^2\cdot(2 - d) = 0.$ Since $x=0$ implies $y=0$, we see that $d = 2$. 

Using this information, we now have the simpler equation
\begin{equation}\label{eq:3dskein3}
 \begin{tikzpicture}[baseline]
\begin{knot}[clip width = 4]
	\strand (-.5,.5) -- (.5,-.5);
	\strand (-.5,-.5) -- (.5,.5);
\end{knot}
\end{tikzpicture}
 ~ = ~ x\begin{tikzpicture}[baseline, rounded corners = 5mm]
    \draw (.5,.5) -- (.1,0) -- (.5,-.5);
    \draw (-.5,.5) -- (-.1,0) -- (-.5,-.5);
\end{tikzpicture} ~ - ~ x\begin{tikzpicture}[baseline, rounded corners = 5mm]
    \draw (-.5,-.5) -- (0,.1) -- (.5,-.5);
    \draw (-.5,.5) -- (0,-.1) -- (.5,.5);
\end{tikzpicture} ~ + ~ z\begin{tikzpicture}[baseline]
    \draw (-.5,-.5)-- (.5,.5);
    \draw (.5,-.5) -- (-.5,.5);
    \draw[fill=white] (0,0) circle (.1cm);
\end{tikzpicture}

\end{equation} 
with $z^2 = 1 + x^2$. By capping (\ref{eq:3dskein3}) on the top, we see that $-x+z = a^{-1}$, and by capping (\ref{eq:3dskein3}) on the left we obtain the relation $x+z = a.$ This, gives a one parameter family of possible skein theories in terms of $a$ where (\ref{eq:3dskein3}) can be re-parameterized to 
\begin{equation}\label{eq:RepO2}
\begin{tikzpicture}[baseline]
\begin{knot}[clip width = 4]
	\strand (-.5,.5) -- (.5,-.5);
	\strand (-.5,-.5) -- (.5,.5);
\end{knot}
\end{tikzpicture}
 = \dfrac{a^{-1} - a}{2} ~ \begin{tikzpicture}[baseline, rounded corners = 5mm]
    \draw (.5,.5) -- (.1,0) -- (.5,-.5);
    \draw (-.5,.5) -- (-.1,0) -- (-.5,-.5);
\end{tikzpicture} - \dfrac{a^{-1}-a}{2} ~ \begin{tikzpicture}[baseline, rounded corners = 5mm]
    \draw (-.5,-.5) -- (0,.1) -- (.5,-.5);
    \draw (-.5,.5) -- (0,-.1) -- (.5,.5);
\end{tikzpicture} + \dfrac{a^{-1}+a}{2} ~ \begin{tikzpicture}[baseline]
    \draw (-.5,-.5)-- (.5,.5);
    \draw (.5,-.5) -- (-.5,.5);
    \draw[fill=white] (0,0) circle (.1cm);
\end{tikzpicture}
 
\end{equation}
By Proposition \ref{sixbox}, we know that $\wedge^3 = 0$ in $\mathcal{V}$. Thus, by Proposition \ref{prop:o2} $\mathcal{V}$ is isomorphic to $\Rep(O(2),a).$ As we have exhausted all possible relations, we know that $\mathcal{V}$ must be isomorphic to the virtual TLJ planar algebra or $\Rep(O(2),a),$ as desired. 
\end{proof}
We can now prove the results of Theorem \ref{thm:virtual}:

\begin{proof}[Proof of Theorem \ref{thm:virtual}]
Let $\mathcal{V}$ be a non-fully-flat quotient of the planar algebra of unoriented virtual tangles with $\dim V_0 = \dim V_2 = 1$ and $\dim V_4 \leq 3$. Then Proposition \ref{casea} states that the cupcap and identity must be linearly independent. If the cupcap and identity are linearly independent but the virtual crossing, identity and cupcap are linearly dependent, Corollary \ref{cor:caseb} tells us that $\mathcal{V}$ is isomorphic to the virtual TLJ planar algebra at $A=\pm 1$ or the Kauffman polynomial planar algebra at $d = -2$ and $a = \pm 1$.

If all three are linearly independent, then there must be a relation of the form
\begin{equation}
 \begin{tikzpicture}[baseline]
\begin{knot}[clip width = 4]
	\strand (-.5,.5) -- (.5,-.5);
	\strand (-.5,-.5) -- (.5,.5);
\end{knot}
\end{tikzpicture}
 ~ = ~ x\begin{tikzpicture}[baseline, rounded corners = 5mm]
    \draw (.5,.5) -- (.1,0) -- (.5,-.5);
    \draw (-.5,.5) -- (-.1,0) -- (-.5,-.5);
\end{tikzpicture} ~ + ~ y\begin{tikzpicture}[baseline, rounded corners = 5mm]
    \draw (-.5,-.5) -- (0,.1) -- (.5,-.5);
    \draw (-.5,.5) -- (0,-.1) -- (.5,.5);
\end{tikzpicture} ~ + ~ z\begin{tikzpicture}[baseline]
    \draw (-.5,-.5)-- (.5,.5);
    \draw (.5,-.5) -- (-.5,.5);
    \draw[fill=white] (0,0) circle (.1cm);
\end{tikzpicture}

\end{equation}
as we chose $\dim V_4 \leq 3.$ Thus, Theorem \ref{thm:dimthree} tells us that the planar algebra is either the virtual TLJ planar algebra or $\Rep(O(2),a).$ Since we have exhausted all possibilities, we see that, up to isomorphism, $\mathcal{V}$ must be the Kauffman polynomial planar algebra at $d=-2$ and $a = \pm 1$, the virtual TLJ planar algebra, or or $\Rep(O(2),a),$ which completes our proof.
\end{proof}

\subsection{Link invariants}

Since all of the above planar algebras are evaluable, they also provide virtual link invariants. The Kauffman polynomial planar algebra and virtual TLJ planar algebra are well-understood and the invariants they give can be found in \cite{Kau90} and \cite{Kau99}.

The $\Rep(O_{2},a)$ case could potentially give an interesting invariant. In fact, though, the invariant is rather boring after we normalize away the writhe:
\begin{proposition}
Let $K$ be a virtual knot and $\text{wr}(K)$ be the writhe of some fixed orientation of $K$. Then the knot invariant given by $\Rep(O(2),a)$ assigns $a^{\text{wr}(K)} + a^{-\text{wr}(K)}$ to $K$.
\end{proposition}

\begin{proof}
Let $K$ be a virtual knot. Thus it is an element of $V_0$ of $\Rep(O(2),a).$ Since $\Rep(O(2),a)$ is isomorphic to a subplanar algebra of $\Rep(S^1,a)$, $K$ can be thought of as the sum over all possible orientations of it. For knots, note that there are only two possible orientations on the strands. Moreover, the writhes of these knots are negatives of each other.

Fix an orientation on $K$ and let $\text{wr}(K)$ be its writhe. Then by using the relations of $\Rep(S^1,a),$ we perform the following algorithm. First, we resolve all virtual crossings by using the relation
\begin{equation*}
    \virtuala ~=~ \identitya
\end{equation*}
Note that this does not involve any scalars. Next, we resolve all actual crossings by using the relations
\begin{center}
\begin{tabular}{cccc}
$\overcrossinga ~=~ a\cdot\identitya$ &&& 
$\undercrossinga ~=~ a^{-1}\cdot\identitya$
\end{tabular}
\end{center}
The resulting diagram is some number of unknots multiplied by $a^{\text{wr}(K)}.$ Since the circle parameters are 1 in $\Rep(S^1,a),$ this diagram is equivalent to $a^{\text{wr}(K)}$ times the empty diagram. If we were to take the other orientation, we would obtain $a^{-\text{wr}(K)}$ times the empty diagram. Thus the knot $K$ is equivalent to $a^{\text{wr}(K)} + a^{-\text{wr}(K)}$, as desired.
\end{proof}

After normalizing away the writhe, we see that the invariant is trivial for all virtual knots. 

\section{Classification of small quotients of the planar algebra of oriented virtual tangles}

We will now classify all spherical quotients of the planar algebra of oriented virtual tangles, $\mathcal{V}$, such that $\dim V_0 = 1,$ the dimension of all vector spaces with 2 endpoints is 1, and the dimension of all vector spaces with 4 endpoints is 1 or 2. Like for the unoriented case, these conditions immediately imply some relations. The fact that $V_0$ and any vector space with two endpoints are one-dimensional implies the following two sets of relations
\begin{center}
\begin{tabular}{cccccccccc}
$\begin{tikzpicture}[baseline={([yshift=-.5ex]current bounding box.center)}]
\draw[->] (0,0) arc (90:-270:.5);
\end{tikzpicture}
~=~ d_1$
&&
$\begin{tikzpicture}[baseline={([yshift=-.5ex]current bounding box.center)}]
\draw[<-] (0,0) arc (90:-270:.5);
\end{tikzpicture}
~=~ d_2$ &&&
$\begin{tikzpicture}[baseline=.5cm, rounded corners = 5mm]
\draw (0,1)--(1,0) -- (0,0)--(.4,.4);
\draw[->] (.6,.6)--(1,1);
\end{tikzpicture}
~=~
a_1 \cdot 
\begin{tikzpicture}[baseline = .5 cm, rounded corners = 6mm]
    \draw[->] (0,1)--(.5,-.1) -- (1,1);
\end{tikzpicture}$
&&
$\begin{tikzpicture}[baseline=.5cm, rounded corners = 5mm]
\draw[<-] (0,1)--(1,0) -- (0,0)--(.4,.4);
\draw(.6,.6)--(1,1);
\end{tikzpicture}
~=~
a_2 \cdot 
\begin{tikzpicture}[baseline = .5 cm, rounded corners = 6mm]
    \draw[<-] (0,1)--(.5,-.1) -- (1,1);
\end{tikzpicture}$
\\
\multicolumn{3}{c}{Circle parameters} &&&
\multicolumn{3}{c}{Twist Relations}
\end{tabular}
\end{center}
For an arbitrary planar algebra, these parameters can all be different. Since we have assumed the quotient is spherical, however, this implies that $d_1 = d_2$ and $a_1 = a_2.$ Thus, we will use the letters $d$ and $a$ to represent them. $GL_t^a$ and $\Rep(S^1,a)$ from Examples \ref{def:s1} and \ref{ex:glt} are examples of such a quotient. 

\setcounter{subsection}{1}

\begin{theorem}
Let $\mathcal{V}$ be a spherical quotient of the planar algebra of oriented virtual tangles with $\dim V_0 = 1,$ the dimension of all non-empty vector spaces with 2 endpoints is 1, and the dimension of all non-empty vector spaces with 4 endpoints is 1 or 2. Then $\mathcal{V}$ is isomorphic to one of the following:

\begin{enumerate}[label = \roman{enumi}.]
\item Deligne's $GL_{t}^a$ where
\begin{equation*}
    \overcrossinga ~=~ a\cdot\virtuala
\end{equation*}
\item The non-degenerate quotient of $GL_{\pm 1}^a,$ where the virtual crossing is equal to plus or minus the identity.
\end{enumerate}
\end{theorem}

\begin{proof}
Let $\mathcal{V}$ be described as above and fix an orientation on the vertices. Without loss of generality assume that orientation is $(+,+,-,-).$ Now suppose that there exists a relation between the virtual crossing and the identity. Because vR2 is satisfied, this implies that
\begin{equation}\label{eq:ovirtual}
\virtuala ~=~ \pm\identitya
\end{equation}
By capping off on the left, we see that $d = \pm 1$ by vR1. Furthermore, by expanding mR3 we see that (\ref{eq:ovirtual}) implies that there is a relation of the form
\begin{equation*}
\overcrossinga ~=~ z\cdot\identitya
\end{equation*}
By capping off on the left, we see that $z = \pm a$ depending on whether $d = \pm 1.$ It is trivial to check that this formula satisfies the remaining Reidemeister relations. Since these are exactly the skein relations of the non-degenerate quotient of $GL_{\pm 1}^a,$ any spherical quotient of the planar algebra of oriented virtual tangles with a dependence between the virtual crossing and the identity must be isomorphic to one of these planar algebras.

Now, suppose that the identity and the virtual crossing are linearly independent. Thus, we must have the following relations:
\begin{equation}\label{o41}
\overcrossinga ~=~ x_1\cdot\identitya ~+~ z_1\cdot\virtuala
\end{equation}
\begin{equation}\label{o42}
\undercrossinga ~=~ x_2\cdot\identitya ~+~ z_2\cdot\virtuala
\end{equation}
By capping off (\ref{o41}) and (\ref{o42}) above on the left we obtain the equations $a = dx_1 + z_1$ and $a = dx_2 + z_2.$ Since $a\ne 0$, we can always reparameterize the $x_i$ and $z_i$ for any value of $d,$ so without loss of generality assume that $a = 1$. Thus we have $1 = dx_1 +z_1$ and $1 = dx_2 + z_2.$ By multiplying (\ref{o41}) and (\ref{o42}) together, we see that $x_1x_2+z_1z_2 = 1$ and $x_1z_2 + x_2z_1 = 0.$ 

If we expand R3 with the above formula, we see that $x_1 = 0$ or $z_1 =0$ by the independence of the identity and virtual crossing. Solving this system of equations gives us 2 possible cases: The case where $x_1 = x_2 = 0$ and $z_1 = z_2 = 1$ or the case where $d = \pm 1.$ By mR3, the latter case would imply a relation between the virtual crossing and the identity, a contradiction to their independence. The former case is exactly the skein theory of $GL_{t}.$ Remembering that we made a choice of parameterization for $a$, we see that when the identity and virtual crossing are independent we obtain $GL_{t}^a,$ where $a$ represents the case where the crossing is equal to $a$ times the virtual crossing. This completes our proof. 
\end{proof}

\section{Classification of symmetric trivalent categories} \label{tri}
In \cite{MPS15}, a classification of trivalent planar algebras, $\mathcal{V},$ with $\dim V_4 \leq 4$ was given. A full description of that classification is given in Section \ref{trivalent}. In this section, we extend this result to symmetric trivalent planar algebras (See Definition \ref{def:vtri}) with $\dim V_4 \leq 4$. Some of the relations that will explicitly appear later are listed below:
\begin{center}
\begin{tabular}{ccccccccccc}
   $\skeincircle = d = t - 1$  
   &&& 
   $\begin{tikzpicture}[baseline]
 	\draw (0,-.6) -- (0,-.2);
    \draw (0,-.2) to[out=45,in=-90] (.35,.2) to[out=90,in=0] (0,.6)
    			  to[out=180,in=90] (-.35,.2) to[out=-90,in=135] (0,-.2);
\end{tikzpicture}
 ~=~ 0$ 
   &&&
   $\begin{tikzpicture}[baseline]
 	\draw (0,.75) -- (0,.45);
	\draw (0,-.75) -- (0,-.45);
    \draw (0,-.45) .. controls (-.25,0) .. (0,.45);
    \draw (0,-.45) .. controls (.25,0) .. (0,.45);
\end{tikzpicture}
~=~
c_1 \cdot ~
\begin{tikzpicture}[baseline]
 	\draw (0,.75) -- (0,-.75);
\end{tikzpicture}
$
   &&& 
   \includegraphics[valign = c,scale = 1.25]{diagrams/pdf/3boxtrirel.pdf}
   \\
   circle parameter &&& lollipop relation &&& bigon relation &&& triangle relation\\
\end{tabular}
\end{center}
We note that there is a choice of normalization that allows us to rescale $c_2$ appropriately to make $c_1$ any non-zero constant. Thus, we will set $c_1 = 1$.

In order to use the classification of \cite{MPS15}, we will also require that our planar algebra be nondegenerate. There are some degenerate symmetric trivalent planar algebras that fall within our dimension bounds. For example, $S_t$ is degenerate by Theorem \ref{thm:dim} at every $t\in\mathbb{Z}_{\geq 0}.$ Instead, we will consider its nondegenerate quotient at those values of $t$, which by Theorem \ref{thm:dim} is $\Rep(S_t).$  

\subsection{Classification of symmetric trivalent planar algebras}

Since we have imposed $\dim V_4 \leq 4$, it is necessarily true that we have a relation of the form
\begin{equation}\label{trirel}
a_1\stI ~+~ a_2\cdot\stH ~+~ a_3\cdot\begin{tikzpicture}[baseline, rounded corners = 5mm]
    \draw (.5,.5) -- (.1,0) -- (.5,-.5);
    \draw (-.5,.5) -- (-.1,0) -- (-.5,-.5);
\end{tikzpicture} ~+~ a_4\cdot\begin{tikzpicture}[baseline, rounded corners = 5mm]
    \draw (-.5,-.5) -- (0,.1) -- (.5,-.5);
    \draw (-.5,.5) -- (0,-.1) -- (.5,.5);
\end{tikzpicture} ~+~ a_5\cdot\begin{tikzpicture}[baseline]
    \draw (-.5,-.5)-- (.5,.5);
    \draw (.5,-.5) -- (-.5,.5);
    \draw[fill=white] (0,0) circle (.1cm);
\end{tikzpicture}
 ~=~ 0
\end{equation}
with at least two of the $a_i\in\mathbb{C}$ not 0. We have already seen many examples in Section \ref{trivalent} of symmetric trivalent planar algebras. A priori, the above relation could give many different skein relations. As the next two lemmas show, however, we can easily classify all skein theories with $\dim V_4 = 1$ or $2.$

\begin{lemma}\label{tri1dim}
In any symmetric trivalent planar algebra, the cupcap and identity must be linearly independent. In particular, any symmetric trivalent planar algebra must have $\dim V_4 \geq 2$.
\end{lemma}

\begin{proof}
If the identity and cupcap were linearly dependent then, 
\begin{equation*}
\begin{tikzpicture}[baseline = -.05cm, rounded corners = 5mm]
    \draw (0,-.6) -- (0,0)
    	  (0,0) -- (-.5,.5)
    	  (0,0) -- (.5,.5);
\end{tikzpicture} ~=~ \left(~\begin{tikzpicture}[baseline = -.05cm, rounded corners = 5mm]
    \draw (0,-.6) -- (0,0)
    	  (0,0) -- (-.5,.5)
    	  (0,0) -- (.5,.5);
\end{tikzpicture}~\right) \cdot \left(~\begin{tikzpicture}[baseline, rounded corners = 5mm]
    \draw (.5,.5) -- (.1,0) -- (.5,-.5);
    \draw (-.5,.5) -- (-.1,0) -- (-.5,-.5);
\end{tikzpicture}~\right)  ~=~ 
\left(~\begin{tikzpicture}[baseline = -.05cm, rounded corners = 5mm]
    \draw (0,-.6) -- (0,0)
    	  (0,0) -- (-.5,.5)
    	  (0,0) -- (.5,.5);
\end{tikzpicture}~\right) \cdot \left(~x\cdot\begin{tikzpicture}[baseline, rounded corners = 5mm]
    \draw (-.5,-.5) -- (0,.1) -- (.5,-.5);
    \draw (-.5,.5) -- (0,-.1) -- (.5,.5);
\end{tikzpicture}~\right) ~=~ 0
\end{equation*}
which contradicts our assumption that $\dim V_3 = 1$. Thus, the cupcap and identity must be linearly independent, which implies that $\dim V_4 \geq 2$ for any symmetric trivalent planar algebra.
\end{proof}
\begin{lemma}\label{tri2dim}
There are no symmetric trivalent planar algebras with $\dim V_4 = 2$. In addition, the cupcap, identity, and $I$ must be linearly independent.
\end{lemma}

\begin{proof}
Assume that we have a symmetric trivalent planar algebra with $\dim V_4 = 2$. By Lemma \ref{tri1dim}, we know that the cupcap and the identity must be linearly independent and thus span $V_4$. Hence we must a relation of the form:
\begin{equation*}
\stI ~ = ~ A\begin{tikzpicture}[baseline, rounded corners = 5mm]
    \draw (.5,.5) -- (.1,0) -- (.5,-.5);
    \draw (-.5,.5) -- (-.1,0) -- (-.5,-.5);
\end{tikzpicture} ~+~ B\begin{tikzpicture}[baseline, rounded corners = 5mm]
    \draw (-.5,-.5) -- (0,.1) -- (.5,-.5);
    \draw (-.5,.5) -- (0,-.1) -- (.5,.5);
\end{tikzpicture} 
\end{equation*}
The solution to this equation is well-known and yields on two possible relations:
\begin{equation*}
\stI ~ = ~ \begin{tikzpicture}[baseline, rounded corners = 5mm]
    \draw (.5,.5) -- (.1,0) -- (.5,-.5);
    \draw (-.5,.5) -- (-.1,0) -- (-.5,-.5);
\end{tikzpicture} ~-~ \dfrac{1\pm\sqrt{5}}{2}~ \begin{tikzpicture}[baseline, rounded corners = 5mm]
    \draw (-.5,-.5) -- (0,.1) -- (.5,-.5);
    \draw (-.5,.5) -- (0,-.1) -- (.5,.5);
\end{tikzpicture} 
\end{equation*}
Multiplying both sides of the above equation by the virtual crossing implies that the virtual crossing is equal to the identity. The virtual crossing is rotationally invariant, however, and the identity is not. This implies a relation between the identity and the cupcap, a contradiction to their independence. Thus, no such skein theories exist. Moreover, given that any dependence amongst the cupcap, identity, and I would yield the same contradiction, they must be linearly independent.
\end{proof}

Thus, there are no symmetric trivalent planar algebras with $\dim V_4 = 1$ or $2.$ This leaves us only with the cases where the $\dim V_4 = 3$ or $4.$

\begin{lemma}\label{thm:ih}
There is no symmetric trivalent planar algebra with a relation of the form:
\begin{equation*}
\stI ~=~ z\cdot\stH
\end{equation*}
\end{lemma}
\begin{proof}
If there were such a relation, capping the top of both sides would imply by nondegeneracy that either $\dim V_3 = 0$ (if $z=0$) or that $\dim V_2 = 0$, both of which are contradictions to our dimension constraints. Thus, there can be no such dependence.
\end{proof}
We will also make use of the following fact from \cite{MPS15}:
\begin{proposition}[Corollary 8.9 in \cite{MPS15}]\label{thm:sym}
The only trivalent planar algebras with a symmetric braiding are $\Rep(SO(3))$, $\Rep(S_3)$, $\Rep(G_2)$, and $\Rep(OSp(1|2))$ and these are the only braidings on those planar algebras.
\end{proposition}

We also have the following fact:

\begin{lemma}\label{thm:St}
Suppose we have a symmetric trivalent planar algebra with $t \ne 0$ or $2.$ Then a planar algebra has a relation of the form
\begin{equation*}
\stI + \stH = z\cdot\left[~\begin{tikzpicture}[baseline, rounded corners = 5mm]
    \draw (.5,.5) -- (.1,0) -- (.5,-.5);
    \draw (-.5,.5) -- (-.1,0) -- (-.5,-.5);
\end{tikzpicture} + \begin{tikzpicture}[baseline, rounded corners = 5mm]
    \draw (-.5,-.5) -- (0,.1) -- (.5,-.5);
    \draw (-.5,.5) -- (0,-.1) -- (.5,.5);
\end{tikzpicture}~\right]
\end{equation*}
if and only if it has a relation of the form
\begin{equation*}
\stI - \stH = z'\cdot\left[~\begin{tikzpicture}[baseline, rounded corners = 5mm]
    \draw (.5,.5) -- (.1,0) -- (.5,-.5);
    \draw (-.5,.5) -- (-.1,0) -- (-.5,-.5);
\end{tikzpicture} - \begin{tikzpicture}[baseline, rounded corners = 5mm]
    \draw (-.5,-.5) -- (0,.1) -- (.5,-.5);
    \draw (-.5,.5) -- (0,-.1) -- (.5,.5);
\end{tikzpicture}~\right]
\end{equation*}
with $z = \dfrac{1}{t}$ and $z' = \dfrac{1}{t-2}.$
\end{lemma}

\begin{proof}
Assume that the first relations holds. By capping off on the top, we see that $z = \dfrac{1}{t}.$ By using the above I=H relation to simplify the triangle, we see that $c_2 = -\dfrac{t-1}{t}.$ By multiplying the first relation by an H we get the following relation:
\begin{equation*}
-\dfrac{t-1}{t}\cdot\stI ~+~ \begin{tikzpicture}[baseline=0cm]
 	\draw (0,.3) -- (-.3,.6);
    \draw (0,-.3) -- (-.3,-.6);
    \draw (.6,.3) -- (.9,.6);
    \draw (.6,-.3) -- (.9,.-.6);
	\draw (0,-.3) rectangle (.6, .3);
\end{tikzpicture} ~=~ \dfrac{1}{t}\cdot\stH ~+~ \dfrac{1}{t}\cdot\begin{tikzpicture}[baseline, rounded corners = 5mm]
    \draw (-.5,-.5) -- (0,.1) -- (.5,-.5);
    \draw (-.5,.5) -- (0,-.1) -- (.5,.5);
\end{tikzpicture}
\end{equation*}
Since the square diagram above is rotationally invariant, by taking the difference of this equation with its rotation we obtain the following relation: 
\begin{equation*}
\stI - \stH = \dfrac{1}{t-2}\cdot\left[~\begin{tikzpicture}[baseline, rounded corners = 5mm]
    \draw (.5,.5) -- (.1,0) -- (.5,-.5);
    \draw (-.5,.5) -- (-.1,0) -- (-.5,-.5);
\end{tikzpicture} - \begin{tikzpicture}[baseline, rounded corners = 5mm]
    \draw (-.5,-.5) -- (0,.1) -- (.5,-.5);
    \draw (-.5,.5) -- (0,-.1) -- (.5,.5);
\end{tikzpicture}~\right],
\end{equation*}
which exists when $t \ne 2$. If we have a relation of the form
\begin{equation*}
\stI - \stH = z'\cdot\left[~\begin{tikzpicture}[baseline, rounded corners = 5mm]
    \draw (.5,.5) -- (.1,0) -- (.5,-.5);
    \draw (-.5,.5) -- (-.1,0) -- (-.5,-.5);
\end{tikzpicture} - \begin{tikzpicture}[baseline, rounded corners = 5mm]
    \draw (-.5,-.5) -- (0,.1) -- (.5,-.5);
    \draw (-.5,.5) -- (0,-.1) -- (.5,.5);
\end{tikzpicture}~\right],
\end{equation*}
We see by capping that $z' = \dfrac{1}{t-2}.$ Again using this I=H relation to simplify the triangle, we see that $c_2 = \dfrac{t-1}{t-2}.$ Multiplying the I=H relation by an H on top tells us that
\begin{equation*}
\dfrac{t-1}{t-2}\cdot\stI ~-~ \begin{tikzpicture}[baseline=0cm]
 	\draw (0,.3) -- (-.3,.6);
    \draw (0,-.3) -- (-.3,-.6);
    \draw (.6,.3) -- (.9,.6);
    \draw (.6,-.3) -- (.9,.-.6);
	\draw (0,-.3) rectangle (.6, .3);
\end{tikzpicture} ~=~ \dfrac{1}{t-2}\cdot\stH ~-~ \dfrac{1}{t-2}\cdot\begin{tikzpicture}[baseline, rounded corners = 5mm]
    \draw (-.5,-.5) -- (0,.1) -- (.5,-.5);
    \draw (-.5,.5) -- (0,-.1) -- (.5,.5);
\end{tikzpicture}
\end{equation*}
If we again take the difference of this relation with its rotation we get the new relation
\begin{equation*}
\stI + \stH = \dfrac{1}{t}\cdot\left[~\begin{tikzpicture}[baseline, rounded corners = 5mm]
    \draw (.5,.5) -- (.1,0) -- (.5,-.5);
    \draw (-.5,.5) -- (-.1,0) -- (-.5,-.5);
\end{tikzpicture} + \begin{tikzpicture}[baseline, rounded corners = 5mm]
    \draw (-.5,-.5) -- (0,.1) -- (.5,-.5);
    \draw (-.5,.5) -- (0,-.1) -- (.5,.5);
\end{tikzpicture}~\right]
\end{equation*}
when $t \ne 0$, which completes our proof. 
\end{proof}

We are now ready to classify all symmetric trivalent planar algebras:

\begin{theorem}\label{thm:class}
Let $\mathcal{V}$ be a symmetric trivalent planar algebra. Then $\mathcal{V}$ is one of the following:
\begin{center}
\begin{tabular}{l | l}
Dimensions & Name \\
\hline 
$1,0,1,1,3,\dots$ & $\Rep(SO(3)),$ $\Rep(OSp(1|2)),$ and $\Rep(S_3)$  \\
$1,0,1,1,4,9,\dots$ & $\Rep(G_2)$ \\
$1,0,1,1,4,10,\dots$ & $\Rep(S_4)$ \\
$1,0,1,1,4,11,\dots$ & $S_t$ for $t\not\in\mathbb{Z}_{\geq 0}$ \\
$1,0,1,1,4,11,\dots$ & $\Rep(S_t)$ for $t\in\mathbb{Z}_{\geq 0}-\{0,1,2,3\}$ \\
\end{tabular}
\end{center}
\end{theorem}

\begin{proof}
By Lemmas \ref{tri1dim} and \ref{tri2dim}, we know that $\dim V_4 \geq 3$ for any symmetric trivalent planar algebra. If the virtual crossing can be written in terms of the other relations, Proposition \ref{thm:sym} gives us the classification of such planar algebras. Thus, we will assume now that the virtual crossing cannot be written in terms of the other generators. This assumption together with Lemmas \ref{tri2dim} and \ref{thm:ih} tell us that any such symmetric trivalent planar algebra must have $\dim V_4 = 4$. This implies that we have a relation of the following form:
\begin{equation*}
\stI ~=~ x~\stH ~+~ y ~\begin{tikzpicture}[baseline, rounded corners = 5mm]
    \draw (.5,.5) -- (.1,0) -- (.5,-.5);
    \draw (-.5,.5) -- (-.1,0) -- (-.5,-.5);
\end{tikzpicture} ~+~ z~\begin{tikzpicture}[baseline, rounded corners = 5mm]
    \draw (-.5,-.5) -- (0,.1) -- (.5,-.5);
    \draw (-.5,.5) -- (0,-.1) -- (.5,.5);
\end{tikzpicture}
\end{equation*}
By rotating the above equation and solving for I again, we obtain two possible relations:
\begin{equation*}
\stI - \stH = \dfrac{1}{t-2}\left[\begin{tikzpicture}[baseline, rounded corners = 5mm]
    \draw (.5,.5) -- (.1,0) -- (.5,-.5);
    \draw (-.5,.5) -- (-.1,0) -- (-.5,-.5);
\end{tikzpicture} - \begin{tikzpicture}[baseline, rounded corners = 5mm]
    \draw (-.5,-.5) -- (0,.1) -- (.5,-.5);
    \draw (-.5,.5) -- (0,-.1) -- (.5,.5);
\end{tikzpicture}\right]
\end{equation*}
or
\begin{equation*}
\stI + \stH = \dfrac{1}{t}\left[\begin{tikzpicture}[baseline, rounded corners = 5mm]
    \draw (.5,.5) -- (.1,0) -- (.5,-.5);
    \draw (-.5,.5) -- (-.1,0) -- (-.5,-.5);
\end{tikzpicture} + \begin{tikzpicture}[baseline, rounded corners = 5mm]
    \draw (-.5,-.5) -- (0,.1) -- (.5,-.5);
    \draw (-.5,.5) -- (0,-.1) -- (.5,.5);
\end{tikzpicture}\right]
\end{equation*}
By Lemma \ref{thm:St}, however, we know that any skein theory has one relation if and only if it has the other when $t \ne 0$ or $2.$ At $t=0,$ note that the second relation is undefined. At $t = 2$, the first relation is undefined and the second relation implies that I and H are equal, which is impossible by Lemma \ref{thm:ih}. Thus $t \ne 2$ and so we will be justified in assuming the top relation holds.
 
The planar algebra with this I=H relation is exactly the skein theory of Deligne's $S_t,$ described in Example \ref{ex:St}. However, at $t\in\mathbb{Z}_{\geq 0}$ the planar algebra is degenerate. Taking a non-degenerate quotient at these values gives $\Rep(S_t)$. When $t=0$ this planar algebra is isomorphic to $OSp(1|2)$. Since we have now exhausted all possible skein relations, this must be a complete list.
\end{proof}

\subsection{Sub-braidings of symmetric trivalent planar algebras}

We wish to classify all sub-braidings (See Definition \ref{subbraid}) on these planar algebras. In particular, we want to find all relations of the form
\begin{equation}\label{eq:crossing}
\begin{tikzpicture}[baseline]
\begin{knot}[clip width = 4]
	\strand (-.5,.5) -- (.5,-.5);
	\strand (-.5,-.5) -- (.5,.5);
\end{knot}
\end{tikzpicture}
 ~=~ x\begin{tikzpicture}[baseline, rounded corners = 5mm]
    \draw (.5,.5) -- (.1,0) -- (.5,-.5);
    \draw (-.5,.5) -- (-.1,0) -- (-.5,-.5);
\end{tikzpicture} ~+~ y\begin{tikzpicture}[baseline, rounded corners = 5mm]
    \draw (-.5,-.5) -- (0,.1) -- (.5,-.5);
    \draw (-.5,.5) -- (0,-.1) -- (.5,.5);
\end{tikzpicture} ~+~ z\begin{tikzpicture}[baseline]
    \draw (-.5,-.5)-- (.5,.5);
    \draw (.5,-.5) -- (-.5,.5);
    \draw[fill=white] (0,0) circle (.1cm);
\end{tikzpicture}
 ~+~ w\stI
\end{equation}
When looking for sub-braidings, however, we again note that while there may be an element of the 4-box space which satisfies all of the Reidemeister moves, it may not meet the naturality condition described in Definition \ref{def:naturality}. We could impose those requirements outright. Namely,
\begin{equation}\label{R25}
\includegraphics[valign = c,scale = 1.25]{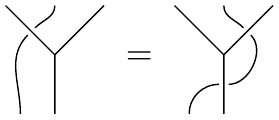}
\end{equation}
We will call this a Type 2.5 Reidemeister move. If a planar algebra satisfies R2.5, then we will say it is fully sub-braided. If we instead asked that
\begin{equation}
\includegraphics[valign = c,scale = .75]{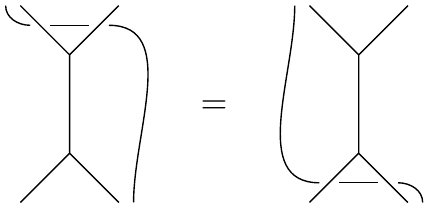}
\end{equation}
be satisfied but not necessarily R2.5, we call this an even sub-braiding. Note that every full sub-braiding gives rise to an even sub-braiding as evidenced by the following lemma: 

\begin{lemma}\label{half}
Every full sub-braiding of a symmetric trivalent planar algebra is equivalent to an even sub-braiding.
\end{lemma}

\begin{proof}
Suppose that we have a full sub-braiding denoted by \begin{tikzpicture}[baseline = -.15cm, scale = .65]
\begin{knot}[clip width = 4]
	\strand (-.5,.5) -- (.5,-.5);
	\strand (-.5,-.5) -- (.5,.5);
\end{knot}
\end{tikzpicture}
. Then it satisfies (\ref{R25}). We create a new element of the planar algebra called \begin{tikzpicture}[baseline = -.15cm, scale = .65]
\begin{knot}[clip width = 4]
	\strand[dashed] (-.5,.5) -- (.5,-.5);
	\strand[dashed] (-.5,-.5) -- (.5,.5);
\end{knot}
\end{tikzpicture} with $\begin{tikzpicture}[baseline = -.15cm, scale = .65]
\begin{knot}[clip width = 4]
	\strand[dashed] (-.5,.5) -- (.5,-.5);
	\strand[dashed] (-.5,-.5) -- (.5,.5);
\end{knot}
\end{tikzpicture} ~=~ -\begin{tikzpicture}[baseline = -.15cm, scale = .65]
\begin{knot}[clip width = 4]
	\strand (-.5,.5) -- (.5,-.5);
	\strand (-.5,-.5) -- (.5,.5);
\end{knot}
\end{tikzpicture}
.$ Clearly, \begin{tikzpicture}[baseline = -.15cm, scale = .65]
\begin{knot}[clip width = 4]
	\strand[dashed] (-.5,.5) -- (.5,-.5);
	\strand[dashed] (-.5,-.5) -- (.5,.5);
\end{knot}
\end{tikzpicture} satisfies the virtual Reidemeister moves because we assumed that \begin{tikzpicture}[baseline = -.15cm, scale = .65]
\begin{knot}[clip width = 4]
	\strand (-.5,.5) -- (.5,-.5);
	\strand (-.5,-.5) -- (.5,.5);
\end{knot}
\end{tikzpicture}
 did. By substituting \begin{tikzpicture}[baseline = -.15cm, scale = .65]
\begin{knot}[clip width = 4]
	\strand[dashed] (-.5,.5) -- (.5,-.5);
	\strand[dashed] (-.5,-.5) -- (.5,.5);
\end{knot}
\end{tikzpicture} for \begin{tikzpicture}[baseline = -.15cm, scale = .65]
\begin{knot}[clip width = 4]
	\strand (-.5,.5) -- (.5,-.5);
	\strand (-.5,-.5) -- (.5,.5);
\end{knot}
\end{tikzpicture}
 in (\ref{R25}), we see that the left-hand side is negative, while the right-hand side remains positive. Hence, \begin{tikzpicture}[baseline = -.15cm, scale = .65]
\begin{knot}[clip width = 4]
	\strand[dashed] (-.5,.5) -- (.5,-.5);
	\strand[dashed] (-.5,-.5) -- (.5,.5);
\end{knot}
\end{tikzpicture} is an even sub-braiding and every full sub-braiding corresponds to an even sub-braiding, as desired. 
\end{proof}

A priori, there could be even sub-braidings that have no corresponding full sub-braiding. By inspection, however, this is not the case. To begin to classify the full sub-braidings of our planar algebras, we note that Theorem \ref{thm:sym} tells us the only full sub-braidings (which are also actual braidings) of $\Rep(SO(3)),$ $\Rep(OSp(1|2)),$ $\Rep(G_2),$ and $\Rep(S_3)$ are fully-flat. The only remaining sub-braidings to classify then are the ones for $S_t.$

\begin{proposition}\label{stbraid}
The only non-fully-flat sub-braidings of $S_t$ for generic $t$ and $\Rep(S_t)$ for $t\in\mathbb{Z}_{\geq 0}-\{0,1,2,3\}$ is the sub-braiding inherited from $SO(3)_{q}$:
\begin{equation*}
\begin{tikzpicture}[baseline]
\begin{knot}[clip width = 4]
	\strand (-.5,.5) -- (.5,-.5);
	\strand (-.5,-.5) -- (.5,.5);
\end{knot}
\end{tikzpicture}
 ~=~ (q^2-1)\cdot\begin{tikzpicture}[baseline, rounded corners = 5mm]
    \draw (.5,.5) -- (.1,0) -- (.5,-.5);
    \draw (-.5,.5) -- (-.1,0) -- (-.5,-.5);
\end{tikzpicture} ~+~ q^{-2}\begin{tikzpicture}[baseline, rounded corners = 5mm]
    \draw (-.5,-.5) -- (0,.1) -- (.5,-.5);
    \draw (-.5,.5) -- (0,-.1) -- (.5,.5);
\end{tikzpicture} ~-~(q^2+q^{-2})\stI
\end{equation*}
for $t = q^2 + 2 + q^{-2}.$ At $t=0$ and $3$ all sub-braidings are fully-flat.
\end{proposition}

\begin{proof}
At $t=1$ and $2$, $S_t$ is not a symmetric trivalent planar algebra, so they are excluded from our classification. Because of the laborious nature of determining all possible sub-braidings on the other values of $t$, a Mathematica program was written to find all such sub-braidings. The Mathematica code used for this is included in the materials uploaded to the Arxiv but follows this algorithm: 
\begin{enumerate}[label = \roman{enumi}.]
\item  We begin by expanding both sides of the relations implied by the Reidemeister moves and R2.5 using (\ref{eq:crossing}). For more detailed examples of this type of calculation see Section \ref{virtualskein}.
\item After removing the crossings, we can now think of these diagrams as graphs with only trivalent and univalent vertices. Because of the naturality relations we assumed between the virtual crossing and trivalent crossing, we can think of the virtual crossing as a crossing of edges on our graph.
\item Next, we can use the I=H relation to simplify complicated diagrams into our chosen basis according to the following process:
\begin{enumerate}
\item Remove all faces from any diagram.
\item Turn any tree (in the graph-theoretic sense) of $n$ univalent vertices into a specified tree of the same number of univalent vertices.
\end{enumerate}
\item Now that all of our diagrams are in terms of our chosen basis, we can solve the equations given by the Reidemeister moves (including R2.5) for generic values of $t$.
\item At $t=0,$ $3,$ $4,$ or $5$, there is a dimension drop in the $6$-box space or lower, and so one must calculate the null space of the inner product matrix of the appropriate $n$-box spaces to transform our chosen generic spanning set into the specific basis at those values.  
\end{enumerate}

The results of the calculation are the ones given above. Except when $t=0$ or $3$, the only non-trivial sub-braidings were the $(SO(3))_q$ sub-braidings. At $t = 0$ or $3,$ $\dim V_4 = 3$, so there is a relation giving I as a linear combination of the identity, cupcap, and virtual crossing. When one makes this substitution, the other terms cancel and we see that the crossing is equal to the symmetric braiding, giving a fully-flat braiding. 
\end{proof}

The work from \cite{MPS15} and the above proposition give us the following theorem:

\begin{corollary}
The only symmetric trivalent planar algebras exhibiting a non-fully-flat sub-braiding are $S_t$ for generic $t$ and $\Rep(S_t)$ with $t\in\mathbb{Z}_{\geq 3}.$
\end{corollary}

\begin{proof}
 By Theorem \ref{thm:sym}, the only planar algebra that could possibly exhibit a non-trivial sub-braiding is $S_t.$ The classification of such sub-braidings is proven in Proposition \ref{stbraid}, where we see the only such sub-braiding is the one inherited from $(SO(3))_q.$
\end{proof}

By the note in Example \ref{ex:St}, when $S_t$ is equipped with this braiding it is also called the second-colored Jones polynomial planar algebra.

\begin{corollary}
Every even sub-braiding of a symmetric trivalent planar algebra corresponds to a full sub-braiding.
\end{corollary}

\begin{proof}
When all possible even sub-braidings were classified in Mathematica, the formula for the crossing for every even sub-braiding could be obtained from a full sub-braiding via the method described in Lemma \ref{half}. Thus, by exhaustion no other even sub-braidings exist.  
\end{proof}

\bibliographystyle{alpha}
\bibliography{text/DissBib}

\end{document}